\UseRawInputEncoding
\documentclass[12pt]{article}
\usepackage{mathrsfs}
\usepackage{amssymb}
\usepackage{amsmath}
\usepackage{amsfonts}
\usepackage{amstext}
\usepackage{amsthm}
\usepackage{graphicx}

\usepackage{pict2e}
\newcommand{\PPWaxislinewidth}{0.96}
\newcommand{\PPWarrowlength}{6.5}
\newcommand{\PPWarrowhalfwidth}{2.5}
\newcommand{\PPWarrowinset}{5.2}
\newcommand{\PPWxaxis}[3]{%
  \begingroup\linethickness{\PPWaxislinewidth bp}%
  \Line(#1,#2)(#3,#2)%
  \put(#3,#2){\polygon*(0,0)(-\PPWarrowlength,\PPWarrowhalfwidth)%
    (-\PPWarrowinset,0)(-\PPWarrowlength,-\PPWarrowhalfwidth)}%
  \endgroup}
\newcommand{\PPWyaxis}[3]{%
  \begingroup\linethickness{\PPWaxislinewidth bp}%
  \Line(#1,#2)(#1,#3)%
  \put(#1,#3){\polygon*(0,0)(-\PPWarrowhalfwidth,-\PPWarrowlength)%
    (0,-\PPWarrowinset)(\PPWarrowhalfwidth,-\PPWarrowlength)}%
  \endgroup}
\usepackage{subfig}
\usepackage{placeins}

\usepackage{indentfirst,latexsym,bm}
\usepackage{mathrsfs}
\usepackage{cite}
\usepackage[all]{xy}
\usepackage{newlfont}

\usepackage{enumerate}
\numberwithin{equation}{section}

\newtheorem{theorem}{Theorem}[section]

\newtheorem{proposition}[theorem]{Proposition}
\numberwithin{equation}{section}
\allowdisplaybreaks[1]

\begin{document}

\title{{Upper bounds for the ratio of the first two eigenvalues of the Robin Laplacian}
\thanks{Research supported by NNSF of China (No. 12371110).}}

\author{Guowei Dai\thanks{Corresponding author.
\newline
School of Mathematical Sciences, Dalian University of Technology, Dalian, 116024, P.R. China
\newline
\text{~~~~ E-mail}: daiguowei@dlut.edu.cn.}, Yingxin Sun}
\date{}
\maketitle

\renewcommand{\abstractname}{Abstract}

\begin{abstract}
We sharpen upper bounds for the ratio of the first two Robin eigenvalues
by matching the first eigenvalue on a comparison ball with an independently
chosen boundary parameter. The bound depends on the interior superlevel
sets of the first eigenfunction and improves the fixed-parameter
comparison. For the conditional ball inequality, we also
prove that equality holds only for a ball. On balls in every dimension
$N\geq2$, we establish strict increase of the first spectral gap and
strict decrease of the ratio of the first two eigenvalues with respect to the positive
Robin parameter, verifying the ball case of Laugesen's spectral-ratio
monotonicity conjecture. Both statements follow from a comparison of the
lowest radial eigenvalues of any two spherical harmonic degrees.
Finally, building on the known disconnected examples, we give an explicit
construction showing that the ratio is unbounded among smooth connected
domains of prescribed volume at every fixed positive parameter; in two
dimensions the domains can be chosen simply connected. The equality
case yields a conditional circularity criterion for elastically
supported membranes.
\end{abstract}




\section{Introduction}\label{sec:introduction}
\bigskip
\quad\, Let $\lambda_n$ with $n\in \mathbb{N}$ denote the eigenvalues (counting multiplicities and arranged in nondecreasing order) of the Dirichlet eigenvalue problem
\begin{equation*}
\left\{
\begin{array}{ll}
-\Delta u=\lambda u\,\,\, &\text{in}\,\,\, \Omega,\\
u=0 &\text{on}\,\,\, \partial \Omega,
\end{array}
\right.
\end{equation*}
where $\Omega\subseteq \mathbb{R}^N$ with $N\geq2$ is a bounded domain (connected open set).

Payne, P\'{o}lya and Weinberger \cite{PPW} proved in 1956 that
$\lambda_2/\lambda_1\leq3$ in two dimensions. Their inequalities for
consecutive eigenvalues were subsequently strengthened by Hile and
Protter \cite{HileProtter} and Yang \cite{Yang}; for the first two
planar eigenvalues, these estimates give the same constant $3$.
Payne, P\'{o}lya and Weinberger conjectured the sharp bound
\begin{equation*}
\frac{\lambda_2}{\lambda_1}\leq \frac{\lambda_2}{\lambda_1}\Big|_{\text{disk}}=\frac{j_{1,1}^2}{j_{0,1}^2}\approx2.5387,
\end{equation*}
where $j_{\nu,k}$ denotes the
$k$th positive zero of the Bessel function
$J_\nu(x)$; see \cite{Abramowitz}.
Yau \cite[Problem~77]{Yau} noted that this bound would characterize
circular drums by their first two tones.

Brands \cite{Brands}, de Vries \cite{Vries}, and Chiti \cite{Chiti}
improved the upper bound to $2.686$, $2.658$, and $2.586$, respectively.
The higher-dimensional PPW conjecture is
\begin{equation*}
\frac{\lambda_2}{\lambda_1}\leq \frac{\lambda_2}{\lambda_1}\Big|_{N-\text{dimensional ball}}=\frac{j_{N/2,1}^2}{j_{N/2-1,1}^2}.
\end{equation*}
Ashbaugh and Benguria \cite{AshbaughBenguria91,AshbaughBenguria92,AshbaughBenguria92CMP}
proved this conjecture in the early 1990s.

The Robin eigenvalue problem is
\begin{equation*}
\left\{
\begin{array}{ll}
-\Delta u=\mu u\,\,\, &\text{in}\,\,\, \Omega,\\
\partial_\texttt{n} u+\sigma u=0 &\text{on}\,\,\, \partial \Omega.
\end{array}
\right.
\end{equation*}
\begingroup
Here $\sigma>0$ and $\mathbf n$ is the outward unit normal.
In Sections~\ref{sec:weights}--\ref{sec:ball-branches}, the boundary is
$C^2$, ensuring $u_1\in C^1(\overline\Omega)$ and the applicability of
the pointwise Robin condition and Hopf lemma \cite{GT,Bucur}.
Section~\ref{sec:necessity} also uses variational arguments on Lipschitz
domains. The Robin eigenvalues, counted with multiplicity, are denoted
by $\mu_j(\Omega,\sigma)$, or $\mu_j$ when both arguments are fixed.

For a conducting body, the outward heat flux is proportional to the
temperature relative to the surroundings. After normalizing the
conductivity, $\sigma$ is the heat-transfer coefficient. The eigenvalues
describe the decay rates of separated temperature modes. For an
elastically supported membrane, the restoring force at the edge is
proportional to its displacement. In this model the vibration
frequencies are proportional to $\sqrt{\mu_j}$;
see \cite{Courant,Bucur,Nahon2022}.

The general Robin theory and shape optimization are surveyed in
\cite{Bucur,Henrot}. Related results include estimates on
balls and planar domains \cite{Freitas,FreitasLaugesen,Freitas2021},
optimization on rectangles and their unions \cite{FreitasKennedy,Laugesen},
and the second eigenvalue on curved membranes \cite{Langford}.
The negative-parameter theory has a different behavior; see
\cite{Antunes,Freitas2015} and the formulas in
Appendix.

In 2001, Payne and Schaefer \cite{PayneSchaefer} asked whether the
disk maximizes $\mu_2/\mu_1$ for all positive boundary parameters,
or at least for an appropriate range. Henrot
\cite[Open problem~15]{Henrot0} subsequently restated the parameter
question. In higher dimensions the analogous comparison domain is a
ball. Laugesen \cite[Section~III.D, Conjecture~G]{Laugesen} formulates
the fixed-volume ball-maximization conjecture for convex domains and
records a counterexample of Bucur to the analogous
claim without convexity. The present results concern general connected
domains; the convex-domain conjecture remains separate. The volume
constraint is essential because
\begin{equation}\label{eq:robin-scaling}
\mu_j(\rho\Omega,\sigma)=\rho^{-2}\mu_j(\Omega,\rho\sigma),\,\,\, \rho>0.
\end{equation}
Let $B$ be the unit ball, let $C_N=|B|$, and put $\nu=N/2-1$. We normalize $|\Omega|=C_N$. For $\sigma>0$, denote by $k_{\nu+\ell,1}(\sigma)$ the first positive root of
\begin{equation}\label{eq:ball-root-definition}
xJ_{\nu+\ell+1}(x)-(\sigma+\ell)J_{\nu+\ell}(x)=0,
\,\,\, \ell=0,1,2,\ldots.
\end{equation}
The first two ball eigenvalues are $k_{\nu,1}(\sigma)^2$ and $k_{\nu+1,1}(\sigma)^2$; see Proposition~\ref{prop:ball-first-two} in Appendix. For a fixed $\sigma>0$, let
\[
\alpha=k_{\nu,1}(\sigma),\,\,\,
\beta=k_{\nu+1,1}(\sigma),\,\,\,
R=\frac{\alpha}{\sqrt{\mu_1(\Omega,\sigma)}}.
\]
Choose a positive first eigenfunction $u_1$, and define
\[
u_{1,pM}=\max_{\partial\Omega}u_1,\,\,\,
\mu(t)=|\{x\in\Omega:u_1(x)>t\}|,\,\,\,
\widetilde R=\left(\frac{\mu(u_{1,pM})}{C_N}\right)^{1/N}.
\]
The maximum principle and the Robin condition give $u_1>0$ on $\overline\Omega$ and $u_{1,pM}<\max_{\overline\Omega}u_1$, so $\widetilde R>0$.
\\ \\
\refstepcounter{thm}\label{thm:ball-bound}\textbf{Theorem~\thethm.}
\emph{Let $\Omega\subset\mathbb R^N$, $N\ge2$, be a bounded connected $C^2$ domain with $|\Omega|=|B|$. At any parameter $\sigma>0$ for which $\widetilde R\ge R$,}
\[
\frac{\mu_2(\Omega,\sigma)}{\mu_1(\Omega,\sigma)}
\le \frac{\mu_2(B,\sigma)}{\mu_1(B,\sigma)}
=\frac{\beta^2}{\alpha^2}.
\]
\emph{Equality holds if and only if $\Omega$ is a translate of $B$.}\\

The upper bound in Theorem~\ref{thm:ball-bound} is already stated in the
preprint of Dai and Shi \cite[Corollary~1.3]{DaiShi}. More generally,
\cite[Theorem~1.1]{DaiShi} gives, with $r_0=\min\{R,\widetilde R\}$,
\begin{equation}\label{eq:prior-bound}
\frac{\mu_2}{\mu_1}
\le 1+\frac{R^2}{r_0^2}
\left(\frac{\beta^2}{\alpha^2}-1\right).
\end{equation}

The radius $\widetilde R$ does not depend on the normalization of $u_1$.
Indeed, for every $c>0$,
\[
\{cu_1>\max_{\partial\Omega}(cu_1)\}
=\{u_1>\max_{\partial\Omega}u_1\}.
\]
The condition $\widetilde R\geq R$ ensures that the level sets used
in the comparison do not meet $\partial\Omega$, avoiding an exterior
boundary term.

For each fixed domain, $\widetilde R\geq R$ holds for all sufficiently
large $\sigma$. Proposition~\ref{prop:boundary-maximum} gives
$\widetilde R(\sigma)\to1$ and, for a nonball domain,
$R(\sigma)\to j_{\nu,1}/\sqrt{\lambda_1(\Omega)}<1$.
The large-parameter comparison also follows directly from the strict
Dirichlet PPW inequality and Robin-to-Dirichlet eigenvalue convergence
(Proposition~\ref{prop:endpoint-regimes}). The threshold depends on
$\Omega$.

In dimension two, Theorem~\ref{thm:ball-bound} also has a conditional interpretation
for an elastic drum of prescribed area and known boundary parameter.
Among domains satisfying $\widetilde R\geq R$, the ratio of the
first two frequencies equals the corresponding disk value if and only
if the drum is circular. Equivalently, the relative spectral gap obeys
\[
 \frac{\mu_2-\mu_1}{\mu_1}
 \leq\frac{\mu_2(B,\sigma)-\mu_1(B,\sigma)}{\mu_1(B,\sigma)}.
\]
This conditional circularity criterion is related to Yau's Problem~77
\cite{Yau} and Kac's inverse spectral question \cite{Kac}.

Varying the Robin parameter on the comparison ball gives an upper
bound without the condition $\widetilde R\geq R$. Unlike a bound
depending only on volume and $\sigma$, this estimate retains the
interior superlevel-set information through $\widetilde R$.
\\ \\
\refstepcounter{thm}\label{thm:matched-bound}\textbf{Theorem~\thethm.}
\emph{Let $\Omega\subset\mathbb R^N$, $N\geq2$, be a bounded connected
$C^2$ domain with $|\Omega|=|B|$, and let $\sigma>0$. Then
\[
 \frac{\mu_2}{\mu_1}
 \leq
 \begin{cases}
  \displaystyle\frac{\mu_2(B,\tau)}{\mu_1(B,\tau)},
     &\mu_1\widetilde R^2<j_{\nu,1}^2,\\[3mm]
  \displaystyle\frac{j_{\nu+1,1}^2}{j_{\nu,1}^2},
     &\mu_1\widetilde R^2\geq j_{\nu,1}^2,
 \end{cases}
\]
where, in the first case, $\tau>0$ is the unique parameter satisfying}
\begin{equation}\label{eq:matched-parameter}
 \mu_1(B,\tau)=\mu_1\widetilde R^2.
\end{equation}
\emph{Equality is attained when $\Omega$ is a unit ball, up to
translation; such a domain belongs to the first case.}\\

The parameter $\tau$ is well defined because $\mu_1(B,\tau)$ increases
strictly from $0$ to $j_{\nu,1}^2$. When $\widetilde R<R$, one has
$\tau<\sigma$, and Theorem~\ref{thm:ball-monotonicity} places the new
bound strictly between the ball ratio at $\sigma$ and
\eqref{eq:prior-bound}. For $\widetilde R\geq R$, it recovers the
upper bound in Theorem~\ref{thm:ball-bound}, with a strict improvement
when $\widetilde R>R$. The bounds coincide when $\widetilde R=R$;
the equality characterization and the comparison of the bounds are
proved in Section~\ref{sec:proof-ball-bound}.

Laugesen \cite[Section~II.A, Conjecture~B]{Laugesen} conjectured that,
for every fixed bounded Lipschitz domain $\Omega$, the ratio
$\mu_2(\Omega,\sigma)/\mu_1(\Omega,\sigma)$ is decreasing in $\sigma>0$.
Theorem~\ref{thm:ball-monotonicity} proves this assertion for balls in
every dimension $N\geq2$, with a strictly negative derivative.
\\ \\
\refstepcounter{thm}\label{thm:ball-monotonicity}\textbf{Theorem~\thethm.}
\emph{For $N\geq2$, the functions
\[
 \sigma\longmapsto\mu_2(B,\sigma)-\mu_1(B,\sigma),
 \,\,\,
 \sigma\longmapsto\frac{\mu_2(B,\sigma)}{\mu_1(B,\sigma)}
\]
are continuously differentiable on $(0,+\infty)$, with strictly
positive and strictly negative derivatives, respectively. Moreover,
\[
 \lim_{\sigma\downarrow0}\frac{\mu_2(B,\sigma)}{\mu_1(B,\sigma)}
 =+\infty,
 \,\,\,
 \lim_{\sigma\to+\infty}\frac{\mu_2(B,\sigma)}{\mu_1(B,\sigma)}
 =\frac{j_{\nu+1,1}^2}{j_{\nu,1}^2}.
\]
Consequently, the ratio takes each value in
$\bigl(j_{\nu+1,1}^2/j_{\nu,1}^2,+\infty\bigr)$ exactly once.}\\

The gap assertion is the ball case of Smits' conjecture, recorded as
Conjecture~A in \cite{Laugesen}; its disk case was already noted there.
Proposition~\ref{prop:angular-monotonicity} compares the lowest radial
eigenvalues of any two spherical harmonic degrees, giving both strict
derivative signs by a common argument. These are comparisons between
angular branches, not arbitrary pairs in the ordered spectrum.
The related work \cite{DaiSunBoxes} proves spectral-ratio monotonicity
on rectangular boxes and gives connected domains on which the ratio
is not monotone in the boundary parameter. Those examples concern
parameter variation on a fixed domain, whereas
Theorem~\ref{thm:unbounded-ratio} below varies the domain at a fixed
parameter.

The ratio monotonicity also clarifies the role of the upper bound $3$
in the planar problem. Payne and Schaefer \cite[Section~III]{PayneSchaefer}
proved that $\mu_2/\mu_1\leq3$ for a range of sufficiently large positive
Robin parameters. For the unit disk, Theorem~\ref{thm:ball-monotonicity}
identifies the exact threshold: since $j_{1,1}^2/j_{0,1}^2<3$, there is
a unique $\overline\sigma>0$ such that
\[
 \frac{\mu_2(B,\overline\sigma)}{\mu_1(B,\overline\sigma)}=3,
 \,\,\,
 \frac{\mu_2(B,\sigma)}{\mu_1(B,\sigma)}\leq3
 \quad\Longleftrightarrow\quad \sigma\geq\overline\sigma.
\]
For a planar domain of the same area satisfying $\widetilde R\geq R$,
Theorem~\ref{thm:ball-bound} therefore gives $\mu_2/\mu_1\leq3$ whenever
$\sigma\geq\overline\sigma$, and a strictly smaller upper bound when
$\sigma>\overline\sigma$. More generally, for every
$c>j_{\nu+1,1}^2/j_{\nu,1}^2$, there is a unique $\sigma_c>0$ at which
the ball ratio equals $c$. The same comparison gives
$\mu_2/\mu_1\leq c$ for $\sigma\geq\sigma_c$, provided
$\widetilde R\geq R$ at that parameter. Figure~\ref{fig:ball-ratio}
shows this dependence on $\sigma$, and
Table~\ref{tab:appendix-positive-spectrum} in Appendix
gives numerical examples, including the disk value $3.66722$ at $\sigma=1$.

\begin{figure}[htbp]
\centering

\resizebox{0.48\textwidth}{!}{%
\begingroup%
\setlength{\unitlength}{1bp}%
\fontsize{12}{14.4}\selectfont%
\begin{picture}(260,180)%
\PPWxaxis{50}{44}{245}%
\PPWyaxis{78}{16}{165}%
\linethickness{0.65bp}%
\multiput(78,66.748)(11,0){13}{\line(1,0){7}}%
\linethickness{1.35bp}%
\qbezier(80.8496,140.316)(81.8115,137.8505)(82.7734,133.54)
\qbezier(82.7734,133.54)(83.7353,129.7164)(84.6973,127.0544)
\qbezier(84.6973,127.0544)(86.8616,121.1873)(89.0259,117.1808)
\qbezier(89.0259,117.1808)(93.3546,109.478)(97.6832,104.5968)
\qbezier(97.6832,104.5968)(106.3405,95.2013)(114.9978,90.0284)
\qbezier(114.9978,90.0284)(132.072,80.0115)(149.1461,76.1376)
\qbezier(149.1461,76.1376)(166.2202,72.101)(183.2943,70.8136)
\qbezier(183.2943,70.8136)(200.6089,69.4282)(217.9235,69.41)
\put(73,154){\makebox(0,0)[r]{$\displaystyle\frac{k_{\nu+1,1}^{2}}{k_{\nu,1}^{2}}$}}
\put(73,66.748){\makebox(0,0)[r]{$\displaystyle\frac{j_{\nu+1,1}^{2}}{j_{\nu,1}^{2}}$}}
\put(72,32){\makebox(0,0)[r]{$0$}}
\put(242,32){\makebox(0,0)[r]{$\sigma$}}
\end{picture}%
\endgroup%
}%

\caption{{The ratio $k_{\nu+1,1}(\sigma)^2/k_{\nu,1}(\sigma)^2$
on the unit ball and its Dirichlet limit.}}
\label{fig:ball-ratio}
\end{figure}

The condition $\widetilde R\geq R$ cannot be omitted from the ball
comparison for general domains. In Bucur's counterexample, recorded by Laugesen
\cite[Section~III.D, after Conjecture~G]{Laugesen}, a small ball supplies
the first two eigenvalues, while finitely many smaller disjoint balls
supply the remaining volume with higher first eigenvalues. Scaling
makes the ratio exceed that of the equal-volume ball. The open set has
smooth boundary but is disconnected. Laugesen suggests joining its
components without significantly changing the first two eigenvalues.
Nahon \cite[Section~1.2]{Nahon2022} also records the unboundedness of the
ratio without connectedness. Using the same separation of spectral
scales, a cube attached to a long thin tube gives an explicit connected
construction with quantitative bounds and smooth approximation at
exactly fixed volume.
\\\\
\refstepcounter{thm}\label{thm:unbounded-ratio}\textbf{Theorem~\thethm.}
\emph{For every $N\geq2$, $V>0$, $\sigma>0$, and $K>0$, there exists a
bounded connected domain $\Omega\subset\mathbb R^N$ with $C^\infty$ boundary
such that}
\[
 |\Omega|=V,\,\,\,
 \frac{\mu_2(\Omega,\sigma)}{\mu_1(\Omega,\sigma)}>K.
\]
\emph{In dimension two, $\Omega$ can also be chosen simply connected.}
\\

Theorem~\ref{thm:unbounded-ratio} shows that smoothness and connectedness
do not by themselves ensure a ball comparison. The main difficulty in
passing from the disconnected examples to a connected domain is to add
the prescribed volume without introducing further low eigenvalues.
For Robin conditions, a thin attachment also changes the boundary
energy, so its spectral effect cannot be inferred from its width alone.
Our construction controls that energy directly: the added part has a
large transverse energy, while a trial function concentrated in the
main part keeps the first eigenvalue small relative to the second.
This gives quantitative bounds on connected domains, and the strict
inequality survives smoothing and exact volume normalization.

For the positive comparisons in Theorems~\ref{thm:ball-bound}
and~\ref{thm:matched-bound}, the principal difficulty is instead the
nonconstant boundary value of the first eigenfunction. The methods of
Ashbaugh and Benguria
\cite{AshbaughBenguria91,AshbaughBenguria92,AshbaughBenguria92CMP}
and Chiti \cite{Chiti1982,Chiti19822} cannot be applied on every
superlevel set: a set meeting $\partial\Omega$ contributes an exterior
boundary term. The boundary decomposition used in Robin rearrangement
arguments \cite{DaiShi,Alvino,Alvino1} removes this term above
$u_{1,pM}$, but leaves only a restricted comparison interval.
Moreover, the comparison eigenfunction is positive at its boundary,
so the usual Dirichlet crossing argument must allow a no-crossing
case as well. Following the Robin comparison of Dai and Shi
\cite[Theorem~3.1]{DaiShi}, a maximum-splicing argument handles both
possibilities while retaining the endpoint contribution.

The refinement in Theorem~\ref{thm:matched-bound} comes from varying
the Robin parameter on the comparison ball, rather than keeping its
unit-ball value fixed. The parameter is optimized under the first-eigenvalue matching
condition and the restriction on the comparison radius. This improves
the upper bound without extending the rearrangement inequality beyond
its valid interval. The necessary radial weight inequalities follow from a
Wronskian identity and a Riccati equation adapted to the Robin
boundary condition. For Theorem~\ref{thm:ball-monotonicity}, the
increase of each eigenvalue alone is insufficient to compare their
ratio. A logarithmic-derivative identity, together with a comparison
of the lowest angular branches, gives the required strict derivative
sign. Thus the sharper bound and the parameter monotonicity rely on
different uses of the ball eigenfunctions.

\par\endgroup

\begingroup
\section{Strict monotonicity of $w(r)$ and $B(r)$}\label{sec:weights}

Fix $\sigma>0$ and retain $\alpha$ and $\beta$ from the introduction.
The radial weights used in the eigenvalue comparison are
\begin{equation*}
 w(r)=\frac{J_{\nu+1}(\beta r)}{J_\nu(\alpha r)},
 \,\,\, B(r)=w'(r)^2+(N-1)\frac{w(r)^2}{r^2},
 \,\,\, 0<r\leq1.
\end{equation*}
Since $\alpha<j_{\nu,1}$ and $\beta<j_{\nu+1,1}$, both Bessel
functions in the quotient are positive on this interval.
The power series for $J_\eta$, with $\eta\geq0$, gives
\[
 J_\eta(x)=\frac{(x/2)^\eta}{\Gamma(\eta+1)}
 \left(1-\frac{x^2}{4(\eta+1)}+O(x^4)\right)
 \,\,\, (x\to0^+);
\]
see \cite[Eq.~10.2.2]{DLMF}. Applying this formula to the numerator
and denominator, and using
$\Gamma(\nu+2)=(\nu+1)\Gamma(\nu+1)$, yields
\[
 w(r)=\frac{\beta^{\nu+1}}{2(\nu+1)\alpha^\nu}\,r
 \left[1+\frac14\left(\frac{\alpha^2}{\nu+1}
                    -\frac{\beta^2}{\nu+2}\right)r^2+O(r^4)\right].
\]
Thus $w$ extends smoothly to zero, with
\[
 w(0)=0,\,\,\,
 w'(0)=\lim_{r\to0^+}\frac{w(r)}r
       =\frac{\beta^{\nu+1}}{2(\nu+1)\alpha^\nu}>0.
\]
The apparent singularity of $B$ is also removable, and
\[
 B(0)=Nw'(0)^2
     =\frac{N\beta^{2\nu+2}}{4(\nu+1)^2\alpha^{2\nu}}>0.
\]
In particular, $B$ is bounded near zero, while $w(r)=O(r)$ makes the
radial trial functions $w(\gamma|x|)x_i/|x|$ regular at the origin
for every fixed $\gamma>0$.

Set
\[
 u(r)=r^{-\nu}J_\nu(\alpha r),\,\,\,
 v(r)=r^{-\nu}J_{\nu+1}(\beta r),\,\,\,
 q(r)=\frac{rw'(r)}{w(r)}.
\]
Then $w=v/u$. Taking its logarithmic derivative gives
\[
 q(r)=r\left[\beta\frac{J'_{\nu+1}(\beta r)}{J_{\nu+1}(\beta r)}
             -\alpha\frac{J'_\nu(\alpha r)}{J_\nu(\alpha r)}\right].
\]
The Bessel derivative formulas
\[
 J_\eta'(x)=\frac{\eta}{x}J_\eta(x)-J_{\eta+1}(x),\,\,\,
 J_{\eta+1}'(x)=J_\eta(x)-\frac{\eta+1}{x}J_{\eta+1}(x)
\]
from \cite[Eq.~10.6.2]{DLMF} therefore give
\[
 q(r)=1+\alpha r\frac{J_{\nu+1}(\alpha r)}{J_\nu(\alpha r)}
        -\beta r\frac{J_{\nu+2}(\beta r)}{J_{\nu+1}(\beta r)}.
\]
The expansion of $w$ shows that $w(r)/r$ is analytic in $r^2$ with
a positive constant term. Hence $q$ extends analytically to zero,
and $rw'(r)/w(r)=1+O(r^2)$ gives $q(0)=1$.
At the other endpoint, the defining root equations imply
\[
 \alpha\frac{J_{\nu+1}(\alpha)}{J_\nu(\alpha)}=\sigma,
 \,\,\,
 \beta\frac{J_{\nu+2}(\beta)}{J_{\nu+1}(\beta)}=\sigma+1.
\]
Consequently, $q(1)=1+\sigma-(\sigma+1)=0$ and $w'(1)=0$.
The same identities show that
$u'(1)+\sigma u(1)=v'(1)+\sigma v(1)=0$.
Extending $w$ constantly to $[1,+\infty)$ therefore gives a $C^1$
function, and
\begin{equation}\label{eq:weight-tail}
 B(r)=(N-1)\frac{w(1)^2}{r^2},\,\,\, r\geq1.
\end{equation}

The endpoint behavior of $q$ depends on two comparisons between
$\alpha$ and $\beta$. For these comparisons, the radial Rayleigh
quotient is
\[
 \mathcal Q_\sigma[f]=
 \frac{\displaystyle\int_0^1r^{N-1}|f'(r)|^2\,\mathrm{d}r
                    +\sigma|f(1)|^2}
      {\displaystyle\int_0^1r^{N-1}|f(r)|^2\,\mathrm{d}r}.
\]
The variational principle gives $\alpha^2\leq\mathcal Q_\sigma[f]$
for every nonzero radial $f\in H^1(B)$; the trial function need not
satisfy a boundary condition.
\\\\
\refstepcounter{lem}\label{lem:ball-gap}\textbf{Lemma~\thelem.}
\emph{For every $\sigma>0$,}
\begin{equation}\label{eq:gap-elementary}
 \beta^2-\alpha^2>N-1.
\end{equation}
\emph{In particular, $q'(1)=N-1+\alpha^2-\beta^2<0$.}\\

\noindent\textbf{Proof.}
The Bessel equation gives the radial equations
\[
 -(r^{N-1}u')'=\alpha^2r^{N-1}u,\,\,\,
 -(r^{N-1}v')'+(N-1)r^{N-3}v=\beta^2r^{N-1}v.
\]
The function $v$ is an admissible trial function in
$\mathcal Q_\sigma$: near zero, $v(r)=O(r)$ and $v'(r)=O(1)$, so
its radial extension belongs to $H^1(B)$ for every $N\geq2$.
These estimates also give
$\lim_{r\to0^+}r^{N-1}v(r)v'(r)=0$.
Multiplying the equation for $v$ by $v$ and integrating by parts,
we obtain
\begin{align*}
 &\int_0^1r^{N-1}v'(r)^2\,\mathrm{d}r+\sigma v(1)^2\\
 &=\bigl[r^{N-1}vv'\bigr]_0^1
   -\int_0^1v(r)(r^{N-1}v'(r))'\,\mathrm{d}r+\sigma v(1)^2\\
 &=\beta^2\int_0^1r^{N-1}v(r)^2\,\mathrm{d}r
   -(N-1)\int_0^1r^{N-3}v(r)^2\,\mathrm{d}r.
\end{align*}
The boundary terms at $r=1$ cancel because $v'(1)=-\sigma v(1)$.
Both integrals are finite, also when $N=2$. Therefore
\[
 \alpha^2\leq\mathcal Q_\sigma[v]
 =\beta^2-(N-1)
 \frac{\displaystyle\int_0^1r^{N-3}v(r)^2\,\mathrm{d}r}
      {\displaystyle\int_0^1r^{N-1}v(r)^2\,\mathrm{d}r}
 <\beta^2-(N-1),
\]
since $r^{-2}>1$ and $v(r)>0$ on $(0,1)$.

To obtain the endpoint derivative, substitute $v=uw$ into the
second radial equation and use the first. This gives
\[
 w''+\left(2\frac{u'}u+\frac{N-1}{r}\right)w'
 +\left(\beta^2-\alpha^2-\frac{N-1}{r^2}\right)w=0.
\]
Since $w'(1)=0$ and $w(1)>0$, differentiation of $q=rw'/w$ yields
\[
 q'(1)=\frac{w''(1)}{w(1)}=N-1+\alpha^2-\beta^2<0.
\]
\qed
\\ \\
\refstepcounter{rek}\label{rem:ball-gap}\textbf{Remark~\therek.}
Inequality \eqref{eq:gap-elementary} gives a lower bound for the first
Robin spectral gap on the unit ball, uniform for all $\sigma>0$.
By \eqref{eq:robin-scaling},
\[
 \mu_2(B_\rho,\sigma)-\mu_1(B_\rho,\sigma)
 >\frac{N-1}{\rho^2},\,\,\, \rho>0,\quad\sigma>0.
\]
Taking the Dirichlet limit $\sigma\to+\infty$ gives
\[
 \lambda_2(B)-\lambda_1(B)
 =j_{\nu+1,1}^2-j_{\nu,1}^2\geq N-1.
\]
For comparison, Singer, Wong, S.-T.~Yau and S.~S.-T.~Yau
\cite[Theorem~1.1]{SingerWongYauYau} proved the lower bound
$\lambda_2-\lambda_1\geq\pi^2/(4d^2)$ on smooth bounded convex
domains of diameter $d$. On the unit ball, $d=2$, so the bound
$N-1$ is stronger than $\pi^2/16$ for every $N\geq2$.
The sharp general convex-domain bound is $3\pi^2/d^2$, due to
Andrews and Clutterbuck \cite[Corollary~1.4]{AndrewsClutterbuck}.
On the unit ball it gives $3\pi^2/4$, which exceeds $N-1$ for
$2\leq N\leq8$ and is smaller than $N-1$ for $N\geq9$.
Proposition~\ref{prop:angular-monotonicity} extends the Robin lower
bound to the lowest radial eigenvalues of any two spherical harmonic
degrees.
Theorem~\ref{thm:ball-monotonicity} and the Neumann limit give the
more precise endpoint bound
\[
 \mu_2(B,\sigma)-\mu_1(B,\sigma)>\mu_2(B,0),\qquad\sigma>0,
\]
whose right-hand side is the infimum of the gap over positive
parameters. The elementary estimate \eqref{eq:gap-elementary}
is useful here because its explicit constant controls the endpoint
behavior of the radial weights.
\\ \\
\refstepcounter{lem}\label{lem:origin-comparison}\textbf{Lemma~\thelem.}
\emph{For every $\sigma>0$,}
\begin{equation*}
 \frac{\alpha^2}{\nu+1}-\frac{\beta^2}{\nu+2}<0.
\end{equation*}
\emph{Moreover, $q'(0)=0$ and
$q''(0)=\alpha^2/(\nu+1)-\beta^2/(\nu+2)<0$.}
\\\\
\textbf{Proof.}
Consider the positive radial function
\[
 \phi(r)=\frac{v(r)}r=r^{-\nu-1}J_{\nu+1}(\beta r).
\]
Its power series contains only even powers of $r$, so it is smooth
at the origin and belongs to $H^1(B)$. The derivative formula gives
\[
 \phi'(r)=-\beta r^{-\nu-1}J_{\nu+2}(\beta r),\,\,\,
 \phi'(1)=-(\sigma+1)\phi(1).
\]
In particular, $\phi$ is not a multiple of the first radial
eigenfunction $u$: such a multiple would satisfy
$\phi'(1)+\sigma\phi(1)=0$, whereas the displayed condition gives
$\phi'(1)+\sigma\phi(1)=-\phi(1)<0$. The Rayleigh principle and
simplicity of the first eigenvalue thus imply
$\alpha^2<\mathcal Q_\sigma[\phi]$.

Since $N-1=2\nu+1$, the preceding derivative identity gives
\[
 r^{N-1}\phi'(r)=-\beta r^\nu J_{\nu+2}(\beta r).
\]
Differentiating once more,
\begin{align*}
 \bigl(r^\nu J_{\nu+2}(\beta r)\bigr)'
 &=\nu r^{\nu-1}J_{\nu+2}(\beta r)
   +\beta r^\nu J'_{\nu+2}(\beta r)\\
 &=2(\nu+1)r^{\nu-1}J_{\nu+2}(\beta r)
   -\beta r^\nu J_{\nu+3}(\beta r).
\end{align*}
The recurrence relation \cite[Eq.~10.6.1]{DLMF}
\[
 J_{\nu+1}(\beta r)+J_{\nu+3}(\beta r)
 =\frac{2(\nu+2)}{\beta r}J_{\nu+2}(\beta r)
\]
then yields
\begin{align*}
 -(r^{N-1}\phi')'
 &=2\beta(\nu+1)r^{\nu-1}J_{\nu+2}(\beta r)
   -\beta^2r^\nu J_{\nu+3}(\beta r)\\
 &=\frac{\beta^2(\nu+1)}{\nu+2}\,r^{N-1}\phi
   -\frac{\beta^2}{\nu+2}\,r^\nu J_{\nu+3}(\beta r).
\end{align*}

Put $I=\int_0^1r^{N-1}\phi^2\,\mathrm{d}r>0$.
Since $\phi(r)=O(1)$ and $\phi'(r)=O(r)$ at zero, the endpoint
term there vanishes. At $r=1$, the boundary condition gives
\begin{align*}
 &\int_0^1r^{N-1}\phi'(r)^2\,\mathrm{d}r+\sigma\phi(1)^2\\
 &=-\int_0^1\phi(r)(r^{N-1}\phi'(r))'\,\mathrm{d}r
   +\bigl[r^{N-1}\phi\phi'\bigr]_0^1+\sigma\phi(1)^2\\
 &=-\int_0^1\phi(r)(r^{N-1}\phi'(r))'\,\mathrm{d}r-\phi(1)^2.
\end{align*}
Consequently,
\begin{align*}
 \alpha^2<\mathcal Q_\sigma[\phi]
 &=\frac{\beta^2(\nu+1)}{\nu+2}
   -\frac{\beta^2}{(\nu+2)I}
       \int_0^1\phi(r)r^\nu J_{\nu+3}(\beta r)\,\mathrm{d}r
   -\frac{\phi(1)^2}{I}\\
 &<\frac{\beta^2(\nu+1)}{\nu+2}.
\end{align*}
Here $J_{\nu+3}(\beta r)>0$ on $(0,1]$, since
$\beta<j_{\nu+1,1}<j_{\nu+3,1}$; see
\cite[Section~10.21]{DLMF}. This proves the required comparison.
Finally, taking the logarithmic derivative of the expansion of $w$
at the beginning of this section gives
\[
 q(r)=1+\frac12\left(\frac{\alpha^2}{\nu+1}
                    -\frac{\beta^2}{\nu+2}\right)r^2+O(r^4).
\]
Thus $q'(0)=0$ and
$q''(0)=\alpha^2/(\nu+1)-\beta^2/(\nu+2)<0$.
\qed\\

The two lemmas give $q'(1)<0$ and $q''(0)<0$. Together with the
radial differential equations, these endpoint inequalities yield the
strict weight monotonicity below, strengthening
\cite[Lemma~A.3]{DaiShi}.
\\\\
\refstepcounter{prop}\label{prop:weight-monotonicity}\textbf{Proposition~\theprop.}
\emph{For every $\sigma>0$, $w$ is strictly increasing on $[0,1]$ and
nondecreasing on $[0,+\infty)$. The function $B$ is strictly decreasing
on $[0,+\infty)$, with $B'(r)<0$ for $0<r<1$. Moreover,
$w(r)>0$ for $r>0$ and $B(r)>0$ for $r\geq0$.}\\
\\
\textbf{Proof.}
Put $D=\beta^2-\alpha^2$. The radial equations, the expansions
$u(r)=u_0+O(r^2)$ and $v(r)=v_1r+O(r^3)$, and the common Robin
condition give
\begin{equation*}
 \begin{gathered}
 W=r^{N-1}(uv'-u'v),\,\,\, W(0)=W(1)=0,\\
 W'=r^{N-1}uv\left(\frac{N-1}{r^2}-D\right).
 \end{gathered}
\end{equation*}
By Lemma~\ref{lem:ball-gap}, $D>N-1$. Hence $W'$ changes sign exactly once,
from positive to negative. Since $W$ vanishes at both endpoints,
we obtain $W>0$ on $(0,1)$. Consequently,
\[
 w'(r)=\frac{W(r)}{r^{N-1}u(r)^2}>0\,\,\,(0<r<1).
\]

The identity
\[
 B(r)=\bigl(q(r)^2+N-1\bigr)\left(\frac{w(r)}r\right)^2
\]
reduces the monotonicity of $B$ to estimates on $q$. Set
\[
 h(r)=-\frac{u'(r)}{u(r)}
      =\alpha\frac{J_{\nu+1}(\alpha r)}{J_\nu(\alpha r)}.
\]
The positivity of $w'$ gives $q>0$ on $(0,1)$, and the endpoint
values are $q(0)=1$ and $q(1)=0$. Since $\alpha<j_{\nu,1}$,
taking the logarithmic derivative of the infinite product for
$J_\nu$ in \cite[Eq.~10.21.15]{DLMF} gives
\begin{equation*}
 h(r)=\sum_{m=1}^{+\infty}
 \frac{2\alpha^2r}{j_{\nu,m}^2-\alpha^2r^2}.
\end{equation*}
The series and its derivatives converge uniformly on $[0,1]$.
Differentiating term by term, we see that $h(r)/r$ is strictly
increasing and $h$ is strictly convex on $(0,1)$.

Differentiating $q=rw'/w$ and using the equation for $w$ gives
\begin{align*}
 q'(r)
 &=\frac{q(r)}r+r\frac{w''(r)}{w(r)}-\frac{q(r)^2}r\\
 &=-Dr+\frac{(1-q(r))(q(r)+N-1)}r+2q(r)h(r).
\end{align*}
Write this as $q'(r)=F(r,q(r))$, where
\[
 F(r,c)=-Dr+\frac{(1-c)(c+N-1)}r+2c\,h(r).
\]
If $q>1$ somewhere, choose $c>1$ below its maximum.
Since $q(0)=1$ and $q(1)=0$, there is an interval $(r_1,r_2)$ on
which $q>c$, with $q(r_1)=q(r_2)=c$. Hence $q'(r_1)\geq0$ and
$q'(r_2)\leq0$. On the other hand,
\[
 \frac{F(r,c)}r=-D+\frac{(1-c)(c+N-1)}{r^2}
                  +2c\frac{h(r)}r
\]
is strictly increasing. This contradicts
$F(r_1,c)/r_1\geq0\geq F(r_2,c)/r_2$, and therefore $0\leq q\leq1$.

Suppose that $q'(r_2)>0$ for some
$r_2\in(0,1)$, and let $c=q(r_2)$. The bounds just proved imply
$c\in(0,1)$. Since $q(0)=1$ and $q(1)=0$, the function must pass
from above $c$ to below $c$ both before and after $r_2$. Thus there
are points $r_1<r_2<r_3$ such that $q(r_i)=c$ and
$q'(r_1),q'(r_3)\leq0$. For this fixed $c$, the function
$F(\cdot,c)$ is strictly convex, since $(1-c)(c+N-1)>0$ and $h$ is
convex. Its values at $r_1$ and $r_3$ are nonpositive, so strict
convexity gives $F(r_2,c)<0$. This contradicts
$F(r_2,c)=q'(r_2)>0$, and proves that $q'\leq0$.

By Lemma~\ref{lem:origin-comparison}, $q(r)<1$ for all sufficiently small positive $r$.
Since $q'\leq0$, this strict inequality holds throughout $(0,1]$.
Finally, differentiating $B$ gives
\begin{equation*}
 B'(r)=2\left(\frac{w(r)}r\right)^2
 \left[qq'+\frac{q-1}{r}(q^2+N-1)\right]<0,
 \,\,\, 0<r<1.
\end{equation*}
The last inequality follows from $qq'\leq0$ and $q<1$.
The positivity of $w$ and $B$ follows from their definitions and
the expansion at the origin. Together with continuity and
\eqref{eq:weight-tail}, these facts prove the desired conclusions
on $[0,+\infty)$.
\qed\\

Figure~\ref{fig:radial-weights} illustrates these two functions.
\begin{figure}[htbp]
\centering
\resizebox{0.48\textwidth}{!}{%
\begingroup%
\setlength{\unitlength}{1bp}%
\fontsize{12}{14.4}\selectfont%
\begin{picture}(260,180)%
\PPWxaxis{87}{44}{248}%
\PPWyaxis{102}{16}{165}%
\linethickness{0.65bp}%
\multiput(102,125)(9,0){14}{\line(1,0){5.5}}%
\multiput(224,44)(0,9){9}{\line(0,1){5.5}}%
\linethickness{1.35bp}%
\qbezier(102,44)(104.4989,48.638)(106.9978,52.3538)
\qbezier(106.9978,52.3538)(108.9781,55.4212)(110.9584,58.3971)
\qbezier(110.9584,58.3971)(114.6361,63.7061)(118.3138,68.0092)
\qbezier(118.3138,68.0092)(125.9522,77.1629)(133.5905,84.2356)
\qbezier(133.5905,84.2356)(148.5842,97.9411)(163.578,106.3626)
\qbezier(163.578,106.3626)(193.789,122.8484)(224,125)
\put(96,157){\makebox(0,0)[r]{$w(r)$}}
\put(96,125){\makebox(0,0)[r]{$\displaystyle\frac{J_{\nu+1}(\beta)}{J_\nu(\alpha)}$}}
\put(96,32){\makebox(0,0)[r]{$0$}}
\put(224,32){\makebox(0,0)[c]{$1$}}
\put(245,32){\makebox(0,0)[c]{$r$}}
\end{picture}%
\endgroup%
}
\resizebox{0.48\textwidth}{!}{%
\begingroup%
\setlength{\unitlength}{1bp}%
\fontsize{12}{14.4}\selectfont%
\begin{picture}(260,180)%
\PPWxaxis{87}{44}{248}%
\PPWyaxis{102}{16}{165}%
\linethickness{0.65bp}%
\multiput(102,78)(9,0){14}{\line(1,0){5.5}}%
\multiput(224,44)(0,9){4}{\line(0,1){5.5}}%
\linethickness{1.35bp}%
\qbezier(102,125)(164.5641,125)(224,78)%
\put(96,157){\makebox(0,0)[r]{$B(r)$}}
\put(96,125){\makebox(0,0)[r]{{$\displaystyle\frac{2(\nu+1)\beta^{2\nu+2}}{4(\nu+1)^2\alpha^{2\nu}}$}}}
\put(96,78){\makebox(0,0)[r]{$(2\nu+1)w(1)^2$}}
\put(96,32){\makebox(0,0)[r]{$0$}}
\put(224,32){\makebox(0,0)[c]{$1$}}
\put(245,32){\makebox(0,0)[c]{$r$}}
\end{picture}%
\endgroup%
}%
\caption{{The functions $w(r)$ and $B(r)$ on $[0,1]$.}}
\label{fig:radial-weights}
\end{figure}
For another positive Robin parameter $\tau$, write $w_\tau$ and
$B_\tau$ for the same functions with $\sigma$ replaced by $\tau$.

For a fixed domain $\Omega$, the applicability of the ball comparison
also depends on the boundary maximum of its first eigenfunction.
The following estimate shows that this maximum tends to zero as
$\sigma$ increases, so the interior superlevel sets fill the domain
and the condition $\widetilde R\geq R$ eventually holds.
\\\\
\refstepcounter{prop}\label{prop:boundary-maximum}\textbf{Proposition~\theprop.}
\emph{Let $\Omega$ be a fixed bounded connected $C^2$ domain with
$|\Omega|=C_N$. For each $\sigma>0$, let $u_\sigma$ be its positive
first Robin eigenfunction normalized by $\|u_\sigma\|_{L^2(\Omega)}=1$.
There is a constant $C_\Omega$ such that}
\begin{equation}\label{eq:boundary-maximum-decay}
 \max_{\partial\Omega}u_\sigma\leq\frac{C_\Omega}{\sigma},
 \,\,\, \sigma\geq1.
\end{equation}
\emph{Moreover, $u_\sigma$ converges uniformly on $\overline\Omega$
to the positive $L^2$-normalized first Dirichlet eigenfunction, and}
\[
 \widetilde R(\sigma)\longrightarrow1,
 \,\,\,
 R(\sigma)\longrightarrow\frac{j_{\nu,1}}{\sqrt{\lambda_1(\Omega)}}
 \quad\text{as }\sigma\to+\infty.
\]
\emph{Consequently, $\widetilde R(\sigma)\geq R(\sigma)$ for all
sufficiently large $\sigma$, with strict inequality if $\Omega$
is not a ball.}\\

\noindent\textbf{Proof.}
Write $\lambda_1=\lambda_1(\Omega)$. A first Dirichlet eigenfunction
has zero boundary trace, so its Robin Rayleigh quotient equals
$\lambda_1$. Thus $0<\mu_1(\Omega,\sigma)\leq\lambda_1$.
For each fixed $\sigma\geq1$, the regularity
$u_\sigma\in C^1(\overline\Omega)$ makes $u_\sigma^{p-1}$ an
admissible test function for every finite $p\geq2$.
Testing the weak equation and using
$\nabla(u_\sigma^{p/2})=(p/2)u_\sigma^{p/2-1}\nabla u_\sigma$ gives
\begin{align*}
 \frac{4(p-1)}{p^2}
 \int_\Omega\bigl|\nabla u_\sigma^{p/2}\bigr|^2\,\mathrm{d}x
 &=(p-1)\int_\Omega u_\sigma^{p-2}|\nabla u_\sigma|^2\,\mathrm{d}x\\
 &=\mu_1(\Omega,\sigma)\int_\Omega u_\sigma^p\,\mathrm{d}x
   -\sigma\int_{\partial\Omega}u_\sigma^p\,\mathrm{d}S\\
 &\leq\lambda_1\int_\Omega u_\sigma^p\,\mathrm{d}x.
\end{align*}
The boundary term has the favorable sign because $\sigma>0$.
Let $\chi=N/(N-2)$ for $N\geq3$, and $\chi=2$ for $N=2$.
The Sobolev inequality on the fixed domain has the form
\[
 \|f\|_{L^{2\chi}(\Omega)}^2
 \leq S_\Omega\left(\|\nabla f\|_{L^2(\Omega)}^2
                         +\|f\|_{L^2(\Omega)}^2\right),
 \,\,\, f\in H^1(\Omega),
\]
where $S_\Omega\geq1$ is independent of $\sigma$ and $p$.
Applying it to $f=u_\sigma^{p/2}$ and taking the $p$th root yields
\[
 \|u_\sigma\|_{L^{p\chi}(\Omega)}
 \leq\left[S_\Omega\left(1+\frac{\lambda_1p^2}{4(p-1)}\right)\right]^{1/p}
       \|u_\sigma\|_{L^p(\Omega)}.
\]
Set $p_j=2\chi^j$. Iterating from $p_0=2$ and using
$\|u_\sigma\|_{L^2(\Omega)}=1$, we obtain
\[
 \|u_\sigma\|_{L^{p_m}(\Omega)}
 \leq\prod_{j=0}^{m-1}
 \left[S_\Omega\left(1+\frac{\lambda_1p_j^2}{4(p_j-1)}\right)\right]^{1/p_j}.
\]
Since $p^2/[4(p-1)]\leq p/2$ for $p\geq2$ and
\[
 \sum_{j=0}^{+\infty}\frac{1+\log p_j}{p_j}<+\infty,
\]
the product is bounded independently of $m$ and $\sigma$.
Letting $m\to+\infty$ gives
\[
 \|u_\sigma\|_{L^\infty(\Omega)}\leq M,\,\,\, \sigma\geq1,
\]
where $M$ depends only on the fixed domain and the dimension.

Let $\psi$ solve $-\Delta\psi=1$ in $\Omega$ and $\psi=0$ on
$\partial\Omega$. Boundary regularity gives
$\psi\in C^1(\overline\Omega)$; set
$K=\|\partial_{\mathbf n}\psi\|_{L^\infty(\partial\Omega)}$.
The function
\[
 V_\sigma=M\lambda_1\left(\psi+\frac K\sigma\right)
\]
satisfies
\[
 -\Delta V_\sigma=M\lambda_1
       \geq\mu_1(\Omega,\sigma)u_\sigma,
 \,\,\,
 \partial_{\mathbf n}V_\sigma+\sigma V_\sigma
       =M\lambda_1(\partial_{\mathbf n}\psi+K)\geq0.
\]
Subtracting its weak inequality from the equation for $u_\sigma$
and testing with $(u_\sigma-V_\sigma)_+\in H^1(\Omega)$ gives
\[
 \int_\Omega\bigl|\nabla(u_\sigma-V_\sigma)_+\bigr|^2\,\mathrm{d}x
 +\sigma\int_{\partial\Omega}(u_\sigma-V_\sigma)_+^2\,\mathrm{d}S
 \leq0.
\]
Both terms are nonnegative, so $(u_\sigma-V_\sigma)_+=0$ and
$u_\sigma\leq V_\sigma$. Since $\psi=0$ on $\partial\Omega$,
\eqref{eq:boundary-maximum-decay} follows with
$C_\Omega=M\lambda_1K$.

The first Robin eigenvalue converges to $\lambda_1$ as $\sigma\to+\infty$;
see \cite[Proposition~4.5]{Bucur}. Also,
\[
 \int_\Omega|\nabla u_\sigma|^2\,\mathrm{d}x
 +\sigma\int_{\partial\Omega}u_\sigma^2\,\mathrm{d}S
 =\mu_1(\Omega,\sigma)\leq\lambda_1.
\]
Together with the normalization, this bounds $u_\sigma$ in
$H^1(\Omega)$, while
$\|u_\sigma\|_{L^2(\partial\Omega)}^2\leq\lambda_1/\sigma\to0$.
Along any sequence $\sigma_j\to+\infty$, a subsequence therefore
converges weakly in $H^1$ and strongly in $L^2$ to a nonnegative
function in $H^1_0(\Omega)$ with $L^2$ norm one. Weak lower
semicontinuity gives a Dirichlet Rayleigh quotient at most
$\lambda_1$, hence equal to $\lambda_1$. Simplicity identifies the
limit as the positive normalized first Dirichlet eigenfunction
$v_1$, and gives convergence of the whole family in $L^2$.
The uniform bound and interior elliptic estimates upgrade this to
uniform convergence on compact subsets of $\Omega$.
For any $\varepsilon>0$, continuity and $\psi=v_1=0$ on
$\partial\Omega$ give a boundary neighborhood on which
$M\lambda_1\psi+v_1<\varepsilon$. For sufficiently large $\sigma$,
$M\lambda_1K/\sigma<\varepsilon$, so the comparison with $V_\sigma$
implies $|u_\sigma-v_1|\leq u_\sigma+v_1<2\varepsilon$ there.
On its compact complement, convergence is uniform. Hence
$u_\sigma\to v_1$ uniformly on $\overline\Omega$.

On each compact subset of $\Omega$, $v_1$ has a positive minimum,
whereas $\max_{\partial\Omega}u_\sigma\to0$. Thus each such compact
subset lies in $\{u_\sigma>\max_{\partial\Omega}u_\sigma\}$ for all
sufficiently large $\sigma$. Exhausting $\Omega$ by compact subsets
gives
\[
 \left|\left\{u_\sigma>\max_{\partial\Omega}u_\sigma\right\}\right|
 \longrightarrow|\Omega|,
\]
and hence $\widetilde R(\sigma)\to1$. The limit for $R(\sigma)$
follows from the eigenvalue convergence on $\Omega$ and on $B$.
If $\Omega$ is not a ball, strictness in the Dirichlet Faber--Krahn
inequality gives $j_{\nu,1}/\sqrt{\lambda_1(\Omega)}<1$; see
\cite{Chavel}. Thus $\widetilde R(\sigma)>R(\sigma)$ eventually.
On a unit ball both radii equal one for every $\sigma>0$.
\qed

\par\endgroup

\begingroup
\section{Chiti type comparison results}\label{sec:comparison}

The weight inequalities from Section~\ref{sec:weights} will be used
to compare integrals involving the first eigenfunction on $\Omega$
with those on a ball. For this purpose, the two eigenfunctions must
have the same first eigenvalue, and their decreasing rearrangements
must have a controlled relative ordering. Following the Robin
comparison of Dai and Shi \cite[Lemma~3.2 and Theorem~3.1]{DaiShi},
we establish that ordering on balls of radius at most $\widetilde R$.
The Robin parameter on the comparison ball may be chosen independently
of $\sigma$, as required by Theorem~\ref{thm:matched-bound}. Since
both eigenfunctions remain positive at their boundaries, the
comparison allows either no crossing or a single change of ordering.

Let $\Omega\subset\mathbb R^N$ be a bounded connected $C^2$ domain,
where $N\geq2$ and $|\Omega|=C_N$. Throughout this section, $\sigma>0$
is fixed, and we write $u=u_1$ and $\mu_1=\mu_1(\Omega,\sigma)$.
The first eigenfunction belongs to $C^1(\overline\Omega)$ and is positive
in $\Omega$. It is also strictly positive on $\partial\Omega$. Indeed,
if it vanished at a boundary point, the Robin condition would give
a zero normal derivative there, contrary to the Hopf lemma.
By the strong minimum principle, its minimum is attained on the
boundary. On the other hand, $\partial_{\texttt n}u=-\sigma u<0$
shows that its maximum is attained in the interior. Set
\[
 m=\min_{\overline\Omega}u,\,\,\,
 u_{1,M}=\max_{\overline\Omega}u.
\]
Then $0<m\leq u_{1,pM}<u_{1,M}$. For $t\geq0$, let
\[
 U_t=\{x\in\Omega:u(x)>t\},\,\,\, |U_t|=\mu(t).
\]
The decreasing rearrangement of $u$ is defined by
\[
 u^*(s)=\inf\{t\geq0:\mu(t)\leq s\},\,\,\, 0<s<|\Omega|,
\]
with endpoint values $u^*(0)=u_{1,M}$ and $u^*(|\Omega|)=m$.
Its spherically symmetric decreasing rearrangement is
$u^\star(x)=u^*(C_N|x|^N)$ on $\Omega^*=B$.
Since $u$ and its rearrangements have the same distribution function,
we have
\begin{equation*}
 \int_\Omega u^2\,\mathrm{d}x
 =\int_0^{|\Omega|}u^*(s)^2\,\mathrm{d}s
 =\int_{\Omega^*}(u^\star)^2\,\mathrm{d}x.
\end{equation*}
The function $u$ is nonconstant and real analytic in $\Omega$, so
each of its level sets has zero $N$-dimensional measure. Moreover,
since $\Omega$ is connected and $u$ is continuous, $\mu$ is strictly
decreasing between $m$ and $u_{1,M}$. It follows that $u^*$ is
continuous on $[0,|\Omega|]$. Put
\[
 s_\partial=\mu(u_{1,pM})=C_N\widetilde R^N>0,
 \,\,\, p(s)=N^2C_N^{2/N}s^{2-2/N}.
\]

For almost every $t\in(m,u_{1,M})$, the boundary of $U_t$ consists
of the interior part $\Gamma_t=\{u=t\}\cap\Omega$ and the exterior
part $E_t=\partial\Omega\cap\{u>t\}$, apart from sets of surface
measure zero. By the coarea formula and the equation for $u$, we have
\[
 -\mu'(t)=\int_{\Gamma_t}\frac1{|\nabla u|}\,\mathrm{d}S,
 \,\,\,
 \mu_1\int_{U_t}u\,\mathrm{d}x
 =\int_{\Gamma_t}|\nabla u|\,\mathrm{d}S
   +\sigma\int_{E_t}u\,\mathrm{d}S.
\]
Here $P(U_t)$ denotes the perimeter of $U_t$ in $\mathbb R^N$.
Applying the Cauchy--Schwarz inequality on the whole boundary of
$U_t$, with weight $|\nabla u|$ on $\Gamma_t$ and $\sigma u$ on
$E_t$, gives
\[
 P(U_t)^2\leq
 \mu_1\int_0^{\mu(t)}u^*(s)\,\mathrm{d}s
 \left(-\mu'(t)+\int_{E_t}\frac1{\sigma u}\,\mathrm{d}S\right).
\]
Combining this with the isoperimetric inequality, we obtain
\begin{equation}\label{Inequality 3.2}
 N^2C_N^{2/N}\mu(t)^{2-2/N}
 \leq\mu_1\int_0^{\mu(t)}u^*(s)\,\mathrm{d}s
 \left(-\mu'(t)+\int_{E_t}\frac1{\sigma u}\,\mathrm{d}S\right).
\end{equation}
For $t>u_{1,pM}$, the set $E_t$ is empty, giving the restricted
rearrangement inequality of \cite[Lemma~3.2]{DaiShi}.
\\ \\
\refstepcounter{lem}\label{lem:restricted-rearrangement}\textbf{Lemma~\thelem.}
\emph{The function $u^*$ is positive and continuous on $[0,|\Omega|]$.
It is locally absolutely continuous on $(0,s_\partial)$ and satisfies}
\begin{equation}\label{Inequality 3.3}
 -(u^*)'(s)\leq
 \mu_1N^{-2}C_N^{-2/N}s^{2/N-2}
 \int_0^s u^*(\zeta)\,\mathrm{d}\zeta
 \,\,\,\text{for a.e. }0<s<s_\partial.
\end{equation}

\noindent\textbf{Proof.}
For $t_0>u_{1,pM}$, the function $(u-t_0)_+$ has compact support
in $\Omega$ and belongs to $H^1_0(\Omega)$. Applying the Sobolev
rearrangement theorem to these functions, we obtain the local
absolute continuity of $u^*$ on $(0,s_\partial)$; see
\cite{Kawohl,Talenti}. For almost every regular level $t>u_{1,pM}$,
the boundary integral in (\ref{Inequality 3.2}) vanishes. Letting
$s=\mu(t)$ and using the fact that $u^*$ is the inverse of $\mu$
on this interval, we obtain (\ref{Inequality 3.3}).
Furthermore, since $u^*\leq u_{1,M}$, we have
\[
 0\leq-(u^*)'(s)\leq
 \mu_1N^{-2}C_N^{-2/N}u_{1,M} s^{2/N-1}.
\]
The right-hand side is integrable at zero. Together with the
continuity of $u^*$, this shows that $u^*$ is absolutely continuous
on $[0,S]$ for every $0<S\leq s_\partial$. The same estimate gives
$\int_0^S p(s)|(u^*)'(s)|^2\,\mathrm{d}s<+\infty$.
\qed
\\ \\
\refstepcounter{lem}\label{lem:comparison-radius}\textbf{Lemma~\thelem.}
\emph{The comparison radius $R$ defined in the introduction satisfies
$R\leq1$, or equivalently
$|B_R|\leq|\Omega|$, with equality if and only if $\Omega$ is a unit
ball, up to translation.}
\\ \\
\textbf{Proof.} By the Faber--Krahn inequality for the Robin Laplacian, we have
\[
 \mu_1(\Omega,\sigma)\geq\mu_1(B,\sigma)=\alpha^2,
\]
with equality if and only if $\Omega$ is a ball; see
\cite{Bossel,Daners,BucurDaners,DanersKennedy}.
The desired conclusion follows from the definition of $R$.
\qed\\

For the eigenfunction comparison, the rearrangement inequality is
paired with the radial equation on a ball having the same first
eigenvalue as $\Omega$. Let $0<r\leq\widetilde R$
and $S=|B_r|$, and let $z$ be a positive radial first eigenfunction
of
\[
 -\Delta z=\mu_1z\quad\text{in }B_r,\,\,\,
 z'(r)+\kappa z(r)=0,
\]
for some $\kappa>0$. Write $Z(s)=z^*(s)$ and $U(s)=u^*(s)$.
Integrating the equation for $z$ over concentric balls, we obtain
\begin{equation}\label{Equation 3.3}
 -p(s)Z'(s)=\mu_1\int_0^s Z(\zeta)\,\mathrm{d}\zeta,
 \,\,\, 0<s\leq S.
\end{equation}
In particular, $Z(S)>0$ and
\begin{equation}\label{Comparison Robin endpoint}
 \mu_1\int_0^S Z(\zeta)\,\mathrm{d}\zeta
 =\kappa|\partial B_r|Z(S).
\end{equation}
Since $r\leq\widetilde R$, Lemma~\ref{lem:restricted-rearrangement}
applies throughout $(0,S)$. For $r=R$, this condition is
$\widetilde R\geq R$. The next lemma excludes a return to equality
once $U$ has become larger than $Z$; no normalization of $z$ is
needed for this conclusion.
\\ \\
\refstepcounter{lem}\label{lem:no-return}\textbf{Lemma~\thelem.}
\emph{Assume $0<r\leq\widetilde R$ and
$\mu_1(B_r,\kappa)=\mu_1(\Omega,\sigma)$ with $\kappa>0$, and use
$U$, $Z$, and $S$ as above. If $U(a)>Z(a)$ for some $a\in[0,S)$,
then $U(s)>Z(s)$ for every $s\in(a,S]$.}
\\ \\
\textbf{Proof.} Suppose, for contradiction, that the conclusion fails. By continuity,
there is a point $b\in(a,S]$ such that $U(b)=Z(b)$. Define
\[
 v(s)=
 \begin{cases}
  \max\{U(s),Z(s)\},&0\leq s\leq b,\\
  Z(s),&b\leq s\leq S.
 \end{cases}
\]
The function $v$ is nonincreasing and absolutely continuous,
$v\geq Z$, $v(S)=Z(S)$, and $v\not\equiv Z$.
For $s\leq b$, we apply Lemma~\ref{lem:restricted-rearrangement} where $v=U$ and
(\ref{Equation 3.3}) where $v=Z$. In either case, the integral of
$v$ from zero to $s$ is no smaller than the corresponding integral
of $U$ or $Z$. For $s>b$, we use (\ref{Equation 3.3}) and $v\geq Z$.
It follows that
\[
 -p(s)v'(s)\leq\mu_1\int_0^s v(\zeta)\,\mathrm{d}\zeta
 \quad\text{for a.e. }0<s<S.
\]
The derivative estimate in Lemma~\ref{lem:restricted-rearrangement} and the corresponding identity
for $Z$ give $\int_0^S p(s)|v'(s)|^2\,\mathrm{d}s<+\infty$.
Thus the radial function $g(x)=v(C_N|x|^N)$ belongs to $H^1(B_r)$.
Put $A(s)=\int_0^s v(\zeta)\,\mathrm{d}\zeta$.
Multiplication by $-v'(s)\geq0$ and integration by parts use
$A(0)=0$ and $A'=v$ to give
\begin{align*}
 \int_{B_r}|\nabla g|^2\,\mathrm{d}x
 &=\int_0^S p(s)v'(s)^2\,\mathrm{d}s\\
 &\leq\mu_1\int_0^S v(s)^2\,\mathrm{d}s
       -\mu_1v(S)\int_0^S v(s)\,\mathrm{d}s\\
 &<\mu_1\int_{B_r}g^2\,\mathrm{d}x
       -\kappa|\partial B_r|Z(S)^2.
\end{align*}
For the strict inequality, continuity and $v\geq Z$, $v\not\equiv Z$
imply $\int_0^S v>\int_0^S Z$. Since $v(S)=Z(S)>0$, the endpoint
identity (\ref{Comparison Robin endpoint}) gives
\[
 \mu_1v(S)\int_0^S v(s)\,\mathrm{d}s
 >\mu_1 Z(S)\int_0^S Z(s)\,\mathrm{d}s
 =\kappa|\partial B_r|Z(S)^2.
\]
The boundary values of $g$ and $z$ agree, so the preceding energy
estimate makes the Robin Rayleigh quotient of $g$ strictly smaller
than $\mu_1(B_r,\kappa)=\mu_1$, a contradiction.
\qed\\

Figure~\ref{fig:excluded-return} illustrates the return excluded by
Lemma~\ref{lem:no-return}. The schematic comparison diagrams use radius
$R$; the same conclusions hold for an admissible radius $r$.
\begin{figure}[htbp]
\centering

\resizebox{0.48\textwidth}{!}{%
\begingroup%
\setlength{\unitlength}{1bp}%
\fontsize{12}{14.4}\selectfont%
\begin{picture}(260,180)%
\PPWxaxis{24}{44}{244}%
\PPWyaxis{60}{16}{165}%
\linethickness{1.4bp}
\qbezier(59.585,124.4141)(69.0257,121.4101)(77.8863,117.3678)
\qbezier(77.8863,117.3678)(86.7469,113.3256)(94.8988,108.2947)
\qbezier(94.8988,108.2947)(102.5674,103.5153)(109.5999,98.0901)
\qbezier(109.5999,98.0901)(116.6325,92.665)(123.2811,86.9034)
\qbezier(123.2811,86.9034)(126.8108,83.8534)(130.3416,80.7993)
\qbezier(130.3416,80.7993)(133.8724,77.7453)(137.5136,74.802)
\qbezier(137.5136,74.802)(141.1546,71.8588)(145.0158,69.141)
\qbezier(145.0158,69.141)(148.8768,66.4233)(153.067,64.046)
\qbezier(153.067,64.046)(158.6506,60.8886)(164.7361,58.533)
\qbezier(164.7361,58.533)(170.8217,56.1776)(177.2967,54.665)
\linethickness{1.4bp}
\put(60.2711,94.9273){\rule{1.7041bp}{1.4422bp}}
\put(64.0009,94.7247){\rule{1.7041bp}{1.4422bp}}
\put(67.7283,94.4909){\rule{1.7041bp}{1.4422bp}}
\put(71.4527,94.2259){\rule{1.7041bp}{1.4422bp}}
\put(75.1738,93.9292){\rule{1.7041bp}{1.4422bp}}
\put(78.8912,93.601){\rule{1.7041bp}{1.4422bp}}
\put(82.6044,93.2409){\rule{1.7041bp}{1.4422bp}}
\put(86.3131,92.849){\rule{1.7041bp}{1.4422bp}}
\put(90.0169,92.425){\rule{1.7041bp}{1.4422bp}}
\put(93.7154,91.969){\rule{1.7041bp}{1.4422bp}}
\put(97.4084,91.4821){\rule{1.7041bp}{1.4422bp}}
\put(101.0956,90.9642){\rule{1.7041bp}{1.4422bp}}
\put(104.776,90.4133){\rule{1.7041bp}{1.4422bp}}
\put(108.4491,89.8287){\rule{1.7041bp}{1.4422bp}}
\put(112.1143,89.2097){\rule{1.7041bp}{1.4422bp}}
\put(115.7711,88.5553){\rule{1.7041bp}{1.4422bp}}
\put(119.4186,87.8652){\rule{1.7041bp}{1.4422bp}}
\put(123.0562,87.1384){\rule{1.7041bp}{1.4422bp}}
\put(126.683,86.3743){\rule{1.7041bp}{1.4422bp}}
\put(130.2984,85.5722){\rule{1.7041bp}{1.4422bp}}
\put(133.9015,84.7314){\rule{1.7041bp}{1.4422bp}}
\put(137.4914,83.8513){\rule{1.7041bp}{1.4422bp}}
\put(141.0674,82.9312){\rule{1.7041bp}{1.4422bp}}
\put(144.6284,81.9705){\rule{1.7041bp}{1.4422bp}}
\put(148.1734,80.9682){\rule{1.7041bp}{1.4422bp}}
\put(151.7017,79.9249){\rule{1.7041bp}{1.4422bp}}
\put(155.2127,78.8405){\rule{1.7041bp}{1.4422bp}}
\put(158.7058,77.7154){\rule{1.7041bp}{1.4422bp}}
\put(162.1802,76.5494){\rule{1.7041bp}{1.4422bp}}
\put(165.6352,75.3429){\rule{1.7041bp}{1.4422bp}}
\put(169.0702,74.096){\rule{1.7041bp}{1.4422bp}}
\put(172.4843,72.8091){\rule{1.7041bp}{1.4422bp}}
\put(175.4315,71.6591){\rule{1.7041bp}{1.4422bp}}
\put(120.9795,86.3463){\rule{1.1662bp}{0.1576bp}}
\put(120.5731,86.4917){\rule{1.979bp}{0.1576bp}}
\put(120.3131,86.6371){\rule{2.499bp}{0.1576bp}}
\put(120.1195,86.7826){\rule{2.8864bp}{0.1576bp}}
\put(119.9675,86.928){\rule{3.1902bp}{0.1576bp}}
\put(119.8459,87.0734){\rule{3.4334bp}{0.1576bp}}
\put(119.7473,87.2189){\rule{3.6309bp}{0.1576bp}}
\put(119.671,87.3643){\rule{3.7832bp}{0.1576bp}}
\put(119.6118,87.5098){\rule{3.9015bp}{0.1576bp}}
\put(119.569,87.6552){\rule{3.9873bp}{0.1576bp}}
\put(119.5429,87.8006){\rule{4.0393bp}{0.1576bp}}
\put(119.531,87.9461){\rule{4.0631bp}{0.1576bp}}
\put(119.5331,88.0915){\rule{4.0592bp}{0.1576bp}}
\put(119.5495,88.2369){\rule{4.0261bp}{0.1576bp}}
\put(119.5822,88.3824){\rule{3.961bp}{0.1576bp}}
\put(119.6297,88.5278){\rule{3.8657bp}{0.1576bp}}
\put(119.694,88.6732){\rule{3.7373bp}{0.1576bp}}
\put(119.7785,88.8187){\rule{3.5683bp}{0.1576bp}}
\put(119.883,88.9641){\rule{3.3591bp}{0.1576bp}}
\put(120.0151,89.1096){\rule{3.0952bp}{0.1576bp}}
\put(120.18,89.255){\rule{2.7651bp}{0.1576bp}}
\put(120.3917,89.4004){\rule{2.3419bp}{0.1576bp}}
\put(120.684,89.5459){\rule{1.7575bp}{0.1576bp}}
\put(121.0869,89.6913){\rule{0.9514bp}{0.0967bp}}
\linethickness{1bp}
\qbezier(123.5961,88.0671)(123.5743,88.7982)(123.0004,89.2838)
\qbezier(123.0004,89.2838)(122.4266,89.7696)(121.5626,89.788)
\qbezier(121.5626,89.788)(120.6989,89.7696)(120.1249,89.2838)
\qbezier(120.1249,89.2838)(119.5511,88.7982)(119.5292,88.0671)
\qbezier(119.5292,88.0671)(119.5511,87.3361)(120.1249,86.8504)
\qbezier(120.1249,86.8504)(120.6989,86.3647)(121.5626,86.3463)
\qbezier(121.5626,86.3463)(122.4266,86.3647)(123.0004,86.8504)
\qbezier(123.0004,86.8504)(123.5743,87.3361)(123.5961,88.0671)
\put(59.3453,122.0994){\rule{1.1662bp}{0.1576bp}}
\put(58.939,122.2448){\rule{1.979bp}{0.1576bp}}
\put(58.679,122.3902){\rule{2.4992bp}{0.1576bp}}
\put(58.4854,122.5357){\rule{2.8865bp}{0.1576bp}}
\put(58.3334,122.6811){\rule{3.1904bp}{0.1576bp}}
\put(58.2118,122.8265){\rule{3.4335bp}{0.1576bp}}
\put(58.1132,122.972){\rule{3.6312bp}{0.1576bp}}
\put(58.0368,123.1174){\rule{3.7837bp}{0.1576bp}}
\put(57.9777,123.2629){\rule{3.902bp}{0.1576bp}}
\put(57.9349,123.4083){\rule{3.9877bp}{0.1576bp}}
\put(57.9088,123.5537){\rule{4.0397bp}{0.1576bp}}
\put(57.8969,123.6992){\rule{4.0635bp}{0.1576bp}}
\put(57.8989,123.8446){\rule{4.0596bp}{0.1576bp}}
\put(57.9154,123.99){\rule{4.0265bp}{0.1576bp}}
\put(57.9481,124.1355){\rule{3.9614bp}{0.1576bp}}
\put(57.9956,124.2809){\rule{3.8662bp}{0.1576bp}}
\put(58.0599,124.4263){\rule{3.7377bp}{0.1576bp}}
\put(58.1444,124.5718){\rule{3.5686bp}{0.1576bp}}
\put(58.2489,124.7172){\rule{3.3594bp}{0.1576bp}}
\put(58.381,124.8627){\rule{3.0954bp}{0.1576bp}}
\put(58.5459,125.0081){\rule{2.7655bp}{0.1576bp}}
\put(58.7574,125.1535){\rule{2.3422bp}{0.1576bp}}
\put(59.0494,125.299){\rule{1.7581bp}{0.1576bp}}
\put(59.4525,125.4444){\rule{0.952bp}{0.0968bp}}
\linethickness{1bp}
\qbezier(61.9624,123.8202)(61.9405,124.5513)(61.3664,125.0371)
\qbezier(61.3664,125.0371)(60.7923,125.5228)(59.9284,125.5412)
\qbezier(59.9284,125.5412)(59.0646,125.5228)(58.4908,125.0371)
\qbezier(58.4908,125.0371)(57.917,124.5513)(57.8951,123.8202)
\qbezier(57.8951,123.8202)(57.917,123.0892)(58.4908,122.6035)
\qbezier(58.4908,122.6035)(59.0646,122.1178)(59.9284,122.0994)
\qbezier(59.9284,122.0994)(60.7923,122.1178)(61.3664,122.6035)
\qbezier(61.3664,122.6035)(61.9405,123.0892)(61.9624,123.8202)
\put(59.0019,93.6787){\rule{1.1662bp}{0.1576bp}}
\put(58.5955,93.8241){\rule{1.979bp}{0.1576bp}}
\put(58.3356,93.9696){\rule{2.499bp}{0.1576bp}}
\put(58.1418,94.115){\rule{2.8862bp}{0.1576bp}}
\put(57.99,94.2604){\rule{3.19bp}{0.1576bp}}
\put(57.8684,94.4059){\rule{3.4331bp}{0.1576bp}}
\put(57.7696,94.5513){\rule{3.6306bp}{0.1576bp}}
\put(57.6935,94.6967){\rule{3.7831bp}{0.1576bp}}
\put(57.6343,94.8422){\rule{3.9012bp}{0.1576bp}}
\put(57.5915,94.9876){\rule{3.987bp}{0.1576bp}}
\put(57.5654,95.133){\rule{4.039bp}{0.1576bp}}
\put(57.5535,95.2785){\rule{4.0628bp}{0.1576bp}}
\put(57.5555,95.4239){\rule{4.0589bp}{0.1576bp}}
\put(57.572,95.5694){\rule{4.0258bp}{0.1576bp}}
\put(57.6047,95.7148){\rule{3.9607bp}{0.1576bp}}
\put(57.6522,95.8602){\rule{3.8654bp}{0.1576bp}}
\put(57.7165,96.0057){\rule{3.7371bp}{0.1576bp}}
\put(57.801,96.1511){\rule{3.568bp}{0.1576bp}}
\put(57.9055,96.2965){\rule{3.3589bp}{0.1576bp}}
\put(58.0376,96.442){\rule{3.095bp}{0.1576bp}}
\put(58.2025,96.5874){\rule{2.7651bp}{0.1576bp}}
\put(58.414,96.7328){\rule{2.3418bp}{0.1576bp}}
\put(58.7063,96.8783){\rule{1.7575bp}{0.1576bp}}
\put(59.1093,97.0237){\rule{0.9514bp}{0.0967bp}}
\linethickness{1bp}
\qbezier(61.6183,95.3996)(61.5966,96.1306)(61.0227,96.6164)
\qbezier(61.0227,96.6164)(60.4488,97.102)(59.585,97.1204)
\qbezier(59.585,97.1204)(58.7212,97.102)(58.1474,96.6164)
\qbezier(58.1474,96.6164)(57.5734,96.1306)(57.5517,95.3996)
\qbezier(57.5517,95.3996)(57.5734,94.6685)(58.1474,94.1829)
\qbezier(58.1474,94.1829)(58.7212,93.6971)(59.585,93.6787)
\qbezier(59.585,93.6787)(60.4488,93.6971)(61.0227,94.1829)
\qbezier(61.0227,94.1829)(61.5966,94.6685)(61.6183,95.3996)
\linethickness{0.6bp}
\qbezier(177.6404,72.178)(177.6404,69.1056)(177.6404,66.0333)
\qbezier(177.6404,61.8762)(177.6404,58.8039)(177.6404,55.7316)
\qbezier(177.6404,51.5745)(177.6404,48.5022)(177.6404,45.4298)
\linethickness{0.6bp}
\qbezier(121.5626,87.1217)(121.5626,84.0493)(121.5626,80.977)
\qbezier(121.5626,76.8199)(121.5626,73.7476)(121.5626,70.6753)
\qbezier(121.5626,66.5182)(121.5626,63.4459)(121.5626,60.3735)
\qbezier(121.5626,56.2165)(121.5626,53.1441)(121.5626,50.0718)
\qbezier(121.5626,45.9147)(121.5626,45.3452)(121.5626,44.7757)
\put(241,32){\makebox(0,0)[c]{$s$}}
\put(53,32){\makebox(0,0)[r]{$0$}}
\put(19.3308,122.6206){\makebox(0,0)[l]{\smash{$u_1^\ast(0)$}}}
\put(20.3622,94.2605){\makebox(0,0)[l]{\smash{$z^\ast(0)$}}}
\put(114.3596,28.8142){\makebox(0,0)[l]{\smash{$s_1$}}}
\put(162.4753,28.523){\makebox(0,0)[l]{\smash{$|B_R|$}}}

\end{picture}%
\endgroup%
}%

\caption{{A configuration excluded by Lemma~\ref{lem:no-return}.
The solid curve represents $u_1^*(s)$ and the dotted curve represents $z^*(s)$.}}
\label{fig:excluded-return}
\end{figure}

With the return to equality excluded, normalization determines the
possible orderings of the two rearrangements. Normalize $z$ so that
\begin{equation}\label{Equality 3.5}
 \int_{B_r}z^2\,\mathrm{d}x=\int_\Omega u_1^2\,\mathrm{d}x.
\end{equation}
Equivalently, $\int_0^S Z^2=\int_0^{|\Omega|}U^2$.
\\ \\
\refstepcounter{lem}\label{lem:single-crossing}\textbf{Lemma~\thelem.}
\emph{Keep $0<r\leq\widetilde R$ and
$\mu_1(B_r,\kappa)=\mu_1(\Omega,\sigma)$ with $\kappa>0$.
Assume also (\ref{Equality 3.5}). Then $Z(0)\geq U(0)$,
and exactly one of the following alternatives holds:}
\begin{itemize}
 \item[(i)] \emph{$Z(s)\geq U(s)$ for all $0\leq s\leq S$;}
 \item[(ii)] \emph{there is $s_1\in(0,S)$ such that
 $Z(s)\geq U(s)$ on $[0,s_1]$, $Z(s_1)=U(s_1)$, and
 $Z(s)<U(s)$ on $(s_1,S]$.}
\end{itemize}

\noindent\textbf{Proof.}
Since $S\leq s_\partial\leq|\Omega|$, if $U(0)>Z(0)$,
Lemma~\ref{lem:no-return} with $a=0$ would give $U>Z$ throughout
$[0,S]$. Consequently,
\[
 \int_0^S Z(s)^2\,\mathrm{d}s
 <\int_0^S U(s)^2\,\mathrm{d}s
 \leq\int_0^{|\Omega|}U(s)^2\,\mathrm{d}s
 =\int_0^S Z(s)^2\,\mathrm{d}s,
\]
a contradiction. Thus $Z(0)\geq U(0)$.

If $U\leq Z$ on $[0,S]$, alternative (i) holds. Otherwise, let
\[
 A=\{s\in[0,S]:U(s)>Z(s)\},\,\,\, s_1=\inf A.
\]
The nonempty set $A$ is relatively open in $[0,S]$, excludes $0$,
and contains a point less than $S$. For every $a\in A\cap[0,S)$,
Lemma~\ref{lem:no-return} gives $(a,S]\subset A$. Thus
$A=(s_1,S]$. In particular $s_1<S$,
$Z\geq U$ on $[0,s_1]$, and $Z<U$ on $(s_1,S]$.
Continuity gives $Z(s_1)=U(s_1)$.

If $s_1=0$, the same integral contradiction displayed above applies,
since $U>Z$ on $(0,S]$. Thus $s_1>0$ and alternative (ii) follows.
The normalization gives
\[
 \int_0^{s_1}(Z^2-U^2)\,\mathrm{d}s
 =\int_{s_1}^{S}(U^2-Z^2)\,\mathrm{d}s
   +\int_S^{|\Omega|}U^2\,\mathrm{d}s>0.
\]
Thus $Z>U$ on an interval before $s_1$, which is the unique change
of ordering but not necessarily the first contact.
\qed\\

Figure~\ref{fig:comparison-contacts} illustrates the contact
configurations excluded in the proof. Lemma~\ref{lem:no-return}
rules out panels (a)--(c), (e), and (f), in which $U$ returns to $Z$
after exceeding it. Normalization rules out panel (d), where the
curves agree before $s_1$.
\begin{figure}[htbp]
\centering
\subfloat[{A return to equality}]{

\resizebox{0.315\textwidth}{!}{%
\begingroup%
\setlength{\unitlength}{1bp}%
\fontsize{14.4}{17.28}\selectfont%
\renewcommand{\PPWaxislinewidth}{1.463}%
\renewcommand{\PPWarrowlength}{9.905}%
\renewcommand{\PPWarrowhalfwidth}{3.81}%
\renewcommand{\PPWarrowinset}{7.924}%
%
\endgroup%
}%
}
\caption{{The relative positions of $u_1^*(s)$ (solid)
and $z^*(s)$ (dotted) in the comparison argument.}}
\label{fig:comparison-contacts}
\end{figure}
\\ \refstepcounter{lem}\label{lem:matched-comparison}\textbf{Lemma~\thelem.}
\emph{Let $\tau>0$ and $0<r\leq\widetilde R$ satisfy}
\[
 \frac{\mu_1(B,\tau)}{r^2}=\mu_1(\Omega,\sigma).
\]
\emph{Let $z$ be the positive first eigenfunction on $B_r$ with Robin
parameter $\tau/r$, normalized by (\ref{Equality 3.5}).
Then $u_1^*$ and $z^*$ satisfy one of the two alternatives in
Lemma~\ref{lem:single-crossing}.}
\\ \\
\noindent\textbf{Proof.}
By scaling, $\mu_1(B_r,\tau/r)=\mu_1(B,\tau)/r^2=\mu_1$.
Since $|B_r|\leq\mu(u_{1,pM})$, the inequality in Lemma~\ref{lem:restricted-rearrangement} holds
throughout the interval required for the comparison.
Applying Lemmas~\ref{lem:no-return} and~\ref{lem:single-crossing} with $\kappa=\tau/r$, we obtain the
desired conclusion.
\qed\\

Taking $r=R$ and $\tau=\sigma$ in Lemma~\ref{lem:matched-comparison}
recovers the comparison in \cite[Theorem~3.1]{DaiShi} with $p=2$,
$R^*=1$, $R_\lambda=R$ and $R_M=\widetilde R$ in their notation.
\\ \\
\refstepcounter{thm}\label{thm:chiti}\textbf{Theorem~\thethm.}
\emph{Assume $\widetilde R\geq R$, and let $z$ be the first
eigenfunction on $B_R$ with Robin parameter $\sigma/R$, normalized
by (\ref{Equality 3.5}) with $r=R$. Then either}
\[
 z^*(s)\geq u_1^*(s)\,\,\,(0\leq s\leq |B_R|),
\]
\emph{or there is $s_1\in(0,|B_R|)$ such that}
\[
 z^*(s)\geq u_1^*(s)\quad(0\leq s\leq s_1),\,\,\,
 z^*(s)<u_1^*(s)\quad(s_1<s\leq |B_R|),
\]
\emph{with equality at $s_1$.}
\\ \\
\textbf{Proof.} By \eqref{eq:robin-scaling} and the definition of $R$,
\[
 \mu_1(B_R,\sigma/R)
 =R^{-2}\mu_1(B,\sigma)
 =\frac{\alpha^2}{R^2}
 =\mu_1(\Omega,\sigma).
\]
Since $R\leq\widetilde R$, Lemma~\ref{lem:single-crossing} applies
with $r=R$, $\kappa=\sigma/R>0$, and normalization
(\ref{Equality 3.5}), proving both alternatives.
\qed\\

The two alternatives in Theorem~\ref{thm:chiti} are shown in
Figure~\ref{fig:chiti-comparison}.
\begin{figure}[htbp]
\centering
\subfloat[{$z^*\geq u_1^*$ throughout}]{

\resizebox{0.48\textwidth}{!}{%
\begingroup%
\setlength{\unitlength}{1bp}%
\fontsize{12}{14.4}\selectfont%
\begin{picture}(260,180)%
\PPWxaxis{24}{44}{244}%
\PPWyaxis{60}{16}{165}%
\put(59.1159,121.8934){\rule{1.1662bp}{0.1576bp}}
\put(58.7095,122.0389){\rule{1.9792bp}{0.1576bp}}
\put(58.4494,122.1843){\rule{2.4992bp}{0.1576bp}}
\put(58.2558,122.3298){\rule{2.8864bp}{0.1576bp}}
\put(58.1039,122.4752){\rule{3.1901bp}{0.1576bp}}
\put(57.9824,122.6206){\rule{3.4333bp}{0.1576bp}}
\put(57.8836,122.7661){\rule{3.6307bp}{0.1576bp}}
\put(57.8074,122.9115){\rule{3.7831bp}{0.1576bp}}
\put(57.7483,123.0569){\rule{3.9014bp}{0.1576bp}}
\put(57.7055,123.2024){\rule{3.987bp}{0.1576bp}}
\put(57.6794,123.3478){\rule{4.039bp}{0.1576bp}}
\put(57.6675,123.4932){\rule{4.0628bp}{0.1576bp}}
\put(57.6695,123.6387){\rule{4.0589bp}{0.1576bp}}
\put(57.6861,123.7841){\rule{4.0258bp}{0.1576bp}}
\put(57.7186,123.9296){\rule{3.9607bp}{0.1576bp}}
\put(57.7663,124.075){\rule{3.8655bp}{0.1576bp}}
\put(57.8305,124.2204){\rule{3.737bp}{0.1576bp}}
\put(57.915,124.3659){\rule{3.5679bp}{0.1576bp}}
\put(58.0197,124.5113){\rule{3.3588bp}{0.1576bp}}
\put(58.1515,124.6567){\rule{3.0949bp}{0.1576bp}}
\put(58.3165,124.8022){\rule{2.7649bp}{0.1576bp}}
\put(58.5282,124.9476){\rule{2.3416bp}{0.1576bp}}
\put(58.8203,125.093){\rule{1.7574bp}{0.1576bp}}
\put(59.2233,125.2385){\rule{0.9513bp}{0.0967bp}}
\linethickness{1bp}
\qbezier(61.7323,123.6142)(61.7105,124.3454)(61.1367,124.831)
\qbezier(61.1367,124.831)(60.5628,125.3168)(59.699,125.3352)
\qbezier(59.699,125.3352)(58.8352,125.3168)(58.2614,124.831)
\qbezier(58.2614,124.831)(57.6874,124.3454)(57.6657,123.6142)
\qbezier(57.6657,123.6142)(57.6874,122.8831)(58.2614,122.3975)
\qbezier(58.2614,122.3975)(58.8352,121.9119)(59.699,121.8934)
\qbezier(59.699,121.8934)(60.5628,121.9119)(61.1367,122.3975)
\qbezier(61.1367,122.3975)(61.7105,122.8831)(61.7323,123.6142)
\linethickness{1.4bp}
\qbezier(60.5867,97.6293)(85.7439,96.1358)(110.6643,92.8421)
\qbezier(110.6643,92.8421)(115.6454,92.1893)(120.6064,91.4242)
\qbezier(120.6064,91.4242)(125.5672,90.6589)(130.473,89.7063)
\qbezier(130.473,89.7063)(135.3788,88.7537)(140.1946,87.5385)
\qbezier(140.1946,87.5385)(145.0105,86.3234)(149.7015,84.7707)
\qbezier(149.7015,84.7707)(156.2425,82.5947)(162.3632,79.7511)
\qbezier(162.3632,79.7511)(168.484,76.9076)(174.1174,73.4386)
\qbezier(174.1174,73.4386)(179.7508,69.9695)(184.8298,65.9171)
\qbezier(184.8298,65.9171)(189.9087,61.8646)(194.3662,57.2708)
\linethickness{1.4bp}%
\put(58.8469,122.8931){\rule{1.7041bp}{1.4422bp}}%
\put(62.4416,122.6975){\rule{1.7041bp}{1.4422bp}}%
\put(66.0352,122.4819){\rule{1.7041bp}{1.4422bp}}%
\put(69.6274,122.2462){\rule{1.7041bp}{1.4422bp}}%
\put(73.2183,121.9906){\rule{1.7041bp}{1.4422bp}}%
\put(76.8078,121.7151){\rule{1.7041bp}{1.4422bp}}%
\put(80.3956,121.4196){\rule{1.7041bp}{1.4422bp}}%
\put(83.9818,121.1042){\rule{1.7041bp}{1.4422bp}}%
\put(87.5661,120.7691){\rule{1.7041bp}{1.4422bp}}%
\put(91.1486,120.4141){\rule{1.7041bp}{1.4422bp}}%
\put(94.7290,120.0394){\rule{1.7041bp}{1.4422bp}}%
\put(98.3074,119.6449){\rule{1.7041bp}{1.4422bp}}%
\put(101.8835,119.2308){\rule{1.7041bp}{1.4422bp}}%
\put(105.4572,118.7971){\rule{1.7041bp}{1.4422bp}}%
\put(109.0286,118.3439){\rule{1.7041bp}{1.4422bp}}%
\put(112.5974,117.8711){\rule{1.7041bp}{1.4422bp}}%
\put(116.1636,117.3789){\rule{1.7041bp}{1.4422bp}}%
\put(119.7270,116.8673){\rule{1.7041bp}{1.4422bp}}%
\put(123.2877,116.3364){\rule{1.7041bp}{1.4422bp}}%
\put(126.8454,115.7863){\rule{1.7041bp}{1.4422bp}}%
\put(130.4001,115.2169){\rule{1.7041bp}{1.4422bp}}%
\put(133.9516,114.6284){\rule{1.7041bp}{1.4422bp}}%
\put(137.5000,114.0208){\rule{1.7041bp}{1.4422bp}}%
\put(141.0450,113.3942){\rule{1.7041bp}{1.4422bp}}%
\put(144.5867,112.7487){\rule{1.7041bp}{1.4422bp}}%
\put(148.1249,112.0844){\rule{1.7041bp}{1.4422bp}}%
\put(151.6594,111.4012){\rule{1.7041bp}{1.4422bp}}%
\put(155.1904,110.6994){\rule{1.7041bp}{1.4422bp}}%
\put(158.7175,109.9790){\rule{1.7041bp}{1.4422bp}}%
\put(162.2409,109.2400){\rule{1.7041bp}{1.4422bp}}%
\put(165.7603,108.4825){\rule{1.7041bp}{1.4422bp}}%
\put(169.2757,107.7067){\rule{1.7041bp}{1.4422bp}}%
\put(172.7870,106.9126){\rule{1.7041bp}{1.4422bp}}%
\put(176.2942,106.1003){\rule{1.7041bp}{1.4422bp}}%
\put(179.7971,105.2698){\rule{1.7041bp}{1.4422bp}}%
\put(183.2956,104.4214){\rule{1.7041bp}{1.4422bp}}%
\put(186.7898,103.5550){\rule{1.7041bp}{1.4422bp}}%
\put(190.2795,102.6707){\rule{1.7041bp}{1.4422bp}}%
\put(193.7647,101.7687){\rule{1.7041bp}{1.4422bp}}%
\linethickness{0.6bp}
\qbezier(195.3112,102.3078)(195.3096,99.2354)(195.3082,96.1631)
\qbezier(195.3062,92.006)(195.3046,88.9337)(195.3031,85.8614)
\qbezier(195.301,81.7043)(195.2996,78.632)(195.298,75.5596)
\qbezier(195.296,71.4026)(195.2946,68.3302)(195.293,65.2579)
\qbezier(195.291,61.1008)(195.2895,58.0285)(195.288,54.9561)
\qbezier(195.286,50.7991)(195.2844,47.7267)(195.283,44.6544)
\put(58.815,95.9085){\rule{1.1661bp}{0.1576bp}}
\put(58.4086,96.0539){\rule{1.9791bp}{0.1576bp}}
\put(58.1485,96.1993){\rule{2.499bp}{0.1576bp}}
\put(57.9549,96.3448){\rule{2.8864bp}{0.1576bp}}
\put(57.803,96.4902){\rule{3.1903bp}{0.1576bp}}
\put(57.6814,96.6356){\rule{3.4336bp}{0.1576bp}}
\put(57.5827,96.7811){\rule{3.631bp}{0.1576bp}}
\put(57.5064,96.9265){\rule{3.7835bp}{0.1576bp}}
\put(57.4474,97.072){\rule{3.902bp}{0.1576bp}}
\put(57.4045,97.2174){\rule{3.9877bp}{0.1576bp}}
\put(57.3784,97.3628){\rule{4.0397bp}{0.1576bp}}
\put(57.3665,97.5083){\rule{4.0635bp}{0.1576bp}}
\put(57.3685,97.6537){\rule{4.0596bp}{0.1576bp}}
\put(57.385,97.7991){\rule{4.0266bp}{0.1576bp}}
\put(57.4176,97.9446){\rule{3.9614bp}{0.1576bp}}
\put(57.4652,98.09){\rule{3.8662bp}{0.1576bp}}
\put(57.5293,98.2354){\rule{3.7377bp}{0.1576bp}}
\put(57.6138,98.3809){\rule{3.5687bp}{0.1576bp}}
\put(57.7184,98.5263){\rule{3.3597bp}{0.1576bp}}
\put(57.8503,98.6718){\rule{3.0959bp}{0.1576bp}}
\put(58.0151,98.8172){\rule{2.7661bp}{0.1576bp}}
\put(58.2264,98.9626){\rule{2.3432bp}{0.1576bp}}
\put(58.518,99.1081){\rule{1.7601bp}{0.1576bp}}
\put(58.921,99.2535){\rule{0.954bp}{0.0973bp}}
\linethickness{1bp}
\qbezier(61.432,97.6293)(61.4101,98.3605)(60.836,98.8464)
\qbezier(60.836,98.8464)(60.2619,99.3323)(59.398,99.3508)
\qbezier(59.398,99.3508)(58.5342,99.3323)(57.9604,98.8464)
\qbezier(57.9604,98.8464)(57.3866,98.3605)(57.3647,97.6293)
\qbezier(57.3647,97.6293)(57.3866,96.8983)(57.9604,96.4126)
\qbezier(57.9604,96.4126)(58.5342,95.927)(59.398,95.9085)
\qbezier(59.398,95.9085)(60.2619,95.927)(60.836,96.4126)
\qbezier(60.836,96.4126)(61.4101,96.8983)(61.432,97.6293)
\put(241,32){\makebox(0,0)[c]{$s$}}
\put(53,32){\makebox(0,0)[r]{$0$}}
\put(200.5525,101.9684){\makebox(0,0)[l]{\smash{$z^\ast(s)$}}}
\put(201.7554,55.8654){\makebox(0,0)[l]{\smash{$u_1^\ast(s)$}}}
\put(180.9625,28.8143){\makebox(0,0)[l]{\smash{$|B_R|$}}}

\end{picture}%
\endgroup%
}%
}\hfill
\subfloat[{One change of sign}]{

\resizebox{0.48\textwidth}{!}{%
\begingroup%
\setlength{\unitlength}{1bp}%
\fontsize{12}{14.4}\selectfont%
\begin{picture}(260,180)%
\PPWxaxis{24}{44}{244}%
\PPWyaxis{60}{16}{165}%
\put(59.1164,121.8934){\rule{1.1662bp}{0.1576bp}}
\put(58.71,122.0389){\rule{1.9792bp}{0.1576bp}}
\put(58.4498,122.1843){\rule{2.4993bp}{0.1576bp}}
\put(58.2562,122.3298){\rule{2.8867bp}{0.1576bp}}
\put(58.1042,122.4752){\rule{3.1905bp}{0.1576bp}}
\put(57.9826,122.6206){\rule{3.4337bp}{0.1576bp}}
\put(57.8838,122.7661){\rule{3.6313bp}{0.1576bp}}
\put(57.8076,122.9115){\rule{3.7837bp}{0.1576bp}}
\put(57.7483,123.0569){\rule{3.902bp}{0.1576bp}}
\put(57.7055,123.2024){\rule{3.9877bp}{0.1576bp}}
\put(57.6794,123.3478){\rule{4.0397bp}{0.1576bp}}
\put(57.6675,123.4932){\rule{4.0635bp}{0.1576bp}}
\put(57.6695,123.6387){\rule{4.0596bp}{0.1576bp}}
\put(57.6861,123.7841){\rule{4.0265bp}{0.1576bp}}
\put(57.7188,123.9296){\rule{3.9612bp}{0.1576bp}}
\put(57.7663,124.075){\rule{3.866bp}{0.1576bp}}
\put(57.8306,124.2204){\rule{3.7376bp}{0.1576bp}}
\put(57.9151,124.3659){\rule{3.5684bp}{0.1576bp}}
\put(58.0198,124.5113){\rule{3.3592bp}{0.1576bp}}
\put(58.1518,124.6567){\rule{3.0952bp}{0.1576bp}}
\put(58.317,124.8022){\rule{2.7651bp}{0.1576bp}}
\put(58.5286,124.9476){\rule{2.3418bp}{0.1576bp}}
\put(58.8209,125.093){\rule{1.7574bp}{0.1576bp}}
\put(59.224,125.2385){\rule{0.9513bp}{0.0967bp}}
\linethickness{1bp}
\qbezier(61.733,123.6142)(61.7111,124.3454)(61.1373,124.831)
\qbezier(61.1373,124.831)(60.5635,125.3168)(59.6997,125.3352)
\qbezier(59.6997,125.3352)(58.8358,125.3168)(58.2617,124.831)
\qbezier(58.2617,124.831)(57.6876,124.3454)(57.6657,123.6142)
\qbezier(57.6657,123.6142)(57.6876,122.8831)(58.2617,122.3975)
\qbezier(58.2617,122.3975)(58.8358,121.9119)(59.6997,121.8934)
\qbezier(59.6997,121.8934)(60.5635,121.9119)(61.1373,122.3975)
\qbezier(61.1373,122.3975)(61.7111,122.8831)(61.733,123.6142)
\linethickness{1.4bp}
\qbezier(60.5874,99.3991)(73.6374,99.9345)(86.6392,99.1111)
\qbezier(86.6392,99.1111)(99.6411,98.2878)(112.4549,96.1144)
\qbezier(112.4549,96.1144)(122.3686,94.4222)(132.0539,91.9939)
\qbezier(132.0539,91.9939)(141.7392,89.5656)(151.191,86.5281)
\qbezier(151.191,86.5281)(162.4721,82.9257)(173.3639,78.5669)
\qbezier(173.3639,78.5669)(184.2558,74.208)(194.71,69.1239)
\linethickness{1.4bp}
\put(59.6748,122.9744){\rule{1.7041bp}{1.4422bp}}
\put(63.2889,122.1687){\rule{1.7041bp}{1.4422bp}}
\put(66.8787,121.2882){\rule{1.7041bp}{1.4422bp}}
\put(70.4417,120.3332){\rule{1.7041bp}{1.4422bp}}
\put(73.9758,119.3038){\rule{1.7041bp}{1.4422bp}}
\put(77.4786,118.2007){\rule{1.7041bp}{1.4422bp}}
\put(80.9478,117.024){\rule{1.7041bp}{1.4422bp}}
\put(84.3814,115.7746){\rule{1.7041bp}{1.4422bp}}
\put(87.7745,114.4484){\rule{1.7041bp}{1.4422bp}}
\put(91.1219,113.0413){\rule{1.7041bp}{1.4422bp}}
\put(94.4171,111.549){\rule{1.7041bp}{1.4422bp}}
\put(97.6537,109.9673){\rule{1.7041bp}{1.4422bp}}
\put(100.8239,108.2922){\rule{1.7041bp}{1.4422bp}}
\put(103.9197,106.5201){\rule{1.7041bp}{1.4422bp}}
\put(106.9319,104.6477){\rule{1.7041bp}{1.4422bp}}
\put(109.8511,102.6727){\rule{1.7041bp}{1.4422bp}}
\put(112.667,100.5932){\rule{1.7041bp}{1.4422bp}}
\put(115.3695,98.4084){\rule{1.7041bp}{1.4422bp}}
\put(117.9479,96.1187){\rule{1.7041bp}{1.4422bp}}
\put(120.3967,93.7295){\rule{1.7041bp}{1.4422bp}}
\put(122.7407,91.2658){\rule{1.7041bp}{1.4422bp}}
\put(125.023,88.7609){\rule{1.7041bp}{1.4422bp}}
\put(127.2903,86.246){\rule{1.7041bp}{1.4422bp}}
\put(129.5894,83.7522){\rule{1.7041bp}{1.4422bp}}
\put(131.9662,81.3111){\rule{1.7041bp}{1.4422bp}}
\put(134.4626,78.9574){\rule{1.7041bp}{1.4422bp}}
\put(137.0903,76.7083){\rule{1.7041bp}{1.4422bp}}
\put(139.8361,74.5624){\rule{1.7041bp}{1.4422bp}}
\put(142.6864,72.5166){\rule{1.7041bp}{1.4422bp}}
\put(145.629,70.5667){\rule{1.7041bp}{1.4422bp}}
\put(148.6524,68.7073){\rule{1.7041bp}{1.4422bp}}
\put(151.7465,66.9329){\rule{1.7041bp}{1.4422bp}}
\put(154.9072,65.2449){\rule{1.7041bp}{1.4422bp}}
\put(158.1272,63.6391){\rule{1.7041bp}{1.4422bp}}
\put(161.4033,62.1167){\rule{1.7041bp}{1.4422bp}}
\put(164.7324,60.6793){\rule{1.7041bp}{1.4422bp}}
\put(168.1117,59.3281){\rule{1.7041bp}{1.4422bp}}
\put(171.5379,58.0642){\rule{1.7041bp}{1.4422bp}}
\put(175.0079,56.8891){\rule{1.7041bp}{1.4422bp}}
\put(178.5183,55.8038){\rule{1.7041bp}{1.4422bp}}
\put(182.0663,54.8091){\rule{1.7041bp}{1.4422bp}}
\put(185.648,53.9062){\rule{1.7041bp}{1.4422bp}}
\put(189.2607,53.0956){\rule{1.7041bp}{1.4422bp}}
\put(192.8124,52.3946){\rule{1.7041bp}{1.4422bp}}
\linethickness{0.6bp}
\qbezier(195.311,67.7543)(195.3073,64.682)(195.3037,61.6096)
\qbezier(195.2987,57.4526)(195.2951,54.3802)(195.2914,51.3079)
\qbezier(195.2864,47.1508)(195.2848,45.8177)(195.2832,44.4845)
\put(58.5146,97.6782){\rule{1.1662bp}{0.1576bp}}
\put(58.1083,97.8236){\rule{1.979bp}{0.1576bp}}
\put(57.8483,97.9691){\rule{2.4992bp}{0.1576bp}}
\put(57.6546,98.1145){\rule{2.8865bp}{0.1576bp}}
\put(57.5027,98.2599){\rule{3.1904bp}{0.1576bp}}
\put(57.3811,98.4054){\rule{3.4335bp}{0.1576bp}}
\put(57.2825,98.5508){\rule{3.6312bp}{0.1576bp}}
\put(57.2061,98.6962){\rule{3.7837bp}{0.1576bp}}
\put(57.147,98.8417){\rule{3.902bp}{0.1576bp}}
\put(57.1042,98.9871){\rule{3.9877bp}{0.1576bp}}
\put(57.0781,99.1325){\rule{4.0397bp}{0.1576bp}}
\put(57.0662,99.278){\rule{4.0635bp}{0.1576bp}}
\put(57.0682,99.4234){\rule{4.0596bp}{0.1576bp}}
\put(57.0847,99.5689){\rule{4.0265bp}{0.1576bp}}
\put(57.1173,99.7143){\rule{3.9614bp}{0.1576bp}}
\put(57.1649,99.8597){\rule{3.8662bp}{0.1576bp}}
\put(57.2292,100.0052){\rule{3.7377bp}{0.1576bp}}
\put(57.3137,100.1506){\rule{3.5686bp}{0.1576bp}}
\put(57.4182,100.296){\rule{3.3594bp}{0.1576bp}}
\put(57.5502,100.4415){\rule{3.0953bp}{0.1576bp}}
\put(57.7152,100.5869){\rule{2.7652bp}{0.1576bp}}
\put(57.9269,100.7323){\rule{2.3419bp}{0.1576bp}}
\put(58.219,100.8778){\rule{1.7575bp}{0.1576bp}}
\put(58.622,101.0232){\rule{0.9514bp}{0.0967bp}}
\linethickness{1bp}
\qbezier(61.1317,99.3991)(61.1098,100.1301)(60.5357,100.6159)
\qbezier(60.5357,100.6159)(59.9616,101.1015)(59.0977,101.1199)
\qbezier(59.0977,101.1199)(58.2339,101.1015)(57.6601,100.6159)
\qbezier(57.6601,100.6159)(57.0863,100.1301)(57.0644,99.3991)
\qbezier(57.0644,99.3991)(57.0863,98.668)(57.6601,98.1822)
\qbezier(57.6601,98.1822)(58.2339,97.6966)(59.0977,97.6782)
\qbezier(59.0977,97.6782)(59.9616,97.6966)(60.5357,98.1822)
\qbezier(60.5357,98.1822)(61.1098,98.668)(61.1317,99.3991)
\linethickness{0.6bp}
\qbezier(121.6913,94.9995)(121.6896,91.9271)(121.6879,88.8548)
\qbezier(121.6856,84.6977)(121.6838,81.6254)(121.6821,78.5531)
\qbezier(121.6798,74.396)(121.6781,71.3237)(121.6764,68.2513)
\qbezier(121.6741,64.0943)(121.6724,61.0219)(121.6707,57.9496)
\qbezier(121.6684,53.7925)(121.6667,50.7202)(121.6651,47.6478)
\put(121.1083,92.6244){\rule{1.1662bp}{0.1576bp}}
\put(120.7019,92.7698){\rule{1.979bp}{0.1576bp}}
\put(120.442,92.9153){\rule{2.499bp}{0.1576bp}}
\put(120.2482,93.0607){\rule{2.8862bp}{0.1576bp}}
\put(120.0965,93.2061){\rule{3.19bp}{0.1576bp}}
\put(119.9749,93.3516){\rule{3.4331bp}{0.1576bp}}
\put(119.8761,93.497){\rule{3.6306bp}{0.1576bp}}
\put(119.7999,93.6424){\rule{3.7831bp}{0.1576bp}}
\put(119.7407,93.7879){\rule{3.9012bp}{0.1576bp}}
\put(119.6979,93.9333){\rule{3.987bp}{0.1576bp}}
\put(119.6719,94.0788){\rule{4.039bp}{0.1576bp}}
\put(119.66,94.2242){\rule{4.0628bp}{0.1576bp}}
\put(119.662,94.3696){\rule{4.0589bp}{0.1576bp}}
\put(119.6784,94.5151){\rule{4.0258bp}{0.1576bp}}
\put(119.7111,94.6605){\rule{3.9607bp}{0.1576bp}}
\put(119.7586,94.8059){\rule{3.8654bp}{0.1576bp}}
\put(119.8229,94.9514){\rule{3.7371bp}{0.1576bp}}
\put(119.9074,95.0968){\rule{3.568bp}{0.1576bp}}
\put(120.012,95.2422){\rule{3.3589bp}{0.1576bp}}
\put(120.144,95.3877){\rule{3.095bp}{0.1576bp}}
\put(120.309,95.5331){\rule{2.7651bp}{0.1576bp}}
\put(120.5205,95.6786){\rule{2.3418bp}{0.1576bp}}
\put(120.8127,95.824){\rule{1.7575bp}{0.1576bp}}
\put(121.2157,95.9694){\rule{0.9514bp}{0.0967bp}}
\linethickness{1bp}
\qbezier(123.7248,94.3453)(123.703,95.0763)(123.1292,95.562)
\qbezier(123.1292,95.562)(122.5552,96.0477)(121.6914,96.0661)
\qbezier(121.6914,96.0661)(120.8276,96.0477)(120.2537,95.562)
\qbezier(120.2537,95.562)(119.6799,95.0763)(119.6581,94.3453)
\qbezier(119.6581,94.3453)(119.6799,93.6142)(120.2537,93.1286)
\qbezier(120.2537,93.1286)(120.8276,92.6428)(121.6914,92.6244)
\qbezier(121.6914,92.6244)(122.5552,92.6428)(123.1292,93.1286)
\qbezier(123.1292,93.1286)(123.703,93.6142)(123.7248,94.3453)
\put(241,32){\makebox(0,0)[c]{$s$}}
\put(53,32){\makebox(0,0)[r]{$0$}}
\put(200.5526,51.5022){\makebox(0,0)[l]{\smash{$z^\ast(s)$}}}
\put(200.8967,70.6995){\makebox(0,0)[l]{\smash{$u_1^\ast(s)$}}}
\put(180.7909,28.8141){\makebox(0,0)[l]{\smash{$|B_R|$}}}
\put(111.5389,28.8141){\makebox(0,0)[l]{\smash{$s_1$}}}

\end{picture}%
\endgroup%
}%
}
\caption{{The two alternatives in Theorem~\ref{thm:chiti}.}}
\label{fig:chiti-comparison}
\end{figure}
\\

\noindent\refstepcounter{rek}\label{rem:strict-mass}\textbf{Remark~\therek.}
In alternative (ii), it follows from (\ref{Equality 3.5}) that
\[
 \int_0^{s_1}\bigl(Z(s)^2-U(s)^2\bigr)\,\mathrm{d}s>0.
\]
In alternative (i), if $S<|\Omega|$, then
\[
 \int_0^S\bigl(Z(s)^2-U(s)^2\bigr)\,\mathrm{d}s
 =\int_S^{|\Omega|}U(s)^2\,\mathrm{d}s>0.
\]
\\
\refstepcounter{rek}\label{rem:weighted-comparison}\textbf{Remark~\therek.}
Extend $z$ by zero outside $B_r$. In both alternatives,
$z^2-(u_1^\star)^2$ is almost everywhere nonnegative up to a radius
and nonpositive afterwards, while its integral over $\Omega^*$ is zero.
\par\endgroup

\section{A common comparison estimate}
\label{sec:common-comparison}

Let $\Omega\subset\mathbb R^N$, $N\geq2$, be a bounded connected $C^2$
domain with $|\Omega|=|B|$, fix $\sigma>0$, and write
$\mu_1=\mu_1(\Omega,\sigma)$. The comparison below allows the Robin
parameter on the ball to vary independently of $\sigma$. Its two
choices in Section~\ref{sec:matched-bound} give the bounds in
Theorems~\ref{thm:ball-bound} and~\ref{thm:matched-bound}.

The Rayleigh principle gives
\[
 \mu_2(\Omega,\sigma)
 =\inf_{\substack{\varphi\in H^1(\Omega)\setminus\{0\}\\
       \int_\Omega\varphi u_1\,\mathrm{d}x=0}}
 \frac{\displaystyle\int_\Omega|\nabla\varphi|^2\,\mathrm{d}x
       +\sigma\int_{\partial\Omega}\varphi^2\,\mathrm{d}S}
      {\displaystyle\int_\Omega\varphi^2\,\mathrm{d}x}.
\]
Since $u_1\in C^1(\overline\Omega)$ and
$\min_{\overline\Omega}u_1>0$, multiplication and division by $u_1$
are bounded maps on $H^1(\Omega)$. Put $\varphi=u_1P$. For
$P\in C^\infty(\overline\Omega)$, testing the equation for $u_1$
with $u_1P^2$ yields
\[
 2\int_\Omega u_1P\nabla u_1\cdot\nabla P\,\mathrm{d}x
 +\int_\Omega P^2|\nabla u_1|^2\,\mathrm{d}x
 +\sigma\int_{\partial\Omega}u_1^2P^2\,\mathrm{d}S
 =\mu_1\int_\Omega u_1^2P^2\,\mathrm{d}x.
\]
Expanding $|\nabla(u_1P)|^2$ therefore gives
\begin{equation}\label{eq:ground-state-identity}
 \int_\Omega|\nabla(u_1P)|^2\,\mathrm{d}x
 +\sigma\int_{\partial\Omega}u_1^2P^2\,\mathrm{d}S
 =\mu_1\int_\Omega u_1^2P^2\,\mathrm{d}x
   +\int_\Omega u_1^2|\nabla P|^2\,\mathrm{d}x.
\end{equation}
Smooth approximation in $H^1(\Omega)$ and continuity of the trace
extend this identity to every $P\in H^1(\Omega)$. Thus
\begin{equation}\label{Equality 4.2}
 \mu_2-\mu_1
 =\inf_{\substack{P\in H^1(\Omega)\setminus\{0\}\\
       \int_\Omega Pu_1^2\,\mathrm{d}x=0}}
 \frac{\displaystyle\int_\Omega u_1^2|\nabla P|^2\,\mathrm{d}x}
      {\displaystyle\int_\Omega P^2u_1^2\,\mathrm{d}x}.
\end{equation}
Payne and Schaefer \cite{PayneSchaefer} established this
characterization in two dimensions.

Fix $\tau>0$ satisfying
\begin{equation}\label{eq:auxiliary-admissibility}
 \mu_1(B,\tau)\leq\mu_1\widetilde R^2,
 \qquad
 \rho=\left(\frac{\mu_1(B,\tau)}{\mu_1}\right)^{1/2}.
\end{equation}
Then $0<\rho\leq\widetilde R\leq1$, and scaling gives
\[
 \mu_1\left(B_\rho,\frac{\tau}{\rho}\right)=\mu_1,
 \qquad
 \mu_2\left(B_\rho,\frac{\tau}{\rho}\right)
 =\rho^{-2}\mu_2(B,\tau).
\]
Let $z$ be the positive radial first eigenfunction on this ball,
normalized by
\[
 \int_{B_\rho}z^2\,\mathrm{d}x=\int_\Omega u_1^2\,\mathrm{d}x.
\]
Up to a positive constant factor,
$z(r)=r^{-\nu}J_\nu(\sqrt{\mu_1}\,r)$ for $0<r\leq\rho$,
with its continuous value at zero.

The weights associated with this ball are
\[
 g(r)=w_\tau(r/\rho),\qquad
 \mathcal N_\tau(r)=g'(r)^2+(N-1)\frac{g(r)^2}{r^2}
                  =\rho^{-2}B_\tau(r/\rho).
\]
By Proposition~\ref{prop:weight-monotonicity}, $g$ is nondecreasing
and bounded, while $\mathcal N_\tau$ is positive and strictly decreasing.
The expansion in Section~\ref{sec:weights} gives
$g(r)=\rho^{-1}w_\tau'(0)r+O(r^3)$, so
$P_i(x)=g(|x|)x_i/|x|$ extends regularly across the origin and
$P_iu_1\in H^1(\Omega)$. The finite limit $B_\tau(0)$ also makes
$\mathcal N_\tau$ bounded there. Choosing the origin as in
\cite{AshbaughBenguria91,AshbaughBenguria92}, we have
\[
 \int_\Omega P_i u_1^2\,\mathrm{d}x=0,
 \qquad i=1,\ldots,N.
\]
Applying \eqref{Equality 4.2} to each $P_i$ and summing gives
\begin{equation}\label{Inequality 4.3}
 \mu_2-\mu_1
 \leq\frac{\displaystyle\int_\Omega\mathcal N_\tau(r)u_1^2\,\mathrm{d}x}
          {\displaystyle\int_\Omega g(r)^2u_1^2\,\mathrm{d}x},
 \qquad r=|x|,
\end{equation}
because $\sum_iP_i^2=g^2$ and
$\sum_i|\nabla P_i|^2=g'^2+(N-1)g^2/r^2$.

To compare these integrals with their rearrangements, write
$f_*(s)=f^*(|\Omega|-s)$ and $f_\star(x)=f_*(C_N|x|^N)$ for the
increasing rearrangements of a nonnegative measurable function $f$.
For nonnegative $f$ and $h$, the Hardy--Littlewood inequalities are
\begin{equation}\label{Inequality 4.1}
 \int_B f_\star h^\star\,\mathrm{d}x
 \leq\int_\Omega fh\,\mathrm{d}x
 \leq\int_B f^\star h^\star\,\mathrm{d}x.
\end{equation}
For a nonincreasing radial function $f(|x|)$, its decreasing
rearrangement from $\Omega$ to $B$ is at most $f(r)$. For a
nondecreasing radial function, its increasing rearrangement is at
least $f(r)$. Using these facts, $(u_1^2)^\star=(u_1^\star)^2$, and
\eqref{Inequality 4.1}, we obtain
\begin{align*}
 \int_\Omega\mathcal N_\tau(r)u_1^2\,\mathrm{d}x
 &\leq\int_B\mathcal N_\tau(r)(u_1^\star)^2\,\mathrm{d}x,\\
 \int_\Omega g(r)^2u_1^2\,\mathrm{d}x
 &\geq\int_Bg(r)^2(u_1^\star)^2\,\mathrm{d}x;
\end{align*}
see \cite{Hardy,AshbaughBenguria92,DaiShi}.

Since $\rho\leq\widetilde R$, Lemma~\ref{lem:matched-comparison}
compares $u_1^\star$ and $z$ on $B_\rho$. For the weighted integrals,
extend $z$ by zero for $r>\rho$. In both alternatives of that lemma there is
$r_*\in(0,\rho]$ such that $z^2-(u_1^\star)^2$ is nonnegative for
$0\leq r<r_*$ and nonpositive for $r_*<r\leq1$; in the no-crossing case,
take $r_*=\rho$. Its integral is zero by normalization. Hence, for
every nonnegative continuous nonincreasing radial weight $F$,
\begin{align}
 &\int_{B_\rho}F(r)z^2\,\mathrm{d}x
   -\int_BF(r)(u_1^\star)^2\,\mathrm{d}x\notag\\
 &\quad=NC_N\int_0^1\bigl(F(r)-F(r_*)\bigr)
       \bigl(z(r)^2-u_1^\star(r)^2\bigr)r^{N-1}\,\mathrm{d}r
 \geq0.\label{eq:matched-weight-comparison}
\end{align}
The inequality reverses for nondecreasing $F$. This single identity
covers both possible orderings of the rearranged eigenfunctions.
Taking $F=\mathcal N_\tau$ and $F=g^2$, respectively, and using
\eqref{Inequality 4.3}, we conclude that
\begin{align*}
 \mu_2-\mu_1
 &\leq\frac{\int_\Omega\mathcal N_\tau(r)u_1^2\,\mathrm{d}x}
              {\int_\Omega g(r)^2u_1^2\,\mathrm{d}x}\\
 &\leq\frac{\int_B\mathcal N_\tau(r)(u_1^\star)^2\,\mathrm{d}x}
              {\int_Bg(r)^2(u_1^\star)^2\,\mathrm{d}x}\\
 &\leq\frac{\int_{B_\rho}\mathcal N_\tau(r)z^2\,\mathrm{d}x}
              {\int_{B_\rho}g(r)^2z^2\,\mathrm{d}x}
 =\frac{\mu_2(B,\tau)-\mu_1(B,\tau)}{\rho^2}.
\end{align*}
For the last equality, the functions $z(r)g(r)x_i/r$ are second
eigenfunctions on $B_\rho$ with Robin parameter $\tau/\rho$.
Applying \eqref{eq:ground-state-identity} on that ball to each of
them and summing gives the displayed quotient. The common Robin
condition, equivalently $g'(\rho)=0$, cancels the boundary contribution.
Since $\mu_1\rho^2=\mu_1(B,\tau)$, the common comparison estimate is
\begin{equation}\label{eq:general-auxiliary-bound}
 \frac{\mu_2(\Omega,\sigma)}{\mu_1(\Omega,\sigma)}
 \leq\frac{\mu_2(B,\tau)}{\mu_1(B,\tau)}
 \quad\text{whenever }\tau>0\text{ and }
 \mu_1(B,\tau)\leq\mu_1\widetilde R^2.
\end{equation}

\section{Proofs of Theorems~\ref{thm:ball-bound} and~\ref{thm:matched-bound}}
\label{sec:proof-ball-bound}\label{sec:matched-bound}

\noindent\textbf{Proof of Theorem~\ref{thm:ball-bound}.}
Under the assumption $\widetilde R\geq R$, the choice $\tau=\sigma$
is admissible in \eqref{eq:auxiliary-admissibility}, since
$\mu_1(B,\sigma)=\mu_1R^2\leq\mu_1\widetilde R^2$. It gives
$\rho=R$, and \eqref{eq:general-auxiliary-bound} yields
\[
 \frac{\mu_2(\Omega,\sigma)}{\mu_1(\Omega,\sigma)}
 \leq\frac{\mu_2(B,\sigma)}{\mu_1(B,\sigma)}
 =\frac{\beta^2}{\alpha^2}.
\]
A translate of the unit ball attains equality. If $\Omega$ is not
a ball, Lemma~\ref{lem:comparison-radius} gives $R<1$. With $z$
extended by zero outside $B_R$, the radius $r_*$ in
\eqref{eq:matched-weight-comparison} belongs to $(0,R]$.
On the annulus $R<r<1$, we have $z=0$, $u_1^\star>0$, and
\[
 \mathcal N_\sigma(r)
 =\frac{(N-1)w(1)^2}{r^2}
 <\mathcal N_\sigma(R)\leq\mathcal N_\sigma(r_*).
\]
Thus \eqref{eq:matched-weight-comparison}, with
$F=\mathcal N_\sigma$, is strict:
\begin{align*}
 &\int_{B_R}\mathcal N_\sigma(r)z^2\,\mathrm{d}x
   -\int_B\mathcal N_\sigma(r)(u_1^\star)^2\,\mathrm{d}x\\
 &\quad=NC_N\int_0^1\bigl(\mathcal N_\sigma(r)-\mathcal N_\sigma(r_*)\bigr)
      \bigl(z(r)^2-u_1^\star(r)^2\bigr)r^{N-1}\,\mathrm{d}r>0.
\end{align*}
The integrand is nonnegative everywhere and positive on $(R,1)$.
Together with the denominator comparison, this makes
\eqref{eq:general-auxiliary-bound} strict. Equality therefore holds
if and only if $\Omega$ is a translate of $B$.\qed
\\ \\
\textbf{Proof of Theorem~\ref{thm:matched-bound}.}
The first ball eigenvalue $\mu_1(B,\tau)$ increases strictly from $0$
to $j_{\nu,1}^2$ as $\tau$ increases from $0$ to $+\infty$.
If $\mu_1\widetilde R^2<j_{\nu,1}^2$, there is a unique $\tau>0$
satisfying \eqref{eq:matched-parameter}. It is admissible in
\eqref{eq:auxiliary-admissibility}, with $\rho=\widetilde R$.
Equation~\eqref{eq:general-auxiliary-bound} proves the first case.
The ratio monotonicity in Theorem~\ref{thm:ball-monotonicity} also
shows that this largest admissible parameter gives the smallest
upper bound in \eqref{eq:general-auxiliary-bound}. On a unit ball,
$\widetilde R=1$ and $\tau=\sigma$, so equality holds.

If $\mu_1\widetilde R^2\geq j_{\nu,1}^2$, then
$\mu_1(B,\tau)<j_{\nu,1}^2\leq\mu_1\widetilde R^2$ for every finite
$\tau>0$. Letting $\tau\to+\infty$ in
\eqref{eq:general-auxiliary-bound} and using the Robin-to-Dirichlet
limit on the ball gives
\[
 \frac{\mu_2(\Omega,\sigma)}{\mu_1(\Omega,\sigma)}
 \leq\frac{j_{\nu+1,1}^2}{j_{\nu,1}^2}.
\]
This proves the second case.\qed\\

For comparison with \eqref{eq:prior-bound}, if $\widetilde R<R$,
then $\tau<\sigma$ and Theorem~\ref{thm:ball-monotonicity} gives
\[
 \frac{\mu_2(B,\sigma)}{\mu_1(B,\sigma)}
 <\frac{\mu_2(B,\tau)}{\mu_1(B,\tau)}
 =1+\frac{\mu_2(B,\tau)-\mu_1(B,\tau)}{\mu_1\widetilde R^2}
 <1+\frac{\beta^2-\alpha^2}{\mu_1\widetilde R^2}.
\]
The right-hand side is \eqref{eq:prior-bound} in this case.
If $\widetilde R=R$, then $\tau=\sigma$ and the two ball bounds
coincide. If $\widetilde R>R$, either the matching parameter is
finite and $\tau>\sigma$, or the bound is the Dirichlet ratio.
In both cases, Theorem~\ref{thm:ball-monotonicity} makes it strictly
smaller than $\mu_2(B,\sigma)/\mu_1(B,\sigma)$.

\begingroup
\section{Proof of Theorem~\ref{thm:ball-monotonicity}}
\label{sec:ball-branches}

The first two eigenvalues on the ball belong to the lowest radial
branches of spherical harmonic degrees zero and one. Their monotonicity
follows from a comparison valid for any two such branches.
For each integer $\ell\geq0$, write
\[
 a_\ell(\sigma)=k_{\nu+\ell,1}(\sigma)^2,
 \,\,\, c_\ell=\ell(\ell+N-2),\,\,\, \sigma>0.
\]
Thus $a_\ell(\sigma)$ is the lowest eigenvalue whose angular part is
a spherical harmonic of degree $\ell$. In particular,
$a_0(\sigma)=\mu_1(B,\sigma)$ and $a_1(\sigma)=\mu_2(B,\sigma)$ by
\eqref{eq:appendix-positive-first-two}.
\medskip
\begin{proposition}
\label{prop:angular-monotonicity}
Let $N\geq2$ and let $\ell>i\geq0$ be integers. For every $\sigma>0$,
\[
 a_\ell(\sigma)-a_i(\sigma)>c_\ell-c_i.
\]
The functions $a_\ell-a_i$ and $a_\ell/a_i$ are continuously
differentiable on $(0,+\infty)$ and satisfy
\[
 \frac{\mathrm{d}}{\mathrm{d}\sigma}(a_\ell-a_i)>0,
 \,\,\,
 \frac{\mathrm{d}}{\mathrm{d}\sigma}\left(\frac{a_\ell}{a_i}\right)<0.
\]
\end{proposition}
\medskip
\noindent\textbf{Proof.}
The corresponding positive radial
function can be chosen as
\[
 U_\ell(r)=r^{-\nu}J_{\nu+\ell}\left(\sqrt{a_\ell(\sigma)}\,r\right),
 \,\,\, 0<r\leq1.
\]
It solves
\begin{equation*}
 -\left(r^{N-1}U_\ell'\right)'
 +c_\ell r^{N-3}U_\ell
 =a_\ell(\sigma)r^{N-1}U_\ell,
 \,\,\, U_\ell'(1)+\sigma U_\ell(1)=0,
\end{equation*}
with $U_\ell(r)\sim C_\ell r^\ell$ as $r\to0$,
where $C_\ell>0$. In particular, $U_\ell$ is positive on $(0,1]$.
For each fixed $\ell$, the first eigenvalue of this radial problem is
simple. The multiplicity arising from the spherical harmonics therefore
does not affect its differentiation with respect to $\sigma$.

The admissible radial functions do not depend on $\sigma$, and the quadratic
form contains $\sigma$ only in the boundary term $\sigma|U(1)|^2$.
The usual perturbation formula for a simple eigenvalue shows that
$a_\ell$ is continuously differentiable and gives
\begin{equation}\label{eq:angular-hf}
 a_\ell'(\sigma)=
 \frac{U_\ell(1)^2}{\displaystyle\int_0^1r^{N-1}U_\ell(r)^2\,\mathrm{d}r}>0.
\end{equation}
This formula holds for any normalization of $U_\ell$.

Fix $\ell>i\geq0$ and $\sigma>0$. Write
\[
 u=U_i,\,\,\, v=U_\ell,\,\,\,
 d=a_\ell(\sigma)-a_i(\sigma),\,\,\, c=c_\ell-c_i>0.
\]
Define the Wronskian
\[
 W(r)=r^{N-1}\bigl(v'(r)u(r)-v(r)u'(r)\bigr).
\]
Combining the two radial equations, we obtain
\begin{equation}\label{eq:angular-wronskian}
 W'(r)=r^{N-1}u(r)v(r)\left(\frac{c}{r^2}-d\right).
\end{equation}
The behavior of $u$ and $v$ at the origin gives
$W(r)=O(r^{N-2+i+\ell})$, and therefore $W(0)=0$ because
$N\geq2$ and $i+\ell\geq1$. The common Robin condition gives
$W(1)=0$. Integrating \eqref{eq:angular-wronskian}, we obtain
\begin{equation}\label{eq:angular-gap-lower-bound}
 d=c\,
 \frac{\displaystyle\int_0^1r^{N-3}u(r)v(r)\,\mathrm{d}r}
 {\displaystyle\int_0^1r^{N-1}u(r)v(r)\,\mathrm{d}r}>c.
\end{equation}
Both integrals are finite, including when $N=2$ and $i=0$, because
$u(r)v(r)=O(r^{i+\ell})$.

It follows that $r_0=\sqrt{c/d}$ belongs to $(0,1)$.
Equation \eqref{eq:angular-wronskian} shows that $W$ is strictly
increasing on $(0,r_0)$ and strictly decreasing on $(r_0,1)$.
Since $W$ vanishes at both endpoints, $W(r)>0$ for $0<r<1$.
It follows that
\[
 \left(\frac{v}{u}\right)'=
 \frac{W}{r^{N-1}u^2}>0,
 \,\,\,
 \frac{v(r)}{v(1)}<\frac{u(r)}{u(1)}\quad(0<r<1).
\]
Squaring the last inequality, integrating against $r^{N-1}\mathrm{d}r$, and
using \eqref{eq:angular-hf}, we obtain
\[
 a_\ell'(\sigma)
 =\left(\int_0^1r^{N-1}\frac{v(r)^2}{v(1)^2}\,\mathrm{d}r\right)^{-1}
 >\left(\int_0^1r^{N-1}\frac{u(r)^2}{u(1)^2}\,\mathrm{d}r\right)^{-1}
 =a_i'(\sigma).
\]
Thus $a_\ell-a_i$ is strictly increasing.

The boundary condition gives
\[
 \frac{\sqrt{a_\ell(\sigma)}
 J_{\nu+\ell+1}(\sqrt{a_\ell(\sigma)})}
 {J_{\nu+\ell}(\sqrt{a_\ell(\sigma)})}=\sigma+\ell.
\]
The denominator is nonzero because
$0<\sqrt{a_\ell(\sigma)}<j_{\nu+\ell,1}$.
Put $m=\nu+\ell$, $x=\sqrt{a_\ell(\sigma)}$, and
$F_m(x)=xJ_{m+1}(x)/J_m(x)$. The recurrence relations
\[
 J_m'(x)=\frac m xJ_m(x)-J_{m+1}(x),\,\,\,
 J_{m+1}'(x)=J_m(x)-\frac{m+1}{x}J_{m+1}(x)
\]
give
\begin{align*}
 F_m'(x)
 &=\frac{J_{m+1}}{J_m}
   +x\frac{J_{m+1}'J_m-J_{m+1}J_m'}{J_m^2}\notag\\
 &=x+\frac{F_m(x)^2-2mF_m(x)}x.
\end{align*}
At the root, $F_m(x)=\sigma+\ell$, and therefore
\[
 xF_m'(x)
 =a_\ell(\sigma)+(\sigma+\ell)^2
                   -2(\nu+\ell)(\sigma+\ell)
 =a_\ell(\sigma)+\sigma^2-(N-2)\sigma-c_\ell.
\]
Define the last expression to be $D_\ell(\sigma)$.
The chain rule gives
\[
 1=\frac{\mathrm{d}}{\mathrm{d}\sigma}F_m\bigl(\sqrt{a_\ell(\sigma)}\bigr)
   =F_m'(x)\frac{a_\ell'(\sigma)}{2x}
   =\frac{D_\ell(\sigma)a_\ell'(\sigma)}{2a_\ell(\sigma)}.
\]
Consequently,
\begin{equation}\label{eq:angular-bessel-derivative}
 \frac{a_\ell'(\sigma)}{a_\ell(\sigma)}=\frac{2}{D_\ell(\sigma)},
 \,\,\,
 D_\ell(\sigma)=a_\ell(\sigma)+\sigma^2-(N-2)\sigma-c_\ell.
\end{equation}
Since $a_\ell>0$ and $a_\ell'>0$ by \eqref{eq:angular-hf},
$D_\ell(\sigma)>0$. Combining \eqref{eq:angular-bessel-derivative}
with \eqref{eq:angular-gap-lower-bound}, we conclude that
\begin{align*}
 \frac{\mathrm{d}}{\mathrm{d}\sigma}\log\frac{a_\ell(\sigma)}{a_i(\sigma)}
 &=\frac{2}{D_\ell(\sigma)}-\frac{2}{D_i(\sigma)}\\
 &=-\frac{2\bigl[(a_\ell(\sigma)-a_i(\sigma))-(c_\ell-c_i)\bigr]}
 {D_\ell(\sigma)D_i(\sigma)}<0.
\end{align*}
Thus $a_\ell/a_i$ is strictly decreasing.\qed
\\ \\
\textbf{Proof of Theorem~\ref{thm:ball-monotonicity}.}
Taking $(\ell,i)=(1,0)$ in Proposition~\ref{prop:angular-monotonicity}
gives the asserted differentiability and strict derivative signs for
the gap and ratio of the first two ball eigenvalues. By continuity
with respect to the Robin parameter, as $\sigma\downarrow0$ we have
\[
 a_0(\sigma)\longrightarrow0,\,\,\,
 a_1(\sigma)\longrightarrow\mu_2(B,0)>0.
\]
The positivity of the second limit follows because the unit ball is
connected. Hence the ratio tends to $+\infty$. As
$\sigma\to+\infty$, each radial eigenvalue converges monotonically
to the corresponding Dirichlet eigenvalue, so
\[
 a_\ell(\sigma)\longrightarrow j_{\nu+\ell,1}^2
 \,\,\,(\sigma\to+\infty);
\]
see Appendix. Hence the ratio tends to
$j_{\nu+1,1}^2/j_{\nu,1}^2$, and continuity and strict decrease give
the stated range and uniqueness.
\qed\\

The argument relies on the positivity of the first radial
eigenfunctions and does not extend to higher radial modes, which
have interior zeros. Its conclusion concerns angular branches,
not arbitrary pairs in the ordered spectrum.
\par\endgroup

\begingroup
\section{On the restriction on the domain}
\label{sec:necessity}

To construct the domains in Theorem~\ref{thm:unbounded-ratio}, fix
$\sigma>0$ on the whole boundary. As in Bucur's example
\cite{Laugesen} and the construction in \cite{Nahon2022}, the aim is
to separate the low-energy part from the part carrying the remaining
volume. Here a small cube joined to a thin tube gives a connected
domain. A trial function concentrated in the cube bounds the first
eigenvalue from above, while the transverse Robin energy in the tube
permits a lower bound on the second.
\medskip
\begin{proposition}
\label{prop:thin-appendage}
Fix $N\geq2$, $V>0$, and $\sigma>0$. For
$0<a<V^{1/N}$ and $0<\varepsilon<a$, put
\[
 Q_a=(0,a)\times(-a/2,a/2)^{N-1},\,\,\,
 T_{a,\varepsilon}=(a,a+L)\times(-\varepsilon/2,\varepsilon/2)^{N-1},
 \,\,\, L=\frac{V-a^N}{\varepsilon^{N-1}},
\]
and define
\[
 \Omega_{a,\varepsilon}
 =\operatorname{int}\big(\overline{Q_a}\cup\overline{T_{a,\varepsilon}}\big).
\]
Suppose that
\begin{equation}\label{eq:tail-threshold}
 \min\left\{\varepsilon^{-2},\frac{\sigma}{2\varepsilon}\right\}
 \geq\frac{\pi^2}{a^2}.
\end{equation}
Then, for every $0<\delta<L$,
\begin{align}
 \mu_2(\Omega_{a,\varepsilon},\sigma)
 &\geq\frac{\pi^2}{a^2},\label{eq:tail-mu2}\\
 \mu_1(\Omega_{a,\varepsilon},\sigma)
 &\leq
 \frac{\sigma(2Na^{N-1}-\varepsilon^{N-1})
       +\varepsilon^{N-1}/\delta
       +2(N-1)\sigma\varepsilon^{N-2}\delta/3}
      {a^N+\varepsilon^{N-1}\delta/3}.
 \label{eq:tail-mu1}
\end{align}
\end{proposition}
\medskip
\noindent\textbf{Proof.} The set $\Omega_{a,\varepsilon}$ is a bounded connected Lipschitz domain
of volume $V$; when $N=2$ it is a simply connected polygon. The cube and
the tube meet along
\[
 \Gamma=\{a\}\times(-\varepsilon/2,\varepsilon/2)^{N-1}.
\]
We separate the cube and the tube by imposing Neumann conditions on the
two sides of $\Gamma$, and keep the Robin condition with parameter
$\sigma$ on the remaining faces. Functions on the two pieces no longer
have to agree on $\Gamma$, so the domain of the quadratic form is
enlarged. By the min--max principle, the eigenvalues of the resulting
direct sum are lower bounds for the corresponding eigenvalues on
$\Omega_{a,\varepsilon}$.

On the cube, the quadratic form for these mixed boundary conditions is
at least the Neumann quadratic form, since $\sigma>0$. Thus its second
eigenvalue is at least $\pi^2/a^2$. To estimate the first eigenvalue of the tube,
observe that every $h\in H^1(0,\varepsilon)$ satisfies
\[
 \int_0^\varepsilon |h(s)|^2\,\mathrm{d}s
 \leq2\varepsilon |h(0)|^2
       +\varepsilon^2\int_0^\varepsilon |h'(s)|^2\,\mathrm{d}s.
\]
Indeed, the absolutely continuous representative of $h$ satisfies
\[
 |h(s)|^2\leq2|h(0)|^2
       +2s\int_0^s|h'(\zeta)|^2\,\mathrm{d}\zeta.
\]
After integration in $s$, Fubini's theorem gives
\[
 2\int_0^\varepsilon s\int_0^s|h'(\zeta)|^2\,\mathrm{d}\zeta\,\mathrm{d}s
 =\int_0^\varepsilon(\varepsilon^2-\zeta^2)
           |h'(\zeta)|^2\,\mathrm{d}\zeta
 \leq\varepsilon^2\int_0^\varepsilon|h'|^2\,\mathrm{d}s.
\]
Put $c_\varepsilon=\min\{\varepsilon^{-2},\sigma/(2\varepsilon)\}$.
Since $c_\varepsilon\varepsilon^2\leq1$ and
$2c_\varepsilon\varepsilon\leq\sigma$, multiplying the preceding
$L^2$ estimate by $c_\varepsilon$ yields
\[
 c_\varepsilon\int_0^\varepsilon|h|^2\,\mathrm{d}s
 \leq\int_0^\varepsilon|h'|^2\,\mathrm{d}s+\sigma|h(0)|^2.
\]
Adding the nonnegative term $\sigma|h(\varepsilon)|^2$ on the
right proves
\begin{equation}\label{eq:transverse-elementary}
 \int_0^\varepsilon |h'|^2\,\mathrm{d}s
 +\sigma\big(|h(0)|^2+|h(\varepsilon)|^2\big)
 \geq
 \min\left\{\varepsilon^{-2},\frac{\sigma}{2\varepsilon}\right\}
 \int_0^\varepsilon |h|^2\,\mathrm{d}s.
\end{equation}
Apply \eqref{eq:transverse-elementary} along one transverse coordinate of
$T_{a,\varepsilon}$, and integrate over all other coordinates. The
gradient terms in the remaining directions and the other Robin boundary
terms are nonnegative. Thus the first eigenvalue of the tube with these
mixed boundary conditions is at least the left-hand side of
\eqref{eq:tail-threshold}.
The Neumann spectrum of the cube is
\[
 \frac{\pi^2}{a^2}(m_1^2+\cdots+m_N^2),
 \,\,\, m_1,\ldots,m_N\in\{0,1,2,\ldots\}.
\]
Only its first eigenvalue, the simple eigenvalue zero, lies below
$\pi^2/a^2$. Adding the nonnegative Robin term on the exterior cube
faces cannot lower any eigenvalue. Thus the separated cube has at
most one eigenvalue below $\pi^2/a^2$, counted with multiplicity,
while the tube has none. The spectrum of the direct sum is the
combined list of the two spectra, so its second eigenvalue is at
least $\pi^2/a^2$. The min--max comparison at $\Gamma$ gives
\[
 \mu_2(\Omega_{a,\varepsilon},\sigma)
 \geq\frac{\pi^2}{a^2}.
\]
This proves \eqref{eq:tail-mu2}.

For the first eigenvalue, use the function $f$ equal to $1$ on $Q_a$
and given on the tube, with $s=x_1-a$, by
\[
 f(x)=\max\{1-s/\delta,0\}.
\]
The two definitions agree on $\Gamma$, so
$f\in H^1(\Omega_{a,\varepsilon})$.
Direct integration gives
\[
 \int_{\Omega_{a,\varepsilon}} f^2\,\mathrm{d}x
 =a^N+\frac{\varepsilon^{N-1}\delta}{3},\,\,\,
 \int_{\Omega_{a,\varepsilon}} |\nabla f|^2\,\mathrm{d}x
 =\frac{\varepsilon^{N-1}}{\delta},
\]
and
\[
 \int_{\partial\Omega_{a,\varepsilon}} f^2\,\mathrm{d}S
 =2Na^{N-1}-\varepsilon^{N-1}
       +\frac{2(N-1)\varepsilon^{N-2}\delta}{3}.
\]
There is no contribution from the far end of the tube because
$\delta<L$. Substitution into the Rayleigh quotient proves
\eqref{eq:tail-mu1}.\qed\\
\\
\textbf{Proof of Theorem~\ref{thm:unbounded-ratio}.}
Set $\varepsilon=a^4$ and $\delta=a^2$, and let $a\downarrow0$.
For all sufficiently small $a$, we have $\varepsilon<a$, $\delta<L$,
and \eqref{eq:tail-threshold}. Proposition~\ref{prop:thin-appendage}
then yields, by omitting the negative term
$-\sigma\varepsilon^{N-1}$ in the numerator of
\eqref{eq:tail-mu1} and noting that its denominator is at least $a^N$,
\begin{equation}\label{eq:unbounded-ratio-quantitative}
 \frac{\mu_2(\Omega_{a,a^4},\sigma)}
      {\mu_1(\Omega_{a,a^4},\sigma)}
 \geq
 \frac{\pi^2}
 {2N\sigma a+\left(1+\frac{2(N-1)\sigma}{3}\right)a^{3N-4}}
 \longrightarrow+\infty.
\end{equation}
Here $3N-4>0$ for every $N\geq2$. Each domain in this family is bounded,
although their diameters need not have a common bound.

To obtain smooth boundaries, fix a sufficiently small $a$
for which the ratio in \eqref{eq:unbounded-ratio-quantitative} is greater
than $K$, and write $D=\Omega_{a,a^4}$. By
\cite[Theorem~8.3.1]{DanersSurvey2008}, we can approximate this fixed
Lipschitz domain by smooth domains $D_j$ with uniform Lipschitz control
and with the gradients of the boundary graphs converging strongly in
every finite $L^p$ space. Such an approximation preserves the Robin boundary term.
For this fixed domain, the smoothing can be confined to boundary
graph neighborhoods at a scale smaller than the fixed tube width.
The perturbations can be chosen to leave connected interior cores
of the cube and tube joined through the neck, so connectedness is
preserved. In dimension two, one rounds the vertices along the
simple closed polygonal boundary without creating intersections;
the resulting smooth boundary is again a simple closed curve and
the enclosed domain is simply connected. These geometric choices
are made for the fixed $a$, before sending the smoothing scale to zero.
By
\cite[Theorem~6.5.1 and Remark~8.3.2]{DanersSurvey2008}, the limiting
boundary density is one. Hence the Robin parameter is unchanged and the
resolvents converge in operator norm. In particular,
\begin{equation}\label{eq:smooth-eigenvalue-limit}
 |D_j|\longrightarrow |D|=V,\,\,\,
 \mu_i(D_j,\sigma)\longrightarrow\mu_i(D,\sigma),
 \quad i=1,2.
\end{equation}
To make the volume exactly $V$, put
$\rho_j=(V/|D_j|)^{1/N}$ and $\widehat D_j=\rho_jD_j$.
Then $\rho_j\to1$ and $|\widehat D_j|=V$. The scaling identity is
\[
 \mu_i(\widehat D_j,\sigma)
 =\rho_j^{-2}\mu_i(D_j,\rho_j\sigma).
\]
For every bounded Lipschitz domain $E$ and every $\rho>0$, positivity of
the gradient and boundary terms in the Robin quadratic form gives
\[
 \min\{1,\rho\}\,\mu_i(E,\sigma)
 \leq\mu_i(E,\rho\sigma)
 \leq\max\{1,\rho\}\,\mu_i(E,\sigma).
\]
By the min--max principle, the constants are independent of $E$.
Combining this with
\eqref{eq:smooth-eigenvalue-limit} proves
\[
 \mu_i(\widehat D_j,\sigma)\longrightarrow\mu_i(D,\sigma),
 \,\,\, i=1,2.
\]
It follows that the ratio is still greater than $K$ for all sufficiently
large $j$. The domains $\widehat D_j$ are smooth, bounded, and connected,
and can be chosen simply connected when $N=2$.\qed\\

For example, in dimension two take
\[
 V=\pi,\quad\sigma=1,\quad a=\frac12,\quad
 \varepsilon=\frac1{100},\quad\delta=\frac1{10},\quad
 L=100\left(\pi-\frac14\right).
\]
The transverse lower bound is $50>4\pi^2$, and
\eqref{eq:tail-mu2}--\eqref{eq:tail-mu1} give
\[
 \mu_2\geq4\pi^2,\,\,\, \mu_1\leq\frac{6470}{751},\,\,\,
 \frac{\mu_2}{\mu_1}\geq\frac{1502\pi^2}{3235}
 \approx4.582425.
\]
The positive function $\phi(x)=1-|x|^2/3$ on the unit disk satisfies
$\partial_{\mathbf n}\phi+\phi=0$ and
$-\Delta\phi=4/3\geq(4/3)\phi$. Green's identity against the positive
first Robin eigenfunction gives $\mu_1(B,1)\geq4/3$.
The odd function $v(x)=x_1(1-|x|^2/2)$ is orthogonal to that radial
first eigenfunction. Its Rayleigh quotient is
\[
 \frac{\int_B|\nabla v|^2\,\mathrm{d}x+\int_{\partial B}v^2\,\mathrm{d}S}
      {\int_Bv^2\,\mathrm{d}x}
 =\frac{2\pi/3}{11\pi/96}=\frac{64}{11}.
\]
Therefore
\[
 \frac{\mu_2(B,1)}{\mu_1(B,1)}
 \leq\frac{48}{11}
 <\frac{1502\pi^2}{3235}
 \leq\frac{\mu_2(\Omega_{1/2,1/100},1)}
              {\mu_1(\Omega_{1/2,1/100},1)}.
\]
The strict inequality persists under the smooth approximation and
volume correction above.

The preceding construction varies the domain at a fixed parameter.
On each fixed nonball domain, however, comparison with the equal-volume
ball holds for all sufficiently small and all sufficiently large
positive Robin parameters.
\medskip
\begin{proposition}
\label{prop:endpoint-regimes}
Let $\Omega\subset\mathbb R^N$ be a fixed bounded connected Lipschitz
domain, $N\geq2$, and let $B_\rho$ be a ball of the same volume. If
$\Omega$ is not a ball, there exist constants
$0<\sigma_-(\Omega)<\sigma_+(\Omega)<+\infty$ such that
\begin{equation}\label{eq:endpoint-ball-comparison}
 \frac{\mu_2(\Omega,\sigma)}{\mu_1(\Omega,\sigma)}
 <\frac{\mu_2(B_\rho,\sigma)}{\mu_1(B_\rho,\sigma)}
\end{equation}
whenever $0<\sigma<\sigma_-(\Omega)$ or $\sigma>\sigma_+(\Omega)$.
For a ball the two ratios agree for every $\sigma>0$.
\end{proposition}
\medskip
\noindent\textbf{Proof.}
Write $V=|\Omega|$, $P(\Omega)=\mathcal H^{N-1}(\partial\Omega)$,
and denote by $\kappa_2(\Omega)>0$ the first nonzero Neumann eigenvalue.
The Neumann eigenvalue $0$ is simple since $\Omega$ is connected.
Perturbing the Robin quadratic form at $\sigma=0$, whose normalized
first eigenfunction is the constant $V^{-1/2}$, gives
\[
 \mu_1(\Omega,\sigma)=\frac{P(\Omega)}{V}\sigma+o(\sigma),
 \,\,\,
 \mu_2(\Omega,\sigma)\longrightarrow\kappa_2(\Omega).
\]
See \cite{Bucur,Laugesen}. Hence
\begin{equation}\label{eq:small-robin-ratio-limit}
 \lim_{\sigma\downarrow0}
 \sigma\frac{\mu_2(\Omega,\sigma)}{\mu_1(\Omega,\sigma)}
 =\frac{V\kappa_2(\Omega)}{P(\Omega)}.
\end{equation}
The Szeg\"o--Weinberger inequality and the strict isoperimetric inequality
for a nonball domain imply
\[
 \frac{V\kappa_2(\Omega)}{P(\Omega)}
 \leq\frac{V\kappa_2(B_\rho)}{P(\Omega)}
 <\frac{V\kappa_2(B_\rho)}{P(B_\rho)}.
\]
Applying \eqref{eq:small-robin-ratio-limit} also to the ball proves
\eqref{eq:endpoint-ball-comparison} for sufficiently small $\sigma>0$.

As $\sigma\to+\infty$, monotone convergence of the Robin forms to the
Dirichlet form gives
\[
 \mu_i(\Omega,\sigma)\longrightarrow\lambda_i(\Omega),\,\,\,
 \mu_i(B_\rho,\sigma)\longrightarrow\lambda_i(B_\rho),
 \quad i=1,2;
\]
see \cite[Proposition~4.5]{Bucur}. The strict Dirichlet PPW inequality
\cite{AshbaughBenguria91,AshbaughBenguria92} states that
\[
 \frac{\lambda_2(\Omega)}{\lambda_1(\Omega)}
 <\frac{\lambda_2(B_\rho)}{\lambda_1(B_\rho)}
\]
when $\Omega$ is not a ball. This strict inequality holds for all
sufficiently large finite $\sigma$, which proves the second assertion.
If necessary, decrease the small threshold or increase the large
threshold to ensure $\sigma_-<\sigma_+$.
\qed\\

By Theorem~\ref{thm:unbounded-ratio}, the thresholds in
Proposition~\ref{prop:endpoint-regimes} cannot be uniform over domains
of fixed volume.
\par\endgroup

\begingroup
\appendix
\makeatletter
\renewcommand{\@seccntformat}[1]{\appendixname.\quad}
\makeatother
\section{Robin spectra on balls and intervals}
\label{app:ball-spectrum}

The ball eigenvalues and eigenfunctions used in the main text arise
from separation of variables. The same calculation describes the
positive, zero, and negative spectral branches for $\sigma\in\mathbb R$
and explains the root locations and multiplicities in the figures
and numerical tables; see \cite[Section~5]{Freitas},
\cite[Section~3]{HassannezhadSher}.

For each spherical harmonic degree $\ell\geq0$, let
$\mu_{\ell,m}(\sigma)$ be the $m$th radial eigenvalue, numbered
from the lowest and including a possible nonpositive first value.
Thus $\mu_{\ell,1}(\sigma)=a_\ell(\sigma)$ for $\sigma>0$.
This radial index is distinct from the index in the full spectrum,
where angular multiplicities are included.

For a separated eigenfunction $u(r,\xi)=v(r)G_\ell(\xi)$, where
$G_\ell$ is a spherical harmonic of degree $\ell$, the identity
$\Delta_{\mathbb S^{N-1}}G_\ell=-c_\ell G_\ell$ gives
\[
 v''+\frac{N-1}{r}v'+\left(\mu-\frac{c_\ell}{r^2}\right)v=0,
 \qquad 0<r<1,
 \qquad v'(1)+\sigma v(1)=0.
\]
Here $c_\ell$ is as in Section~\ref{sec:ball-branches}, and the
radial solution is required to be regular at the origin.

For $\mu=k^2>0$, the substitution $v(r)=r^{-\nu}J(kr)$ reduces
the radial equation to the Bessel equation
\[
 J''(x)+\frac1xJ'(x)
       +\left(1-\frac{(\nu+\ell)^2}{x^2}\right)J(x)=0.
\]
The regular radial factor is therefore
$v(r)=r^{-\nu}J_{\nu+\ell}(kr)$, up to a constant. The recurrence
relation $J_\eta'(x)=\eta J_\eta(x)/x-J_{\eta+1}(x)$ gives
\[
 v'(1)+\sigma v(1)
 = (\sigma+\ell)J_{\nu+\ell}(k)-kJ_{\nu+\ell+1}(k).
\]
Since consecutive Bessel functions have no common positive zero,
the boundary condition is equivalent to \eqref{eq:ball-root-definition}.
Write it as $\sigma=\widetilde h_{\nu+\ell}(k)$, where
\[
 \widetilde h_{\nu+\ell}(k)
 =k\frac{J_{\nu+\ell+1}(k)}{J_{\nu+\ell}(k)}-\ell.
\]
The Bessel product representation \cite[Chapter~10]{DLMF} gives
\[
 \widetilde h_{\nu+\ell}(k)+\ell
 =\sum_{n=1}^{+\infty}\frac{2k^2}{j_{\nu+\ell,n}^2-k^2}.
\]
The series and its derivative converge locally away from the poles.
Consequently, $\widetilde h_{\nu+\ell}$ is strictly increasing between
consecutive poles $j_{\nu+\ell,m}$. On the first interval $(0,j_{\nu+\ell,1})$
it increases from $-\ell$ to $+\infty$; on every subsequent
interval it increases from $-\infty$ to $+\infty$.
Accordingly, if $\sigma>-\ell$, all radial eigenvalues are
positive and $\mu_{\ell,m}=k_{\nu+\ell,m}^2$, where
$k_{\nu+\ell,m}$ is the $m$th positive root. If $\sigma\leq-\ell$,
the first radial eigenvalue is nonpositive, so for $m\geq2$
the symbol $k_{\nu+\ell,m}$ denotes the $(m-1)$th positive root,
not the $m$th. In either case,
\[
 j_{\nu+\ell,m-1}<k_{\nu+\ell,m}<j_{\nu+\ell,m},
 \,\,\, m\geq2.
\]
The interlacing of Bessel zeros gives the more precise locations
\[
 \begin{cases}
 j_{\nu+\ell+1,m-1}<k_{\nu+\ell,m}<j_{\nu+\ell,m},
       &\sigma>-\ell,\\
 k_{\nu+\ell,m}=j_{\nu+\ell+1,m-1},&\sigma=-\ell,\\
 j_{\nu+\ell,m-1}<k_{\nu+\ell,m}<j_{\nu+\ell+1,m-1},
       &\sigma<-\ell,
 \end{cases}
 \qquad m\geq2.
\]
Indeed, $\widetilde h_{\nu+\ell}(j_{\nu+\ell+1,m-1})=-\ell$,
so these inequalities follow from its strict increase between poles.
When $\sigma>-\ell$, the first root lies in $(0,j_{\nu+\ell,1})$.

For a negative eigenvalue $\mu=-k^2$, with $k>0$, the regular
radial factor is $v(r)=r^{-\nu}I_{\nu+\ell}(kr)$. Using
$I_\eta'(x)=I_{\eta+1}(x)+\eta I_\eta(x)/x$ in the Robin condition gives
\[
 (\sigma+\ell)I_{\nu+\ell}(k)+kI_{\nu+\ell+1}(k)=0.
\]
Thus the corresponding root equation is
\[
 \sigma=\widehat h_{\nu+\ell}(k)
 :=-k\frac{I_{\nu+\ell+1}(k)}{I_{\nu+\ell}(k)}-\ell,
\]
where $I_{\nu+\ell}$ is the modified Bessel function of the first
kind. Its product representation gives
\[
 \widehat h_{\nu+\ell}(k)+\ell
 =-\sum_{n=1}^{+\infty}\frac{2k^2}{j_{\nu+\ell,n}^2+k^2}.
\]
Hence $\widehat h_{\nu+\ell}$ is strictly decreasing, with value
$-\ell$ at zero. Since $I_{\nu+\ell+1}(k)/I_{\nu+\ell}(k)\to1$
as $k\to+\infty$, its limit is $-\infty$.
Thus exactly one negative eigenvalue occurs for that degree when
$\sigma<-\ell$:
\[
 \mu_{\ell,1}=-\widehat k_{\nu+\ell,1}^{\,2},
 \,\,\, \widehat h_{\nu+\ell}(\widehat k_{\nu+\ell,1})=\sigma.
\]
For $\mu=0$, the radial equation is an Euler equation whose regular
solution is $v(r)=r^\ell$. The boundary condition becomes
$\ell+\sigma=0$. Thus at $\sigma=-\ell$ the first radial eigenvalue
is zero, with regular radial factor $r^\ell$; no other radial
eigenvalue of that degree is nonpositive.
In Figure~\ref{fig:appendix-ball-roots}, $h_{\nu+\ell}$ stands for
$\widetilde h_{\nu+\ell}$ in the positive panels and
$\widehat h_{\nu+\ell}$ in the negative panels. The root labels
there use the radial numbering just specified.

\begin{figure}[htbp]

\centering
\subfloat[{$\ell=0$, $\mu>0$}]{%

\resizebox{0.48\textwidth}{!}{%
\begingroup%
\setlength{\unitlength}{1bp}%
\fontsize{10.95}{13.14}\selectfont%
\begin{picture}(260,180)%
\PPWxaxis{68}{90}{244}%
\PPWyaxis{90}{16}{165}%
\linethickness{0.8bp}
\qbezier(131.5257,154.4818)(131.5257,151.7373)(131.5257,148.9928)
\qbezier(131.5257,145.1318)(131.5257,142.3873)(131.5257,139.6428)
\qbezier(131.5257,135.7818)(131.5257,133.0373)(131.5257,130.2928)
\qbezier(131.5257,126.4318)(131.5257,123.6873)(131.5257,120.9428)
\qbezier(131.5257,117.0818)(131.5257,114.3373)(131.5257,111.5928)
\qbezier(131.5257,107.7318)(131.5257,104.9873)(131.5257,102.2428)
\qbezier(131.5257,98.3818)(131.5257,95.6373)(131.5257,92.8928)
\qbezier(131.5257,89.0318)(131.5257,86.2873)(131.5257,83.5428)
\qbezier(131.5257,79.6818)(131.5257,76.9373)(131.5257,74.1928)
\qbezier(131.5257,70.3318)(131.5257,67.5873)(131.5257,64.8428)
\qbezier(131.5257,60.9818)(131.5257,58.2373)(131.5257,55.4928)
\qbezier(131.5257,51.6318)(131.5257,48.8873)(131.5257,46.1428)
\qbezier(131.5257,42.2818)(131.5257,39.5373)(131.5257,36.7928)
\qbezier(131.5257,32.9318)(131.5257,30.3854)(131.5257,27.839)
\linethickness{0.8bp}
\qbezier(222.9561,153.4588)(222.9561,150.7143)(222.9561,147.9698)
\qbezier(222.9561,144.1088)(222.9561,141.3643)(222.9561,138.6198)
\qbezier(222.9561,134.7588)(222.9561,132.0143)(222.9561,129.2698)
\qbezier(222.9561,125.4088)(222.9561,122.6643)(222.9561,119.9198)
\qbezier(222.9561,116.0588)(222.9561,113.3143)(222.9561,110.5698)
\qbezier(222.9561,106.7088)(222.9561,103.9643)(222.9561,101.2198)
\qbezier(222.9561,97.3588)(222.9561,94.6143)(222.9561,91.8698)
\qbezier(222.9561,88.0088)(222.9561,85.2643)(222.9561,82.5198)
\qbezier(222.9561,78.6588)(222.9561,75.9143)(222.9561,73.1698)
\qbezier(222.9561,69.3088)(222.9561,66.5643)(222.9561,63.8198)
\qbezier(222.9561,59.9588)(222.9561,57.2143)(222.9561,54.4698)
\qbezier(222.9561,50.6088)(222.9561,47.8643)(222.9561,45.1198)
\qbezier(222.9561,41.2588)(222.9561,38.5143)(222.9561,35.7698)
\qbezier(222.9561,31.9088)(222.9561,29.3623)(222.9561,26.8159)
\linethickness{1.4bp}
\qbezier(89.3588,90.1426)(92.7533,91.2619)(95.9374,93.1238)
\qbezier(95.9374,93.1238)(100.0718,95.5605)(103.7189,99.02)
\qbezier(103.7189,99.02)(106.3764,101.5463)(108.7671,104.458)
\qbezier(108.7671,104.458)(111.1577,107.3698)(113.3277,110.5476)
\qbezier(113.3277,110.5476)(115.4681,113.6544)(117.3797,116.9733)
\qbezier(117.3797,116.9733)(119.2914,120.2922)(120.8984,123.8688)
\qbezier(120.8984,123.8688)(123.3554,129.3963)(124.9824,135.3858)
\qbezier(124.9824,135.3858)(127.3165,143.875)(128.0565,152.8098)
\put(179.1262,89.0869){\rule{0.6383bp}{0.143bp}}
\put(178.9088,89.2189){\rule{1.0733bp}{0.143bp}}
\put(178.7727,89.3509){\rule{1.3454bp}{0.143bp}}
\put(178.674,89.4829){\rule{1.5429bp}{0.143bp}}
\put(178.5994,89.6149){\rule{1.6919bp}{0.143bp}}
\put(178.5432,89.7469){\rule{1.8044bp}{0.143bp}}
\put(178.502,89.8789){\rule{1.8868bp}{0.143bp}}
\put(178.4736,90.0109){\rule{1.9436bp}{0.143bp}}
\put(178.4581,90.1429){\rule{1.9745bp}{0.143bp}}
\put(178.4532,90.2749){\rule{1.9843bp}{0.143bp}}
\put(178.4599,90.4069){\rule{1.971bp}{0.143bp}}
\put(178.4778,90.5389){\rule{1.9352bp}{0.143bp}}
\put(178.5074,90.6709){\rule{1.8759bp}{0.143bp}}
\put(178.5513,90.8029){\rule{1.7881bp}{0.143bp}}
\put(178.6105,90.9349){\rule{1.6699bp}{0.143bp}}
\put(178.6885,91.0669){\rule{1.5138bp}{0.143bp}}
\put(178.7918,91.1989){\rule{1.3073bp}{0.143bp}}
\put(178.9369,91.3309){\rule{1.017bp}{0.143bp}}
\put(179.1548,91.4629){\rule{0.5813bp}{0.1097bp}}
\linethickness{1bp}
\qbezier(180.4379,90.33)(180.4273,90.858)(180.1472,91.2086)
\qbezier(180.1472,91.2086)(179.8671,91.5592)(179.4454,91.5726)
\qbezier(179.4454,91.5726)(179.0238,91.5592)(178.7437,91.2086)
\qbezier(178.7437,91.2086)(178.4636,90.858)(178.4529,90.33)
\qbezier(178.4529,90.33)(178.4636,89.8019)(178.7437,89.4511)
\qbezier(178.7437,89.4511)(179.0238,89.1002)(179.4454,89.0869)
\qbezier(179.4454,89.0869)(179.8671,89.1002)(180.1472,89.4511)
\qbezier(180.1472,89.4511)(180.4273,89.8019)(180.4379,90.33)
\put(208.918,117.7526){\rule{0.6383bp}{0.143bp}}
\put(208.7005,117.8846){\rule{1.0733bp}{0.143bp}}
\put(208.5645,118.0166){\rule{1.3453bp}{0.143bp}}
\put(208.4657,118.1486){\rule{1.5428bp}{0.143bp}}
\put(208.3912,118.2806){\rule{1.6919bp}{0.143bp}}
\put(208.335,118.4126){\rule{1.8043bp}{0.143bp}}
\put(208.2938,118.5446){\rule{1.8868bp}{0.143bp}}
\put(208.2653,118.6766){\rule{1.9436bp}{0.143bp}}
\put(208.2498,118.8086){\rule{1.9745bp}{0.143bp}}
\put(208.2449,118.9406){\rule{1.9843bp}{0.143bp}}
\put(208.2516,119.0726){\rule{1.971bp}{0.143bp}}
\put(208.2695,119.2046){\rule{1.9352bp}{0.143bp}}
\put(208.2991,119.3366){\rule{1.876bp}{0.143bp}}
\put(208.343,119.4686){\rule{1.7883bp}{0.143bp}}
\put(208.4021,119.6006){\rule{1.6701bp}{0.143bp}}
\put(208.4801,119.7326){\rule{1.5142bp}{0.143bp}}
\put(208.5832,119.8646){\rule{1.3079bp}{0.143bp}}
\put(208.7281,119.9966){\rule{1.0181bp}{0.143bp}}
\put(208.9459,120.1286){\rule{0.5825bp}{0.1102bp}}
\linethickness{1bp}
\qbezier(210.2296,118.9957)(210.219,119.5238)(209.9389,119.8746)
\qbezier(209.9389,119.8746)(209.6588,120.2255)(209.2371,120.2388)
\qbezier(209.2371,120.2388)(208.8155,120.2255)(208.5354,119.8746)
\qbezier(208.5354,119.8746)(208.2553,119.5238)(208.2447,118.9957)
\qbezier(208.2447,118.9957)(208.2553,118.4677)(208.5354,118.1168)
\qbezier(208.5354,118.1168)(208.8155,117.7659)(209.2371,117.7526)
\qbezier(209.2371,117.7526)(209.6588,117.7659)(209.9389,118.1168)
\qbezier(209.9389,118.1168)(210.219,118.4677)(210.2296,118.9957)
\linethickness{0.8bp}
\qbezier(209.2371,115.2998)(209.2371,112.5553)(209.2371,109.8108)
\qbezier(209.2371,105.9498)(209.2371,103.2053)(209.2371,100.4608)
\qbezier(209.2371,96.5998)(209.2371,93.8553)(209.2371,91.1108)
\linethickness{1.4bp}
\qbezier(179.6474,90.0769)(172.9812,85.2777)(166.4994,80.0669)
\qbezier(166.4994,80.0669)(162.6319,76.993)(158.9077,73.6594)
\qbezier(158.9077,73.6594)(155.1834,70.3257)(151.7704,66.4926)
\qbezier(151.7704,66.4926)(148.1752,62.4403)(145.125,57.763)
\qbezier(145.125,57.763)(142.0748,53.0855)(139.7114,47.8038)
\qbezier(139.7114,47.8038)(137.7794,43.4747)(136.3926,38.8445)
\qbezier(136.3926,38.8445)(135.0057,34.2142)(134.187,29.3569)
\put(118.12,117.2469){\rule{0.6383bp}{0.143bp}}
\put(117.9025,117.3789){\rule{1.0733bp}{0.143bp}}
\put(117.7665,117.5109){\rule{1.3454bp}{0.143bp}}
\put(117.6677,117.6429){\rule{1.5429bp}{0.143bp}}
\put(117.5932,117.7749){\rule{1.6919bp}{0.143bp}}
\put(117.537,117.9069){\rule{1.8044bp}{0.143bp}}
\put(117.4957,118.0389){\rule{1.8868bp}{0.143bp}}
\put(117.4673,118.1709){\rule{1.9436bp}{0.143bp}}
\put(117.4519,118.3029){\rule{1.9745bp}{0.143bp}}
\put(117.4469,118.4349){\rule{1.9843bp}{0.143bp}}
\put(117.4536,118.5669){\rule{1.971bp}{0.143bp}}
\put(117.4715,118.6989){\rule{1.9352bp}{0.143bp}}
\put(117.5011,118.8309){\rule{1.8759bp}{0.143bp}}
\put(117.545,118.9629){\rule{1.7881bp}{0.143bp}}
\put(117.6042,119.0949){\rule{1.6699bp}{0.143bp}}
\put(117.6822,119.2269){\rule{1.5138bp}{0.143bp}}
\put(117.7855,119.3589){\rule{1.3073bp}{0.143bp}}
\put(117.9306,119.4909){\rule{1.017bp}{0.143bp}}
\put(118.1485,119.6229){\rule{0.5813bp}{0.1097bp}}
\linethickness{1bp}
\qbezier(119.4316,118.49)(119.421,119.018)(119.1409,119.3686)
\qbezier(119.1409,119.3686)(118.8608,119.7193)(118.4391,119.7326)
\qbezier(118.4391,119.7326)(118.0175,119.7193)(117.7374,119.3686)
\qbezier(117.7374,119.3686)(117.4573,119.018)(117.4467,118.49)
\qbezier(117.4467,118.49)(117.4573,117.9619)(117.7374,117.6111)
\qbezier(117.7374,117.6111)(118.0175,117.2602)(118.4391,117.2469)
\qbezier(118.4391,117.2469)(118.8608,117.2602)(119.1409,117.6111)
\qbezier(119.1409,117.6111)(119.421,117.9619)(119.4316,118.49)
\linethickness{0.8bp}
\qbezier(141.7839,88.5698)(141.7819,85.8253)(141.78,83.0808)
\qbezier(141.7772,79.2198)(141.7752,76.4753)(141.7732,73.7308)
\qbezier(141.7704,69.8698)(141.7685,67.1253)(141.7664,64.3808)
\qbezier(141.7637,60.5198)(141.7617,57.7753)(141.7598,55.0308)
\put(141.6405,51.17){\rule{0.6386bp}{0.143bp}}
\put(141.423,51.302){\rule{1.0736bp}{0.143bp}}
\put(141.287,51.434){\rule{1.3457bp}{0.143bp}}
\put(141.1883,51.566){\rule{1.5432bp}{0.143bp}}
\put(141.1137,51.698){\rule{1.6923bp}{0.143bp}}
\put(141.0575,51.83){\rule{1.8046bp}{0.143bp}}
\put(141.0163,51.962){\rule{1.8871bp}{0.143bp}}
\put(140.9879,52.094){\rule{1.9438bp}{0.143bp}}
\put(140.9726,52.226){\rule{1.9746bp}{0.143bp}}
\put(140.9677,52.358){\rule{1.9843bp}{0.143bp}}
\put(140.9744,52.49){\rule{1.9709bp}{0.143bp}}
\put(140.9923,52.622){\rule{1.9351bp}{0.143bp}}
\put(141.022,52.754){\rule{1.8758bp}{0.143bp}}
\put(141.0659,52.886){\rule{1.7879bp}{0.143bp}}
\put(141.1251,53.018){\rule{1.6695bp}{0.143bp}}
\put(141.2031,53.15){\rule{1.5135bp}{0.143bp}}
\put(141.3064,53.282){\rule{1.3068bp}{0.143bp}}
\put(141.4515,53.414){\rule{1.0167bp}{0.143bp}}
\put(141.6693,53.546){\rule{0.5811bp}{0.1097bp}}
\linethickness{1bp}
\qbezier(142.9523,52.4126)(142.9417,52.9407)(142.6616,53.2915)
\qbezier(142.6616,53.2915)(142.3815,53.6424)(141.9598,53.6557)
\qbezier(141.9598,53.6557)(141.5382,53.6424)(141.2581,53.2915)
\qbezier(141.2581,53.2915)(140.978,52.9407)(140.9674,52.4126)
\qbezier(140.9674,52.4126)(140.978,51.8846)(141.2581,51.534)
\qbezier(141.2581,51.534)(141.5382,51.1833)(141.9598,51.17)
\qbezier(141.9598,51.17)(142.3815,51.1833)(142.6616,51.534)
\qbezier(142.6616,51.534)(142.9417,51.8846)(142.9523,52.4126)
\linethickness{0.55bp}
\qbezier(89.6048,52.4126)(163.5307,52.4126)(237.4567,52.4126)
\linethickness{0.55bp}
\qbezier(89.6048,118.7426)(163.4297,118.7426)(237.2547,118.7426)
\linethickness{0.8bp}
\qbezier(118.2371,115.5529)(118.2371,112.8084)(118.2371,110.0639)
\qbezier(118.2371,106.2029)(118.2371,103.4584)(118.2371,100.7139)
\qbezier(118.2371,96.8529)(118.2371,94.1084)(118.2371,91.3639)
\linethickness{1.4bp}
\qbezier(179.5859,90.3188)(186.5859,93.8115)(192.8745,99.02)
\qbezier(192.8745,99.02)(196.1052,101.6874)(199.0479,104.8085)
\qbezier(199.0479,104.8085)(201.9904,107.9298)(204.5295,111.5707)
\qbezier(204.5295,111.5707)(207.0578,115.2295)(209.1317,119.3009)
\qbezier(209.1317,119.3009)(211.2056,123.3723)(212.9172,127.7078)
\qbezier(212.9172,127.7078)(215.394,133.9365)(217.212,140.5118)
\qbezier(217.212,140.5118)(218.8061,146.3115)(219.8733,152.2928)
\put(242,79){\makebox(0,0)[c]{$k$}}
\put(134.9423,103.0128){$j_{\nu,1}$}
\put(84,79){\makebox(0,0)[r]{$0$}}
\put(85,157){\makebox(0,0)[r]{$h_\nu(k)$}}
\put(174.7817,77.4048){$j_{\nu+1,1}$}
\put(226.3199,103.5408){$j_{\nu,2}$}
\put(203.8708,77.1408){$k_{\nu,2}$}
\put(144.9548,77.9328){$k_{\nu,2}$}
\put(85,47.1768){\makebox(0,0)[r]{$\sigma<0$}}
\put(85,113.7048){\makebox(0,0)[r]{$\sigma>0$}}
\put(109.7528,77.1408){$k_{\nu,1}$}
\end{picture}%
\endgroup%
}%
}
\hfill
\subfloat[{$\ell=0$, $\mu<0$}]{%

\resizebox{0.48\textwidth}{!}{%
\begingroup%
\setlength{\unitlength}{1bp}%
\fontsize{10.95}{13.14}\selectfont%
\begin{picture}(260,180)%
\PPWxaxis{68}{90}{244}%
\PPWyaxis{90}{16}{165}%
\linethickness{0.55bp}
\qbezier(89.4605,52.5338)(158.1287,52.5175)(226.797,52.5012)
\linethickness{0.55bp}
\qbezier(89.4605,118.8638)(157.9909,118.9135)(226.5212,118.9631)
\linethickness{1.4bp}
\qbezier(89.4122,89.9338)(104.2568,88.3144)(118.8316,85.3581)
\qbezier(118.8316,85.3581)(124.3362,84.2434)(129.783,82.9238)
\qbezier(129.783,82.9238)(135.2299,81.6041)(140.5929,80.0409)
\qbezier(140.5929,80.0409)(145.9561,78.4777)(151.2095,76.6321)
\qbezier(151.2095,76.6321)(156.4631,74.7866)(161.5812,72.62)
\qbezier(161.5812,72.62)(170.4761,68.8446)(178.7618,64.1331)
\qbezier(178.7618,64.1331)(187.0475,59.4215)(194.6386,53.8363)
\qbezier(194.6386,53.8363)(202.2298,48.2509)(209.041,41.8542)
\qbezier(209.041,41.8542)(215.8521,35.4576)(221.7979,28.3119)
\put(196.6666,51.005){\rule{0.8713bp}{0.143bp}}
\put(196.3697,51.137){\rule{1.465bp}{0.143bp}}
\put(196.1842,51.269){\rule{1.8362bp}{0.143bp}}
\put(196.0493,51.401){\rule{2.1058bp}{0.143bp}}
\put(195.9475,51.533){\rule{2.3094bp}{0.143bp}}
\put(195.8708,51.665){\rule{2.4627bp}{0.143bp}}
\put(195.8146,51.797){\rule{2.5754bp}{0.143bp}}
\put(195.7757,51.929){\rule{2.6529bp}{0.143bp}}
\put(195.7546,52.061){\rule{2.6951bp}{0.143bp}}
\put(195.7479,52.193){\rule{2.7085bp}{0.143bp}}
\put(195.757,52.325){\rule{2.6903bp}{0.143bp}}
\put(195.7815,52.457){\rule{2.6415bp}{0.143bp}}
\put(195.8219,52.589){\rule{2.5607bp}{0.143bp}}
\put(195.8818,52.721){\rule{2.4408bp}{0.143bp}}
\put(195.9625,52.853){\rule{2.2794bp}{0.143bp}}
\put(196.0688,52.985){\rule{2.0666bp}{0.143bp}}
\put(196.2097,53.117){\rule{1.7852bp}{0.143bp}}
\put(196.4074,53.249){\rule{1.3896bp}{0.143bp}}
\put(196.7047,53.381){\rule{0.7949bp}{0.1102bp}}
\linethickness{1bp}
\qbezier(198.4569,52.2481)(198.4424,52.7761)(198.0601,53.127)
\qbezier(198.0601,53.127)(197.6778,53.4779)(197.1022,53.4912)
\qbezier(197.1022,53.4912)(196.5267,53.4779)(196.1444,53.127)
\qbezier(196.1444,53.127)(195.7621,52.7761)(195.7476,52.2481)
\qbezier(195.7476,52.2481)(195.7621,51.7201)(196.1444,51.3692)
\qbezier(196.1444,51.3692)(196.5267,51.0184)(197.1022,51.005)
\qbezier(197.1022,51.005)(197.6778,51.0183)(198.0601,51.3692)
\qbezier(198.0601,51.3692)(198.4424,51.7201)(198.4569,52.2481)
\linethickness{0.8bp}
\qbezier(197.1022,89.8129)(197.1022,87.0684)(197.1022,84.3239)
\qbezier(197.1022,80.4629)(197.1022,77.7184)(197.1022,74.9739)
\qbezier(197.1022,71.1129)(197.1022,68.3684)(197.1022,65.6239)
\qbezier(197.1022,61.7629)(197.1022,59.0184)(197.1022,56.2739)
\put(196.6946,52.1435){\rule{0.8152bp}{0.748bp}}
\put(242,79){\makebox(0,0)[c]{$k$}}
\put(84,79){\makebox(0,0)[r]{$0$}}
\put(85,157){\makebox(0,0)[r]{$h_\nu(k)$}}
\put(85,47.298){\makebox(0,0)[r]{$\sigma<0$}}
\put(85,113.826){\makebox(0,0)[r]{$\sigma>0$}}
\put(184.7184,94.29){$k_{\nu,1}$}
\end{picture}%
\endgroup%
}%
}
\par\medskip
\subfloat[{$\ell\geq1$, $\mu>0$}]{%

\resizebox{0.48\textwidth}{!}{%
\begingroup%
\setlength{\unitlength}{1bp}%
\fontsize{10.95}{13.14}\selectfont%
\begin{picture}(260,180)%
\PPWxaxis{68}{90}{244}%
\PPWyaxis{90}{16}{165}%
\linethickness{0.8bp}
\qbezier(131.5232,144.2905)(131.5232,141.546)(131.5232,138.8015)
\qbezier(131.5232,134.9405)(131.5232,132.196)(131.5232,129.4515)
\qbezier(131.5232,125.5905)(131.5232,122.846)(131.5232,120.1015)
\qbezier(131.5232,116.2405)(131.5232,113.496)(131.5232,110.7515)
\qbezier(131.5232,106.8905)(131.5232,104.146)(131.5232,101.4015)
\qbezier(131.5232,97.5405)(131.5232,94.796)(131.5232,92.0515)
\qbezier(131.5232,88.1905)(131.5232,85.446)(131.5232,82.7015)
\qbezier(131.5232,78.8405)(131.5232,76.096)(131.5232,73.3515)
\qbezier(131.5232,69.4905)(131.5232,66.746)(131.5232,64.0015)
\qbezier(131.5232,60.1405)(131.5232,57.396)(131.5232,54.6515)
\qbezier(131.5232,50.7905)(131.5232,48.046)(131.5232,45.3015)
\qbezier(131.5232,41.4405)(131.5232,38.696)(131.5232,35.9515)
\qbezier(131.5232,32.0905)(131.5232,29.346)(131.5232,26.6015)
\qbezier(131.5232,22.7405)(131.5232,20.1941)(131.5232,17.6477)
\linethickness{0.8bp}
\qbezier(222.9485,143.2675)(222.9485,140.523)(222.9485,137.7785)
\qbezier(222.9485,133.9175)(222.9485,131.173)(222.9485,128.4285)
\qbezier(222.9485,124.5675)(222.9485,121.823)(222.9485,119.0785)
\qbezier(222.9485,115.2175)(222.9485,112.473)(222.9485,109.7285)
\qbezier(222.9485,105.8675)(222.9485,103.123)(222.9485,100.3785)
\qbezier(222.9485,96.5175)(222.9485,93.773)(222.9485,91.0285)
\qbezier(222.9485,87.1675)(222.9485,84.423)(222.9485,81.6785)
\qbezier(222.9485,77.8175)(222.9485,75.073)(222.9485,72.3285)
\qbezier(222.9485,68.4675)(222.9485,65.723)(222.9485,62.9785)
\qbezier(222.9485,59.1175)(222.9485,56.373)(222.9485,53.6285)
\qbezier(222.9485,49.7675)(222.9485,47.023)(222.9485,44.2785)
\qbezier(222.9485,40.4175)(222.9485,37.673)(222.9485,34.9285)
\qbezier(222.9485,31.0675)(222.9485,28.323)(222.9485,25.5785)
\qbezier(222.9485,21.7175)(222.9485,19.171)(222.9485,16.6246)
\linethickness{1.4bp}
\qbezier(90.1842,55.1794)(92.8344,60.8613)(95.6822,66.3775)
\qbezier(95.6822,66.3775)(100.2424,75.2299)(105.2375,83.7025)
\qbezier(105.2375,83.7025)(107.3883,87.3615)(109.5836,90.9775)
\qbezier(109.5836,90.9775)(111.7789,94.5935)(113.9759,98.2113)
\qbezier(113.9759,98.2113)(116.0895,101.6666)(118.1091,105.2046)
\qbezier(118.1091,105.2046)(120.1288,108.7426)(121.9064,112.4787)
\qbezier(121.9064,112.4787)(123.3336,115.5067)(124.517,118.6894)
\qbezier(124.517,118.6894)(125.7005,121.8722)(126.5875,125.2055)
\qbezier(126.5875,125.2055)(127.6688,129.3449)(128.2298,133.6288)
\qbezier(128.2298,133.6288)(128.7907,137.9126)(128.8182,142.2665)
\put(150.2532,53.9363){\rule{0.6383bp}{0.143bp}}
\put(150.0357,54.0683){\rule{1.0732bp}{0.143bp}}
\put(149.8998,54.2003){\rule{1.3452bp}{0.143bp}}
\put(149.801,54.3323){\rule{1.5427bp}{0.143bp}}
\put(149.7264,54.4643){\rule{1.6919bp}{0.143bp}}
\put(149.6702,54.5963){\rule{1.8042bp}{0.143bp}}
\put(149.629,54.7283){\rule{1.8867bp}{0.143bp}}
\put(149.6006,54.8603){\rule{1.9435bp}{0.143bp}}
\put(149.5851,54.9923){\rule{1.9744bp}{0.143bp}}
\put(149.5802,55.1243){\rule{1.9842bp}{0.143bp}}
\put(149.5869,55.2563){\rule{1.9709bp}{0.143bp}}
\put(149.6048,55.3883){\rule{1.9351bp}{0.143bp}}
\put(149.6344,55.5203){\rule{1.8759bp}{0.143bp}}
\put(149.6783,55.6523){\rule{1.7882bp}{0.143bp}}
\put(149.7373,55.7843){\rule{1.67bp}{0.143bp}}
\put(149.8153,55.9163){\rule{1.5141bp}{0.143bp}}
\put(149.9184,56.0483){\rule{1.3078bp}{0.143bp}}
\put(150.0633,56.1803){\rule{1.0181bp}{0.143bp}}
\put(150.2811,56.3123){\rule{0.5825bp}{0.1102bp}}
\linethickness{1bp}
\qbezier(151.5648,55.1794)(151.5541,55.7075)(151.2741,56.0583)
\qbezier(151.2741,56.0583)(150.994,56.4092)(150.5723,56.4225)
\qbezier(150.5723,56.4225)(150.1507,56.4092)(149.8706,56.0583)
\qbezier(149.8706,56.0583)(149.5906,55.7075)(149.5799,55.1794)
\qbezier(149.5799,55.1794)(149.5906,54.6514)(149.8706,54.3005)
\qbezier(149.8706,54.3005)(150.1507,53.9497)(150.5723,53.9363)
\qbezier(150.5723,53.9363)(150.994,53.9497)(151.2741,54.3005)
\qbezier(151.2741,54.3005)(151.5541,54.6514)(151.5648,55.1794)
\put(208.9199,107.5613){\rule{0.6383bp}{0.143bp}}
\put(208.7025,107.6933){\rule{1.0732bp}{0.143bp}}
\put(208.5665,107.8253){\rule{1.3452bp}{0.143bp}}
\put(208.4677,107.9573){\rule{1.5427bp}{0.143bp}}
\put(208.3932,108.0893){\rule{1.6919bp}{0.143bp}}
\put(208.3369,108.2213){\rule{1.8042bp}{0.143bp}}
\put(208.2958,108.3533){\rule{1.8867bp}{0.143bp}}
\put(208.2673,108.4853){\rule{1.9435bp}{0.143bp}}
\put(208.2518,108.6173){\rule{1.9744bp}{0.143bp}}
\put(208.2469,108.7493){\rule{1.9842bp}{0.143bp}}
\put(208.2536,108.8813){\rule{1.9709bp}{0.143bp}}
\put(208.2715,109.0133){\rule{1.9351bp}{0.143bp}}
\put(208.3011,109.1453){\rule{1.8759bp}{0.143bp}}
\put(208.345,109.2773){\rule{1.7882bp}{0.143bp}}
\put(208.404,109.4093){\rule{1.67bp}{0.143bp}}
\put(208.482,109.5413){\rule{1.5141bp}{0.143bp}}
\put(208.5851,109.6733){\rule{1.3078bp}{0.143bp}}
\put(208.7301,109.8053){\rule{1.0181bp}{0.143bp}}
\put(208.9479,109.9373){\rule{0.5825bp}{0.1102bp}}
\linethickness{1bp}
\qbezier(210.2315,108.8044)(210.2209,109.3325)(209.9408,109.6833)
\qbezier(209.9408,109.6833)(209.6607,110.0342)(209.2391,110.0475)
\qbezier(209.2391,110.0475)(208.8174,110.0342)(208.5374,109.6833)
\qbezier(208.5374,109.6833)(208.2573,109.3325)(208.2467,108.8044)
\qbezier(208.2467,108.8044)(208.2573,108.2764)(208.5374,107.9255)
\qbezier(208.5374,107.9255)(208.8174,107.5746)(209.2391,107.5613)
\qbezier(209.2391,107.5613)(209.6607,107.5746)(209.9408,107.9255)
\qbezier(209.9408,107.9255)(210.2209,108.2764)(210.2315,108.8044)
\linethickness{1.4bp}
\qbezier(180.4678,79.8744)(173.3167,75.2382)(166.4949,69.8756)
\qbezier(166.4949,69.8756)(162.6043,66.8459)(158.8915,63.4967)
\qbezier(158.8915,63.4967)(155.1787,60.1475)(151.7755,56.3013)
\qbezier(151.7755,56.3013)(148.1874,52.2347)(145.136,47.5634)
\qbezier(145.136,47.5634)(142.0846,42.8922)(139.7085,37.6125)
\qbezier(139.7085,37.6125)(137.7763,33.2835)(136.3909,28.6532)
\qbezier(136.3909,28.6532)(135.0055,24.0229)(134.1931,19.1656)
\put(119.7605,107.4075){\rule{0.6383bp}{0.143bp}}
\put(119.5431,107.5395){\rule{1.0732bp}{0.143bp}}
\put(119.4071,107.6715){\rule{1.3452bp}{0.143bp}}
\put(119.3083,107.8035){\rule{1.5427bp}{0.143bp}}
\put(119.2338,107.9355){\rule{1.6919bp}{0.143bp}}
\put(119.1776,108.0675){\rule{1.8042bp}{0.143bp}}
\put(119.1364,108.1995){\rule{1.8867bp}{0.143bp}}
\put(119.1079,108.3315){\rule{1.9435bp}{0.143bp}}
\put(119.0925,108.4635){\rule{1.9744bp}{0.143bp}}
\put(119.0875,108.5955){\rule{1.9842bp}{0.143bp}}
\put(119.0942,108.7275){\rule{1.9709bp}{0.143bp}}
\put(119.1121,108.8595){\rule{1.9351bp}{0.143bp}}
\put(119.1417,108.9915){\rule{1.8759bp}{0.143bp}}
\put(119.1856,109.1235){\rule{1.7881bp}{0.143bp}}
\put(119.2448,109.2555){\rule{1.6699bp}{0.143bp}}
\put(119.3227,109.3875){\rule{1.514bp}{0.143bp}}
\put(119.4258,109.5195){\rule{1.3078bp}{0.143bp}}
\put(119.5707,109.6515){\rule{1.018bp}{0.143bp}}
\put(119.7885,109.7835){\rule{0.5824bp}{0.1102bp}}
\linethickness{1bp}
\qbezier(121.0721,108.6506)(121.0615,109.1786)(120.7814,109.5295)
\qbezier(120.7814,109.5295)(120.5013,109.8803)(120.0797,109.8937)
\qbezier(120.0797,109.8937)(119.6581,109.8803)(119.378,109.5295)
\qbezier(119.378,109.5295)(119.0979,109.1786)(119.0873,108.6506)
\qbezier(119.0873,108.6506)(119.0979,108.1226)(119.378,107.7717)
\qbezier(119.378,107.7717)(119.6581,107.4208)(120.0797,107.4075)
\qbezier(120.0797,107.4075)(120.5013,107.4208)(120.7814,107.7717)
\qbezier(120.7814,107.7717)(121.0615,108.1226)(121.0721,108.6506)
\put(141.6462,40.9787){\rule{0.6386bp}{0.143bp}}
\put(141.4288,41.1107){\rule{1.0736bp}{0.143bp}}
\put(141.2927,41.2427){\rule{1.3456bp}{0.143bp}}
\put(141.194,41.3747){\rule{1.5432bp}{0.143bp}}
\put(141.1195,41.5067){\rule{1.6922bp}{0.143bp}}
\put(141.0632,41.6387){\rule{1.8045bp}{0.143bp}}
\put(141.0221,41.7707){\rule{1.887bp}{0.143bp}}
\put(140.9937,41.9027){\rule{1.9436bp}{0.143bp}}
\put(140.9783,42.0347){\rule{1.9745bp}{0.143bp}}
\put(140.9735,42.1667){\rule{1.9842bp}{0.143bp}}
\put(140.9802,42.2987){\rule{1.9708bp}{0.143bp}}
\put(140.9981,42.4307){\rule{1.9349bp}{0.143bp}}
\put(141.0278,42.5627){\rule{1.8757bp}{0.143bp}}
\put(141.0717,42.6947){\rule{1.7878bp}{0.143bp}}
\put(141.1309,42.8267){\rule{1.6695bp}{0.143bp}}
\put(141.2089,42.9587){\rule{1.5134bp}{0.143bp}}
\put(141.3121,43.0907){\rule{1.3067bp}{0.143bp}}
\put(141.4572,43.2227){\rule{1.0167bp}{0.143bp}}
\put(141.675,43.3547){\rule{0.581bp}{0.1097bp}}
\linethickness{1bp}
\qbezier(142.958,42.2213)(142.9473,42.7494)(142.6673,43.1002)
\qbezier(142.6673,43.1002)(142.3872,43.4511)(141.9656,43.4644)
\qbezier(141.9656,43.4644)(141.5439,43.4511)(141.2638,43.1002)
\qbezier(141.2638,43.1002)(140.9838,42.7494)(140.9731,42.2213)
\qbezier(140.9731,42.2213)(140.9838,41.6933)(141.2638,41.3427)
\qbezier(141.2638,41.3427)(141.5439,40.992)(141.9656,40.9787)
\qbezier(141.9656,40.9787)(142.3872,40.992)(142.6673,41.3427)
\qbezier(142.6673,41.3427)(142.9473,41.6933)(142.958,42.2213)
\linethickness{0.55bp}
\qbezier(89.6133,42.2213)(163.5352,42.2213)(237.4571,42.2213)
\linethickness{0.55bp}
\qbezier(89.6133,108.5513)(163.4298,108.5513)(237.2463,108.5513)
\linethickness{0.8bp}
\qbezier(120.2731,109.6954)(120.2743,106.9509)(120.2755,104.2064)
\qbezier(120.2772,100.3454)(120.2784,97.6009)(120.2796,94.8564)
\qbezier(120.2813,90.9954)(120.2814,90.715)(120.2815,90.4346)
\linethickness{1.4bp}
\qbezier(179.651,79.3687)(186.3696,83.763)(192.6666,89.0487)
\qbezier(192.6666,89.0487)(195.8914,91.7292)(198.9087,94.7604)
\qbezier(198.9087,94.7604)(201.926,97.7914)(204.5317,101.3794)
\qbezier(204.5317,101.3794)(207.1011,104.9941)(209.16,109.0808)
\qbezier(209.16,109.0808)(211.219,113.1675)(212.9189,117.5165)
\qbezier(212.9189,117.5165)(215.378,123.7535)(217.2048,130.3205)
\qbezier(217.2048,130.3205)(218.8119,136.1202)(219.8658,142.1015)
\linethickness{0.55bp}
\qbezier(90.4302,55.0694)(164.2466,55.0694)(238.0631,55.0694)
\linethickness{0.55bp}
\qbezier(89.8154,72.2737)(163.7372,72.2737)(237.6591,72.2737)
\linethickness{0.8bp}
\qbezier(209.3271,108.6723)(209.3283,105.9278)(209.3295,103.1833)
\qbezier(209.3313,99.3223)(209.3325,96.5778)(209.3338,93.8333)
\put(209.0369,89.5764){\rule{0.5972bp}{0.748bp}}
\put(98.6126,70.9975){\rule{0.6381bp}{0.143bp}}
\put(98.3951,71.1295){\rule{1.073bp}{0.143bp}}
\put(98.2591,71.2615){\rule{1.3451bp}{0.143bp}}
\put(98.1603,71.3935){\rule{1.5427bp}{0.143bp}}
\put(98.0857,71.5255){\rule{1.6919bp}{0.143bp}}
\put(98.0295,71.6575){\rule{1.8043bp}{0.143bp}}
\put(97.9882,71.7895){\rule{1.8868bp}{0.143bp}}
\put(97.9599,71.9215){\rule{1.9436bp}{0.143bp}}
\put(97.9443,72.0535){\rule{1.9746bp}{0.143bp}}
\put(97.9394,72.1855){\rule{1.9844bp}{0.143bp}}
\put(97.9461,72.3175){\rule{1.971bp}{0.143bp}}
\put(97.964,72.4495){\rule{1.9353bp}{0.143bp}}
\put(97.9936,72.5815){\rule{1.876bp}{0.143bp}}
\put(98.0375,72.7135){\rule{1.7882bp}{0.143bp}}
\put(98.0967,72.8455){\rule{1.67bp}{0.143bp}}
\put(98.1747,72.9775){\rule{1.514bp}{0.143bp}}
\put(98.2778,73.1095){\rule{1.3077bp}{0.143bp}}
\put(98.4227,73.2415){\rule{1.0179bp}{0.143bp}}
\put(98.6405,73.3735){\rule{0.5823bp}{0.1102bp}}
\linethickness{1bp}
\qbezier(99.9242,72.2406)(99.9135,72.7686)(99.6334,73.1195)
\qbezier(99.6334,73.1195)(99.3532,73.4704)(98.9317,73.4837)
\qbezier(98.9317,73.4837)(98.5101,73.4704)(98.2299,73.1195)
\qbezier(98.2299,73.1195)(97.9499,72.7686)(97.9392,72.2406)
\qbezier(97.9392,72.2406)(97.9499,71.7126)(98.2299,71.3617)
\qbezier(98.2299,71.3617)(98.5101,71.0109)(98.9317,70.9975)
\qbezier(98.9317,70.9975)(99.3532,71.0108)(99.6334,71.3617)
\qbezier(99.6334,71.3617)(99.9135,71.7126)(99.9242,72.2406)
\linethickness{0.8bp}
\qbezier(150.5723,89.9504)(150.5723,87.2059)(150.5723,84.4614)
\qbezier(150.5723,80.6004)(150.5723,77.8559)(150.5723,75.1114)
\qbezier(150.5723,71.2504)(150.5723,68.5059)(150.5723,65.7614)
\qbezier(150.5723,61.9004)(150.5723,59.365)(150.5723,56.8296)
\put(168.9598,71.2506){\rule{0.6383bp}{0.143bp}}
\put(168.7423,71.3826){\rule{1.0732bp}{0.143bp}}
\put(168.6063,71.5146){\rule{1.3453bp}{0.143bp}}
\put(168.5076,71.6466){\rule{1.5428bp}{0.143bp}}
\put(168.433,71.7786){\rule{1.6919bp}{0.143bp}}
\put(168.3768,71.9106){\rule{1.8043bp}{0.143bp}}
\put(168.3355,72.0426){\rule{1.8867bp}{0.143bp}}
\put(168.3072,72.1746){\rule{1.9435bp}{0.143bp}}
\put(168.2917,72.3066){\rule{1.9744bp}{0.143bp}}
\put(168.2868,72.4386){\rule{1.9842bp}{0.143bp}}
\put(168.2935,72.5706){\rule{1.9709bp}{0.143bp}}
\put(168.3114,72.7026){\rule{1.9351bp}{0.143bp}}
\put(168.341,72.8346){\rule{1.8758bp}{0.143bp}}
\put(168.3849,72.9666){\rule{1.788bp}{0.143bp}}
\put(168.4441,73.0986){\rule{1.6698bp}{0.143bp}}
\put(168.5221,73.2306){\rule{1.5137bp}{0.143bp}}
\put(168.6254,73.3626){\rule{1.3072bp}{0.143bp}}
\put(168.7705,73.4946){\rule{1.0169bp}{0.143bp}}
\put(168.9883,73.6266){\rule{0.5812bp}{0.1097bp}}
\linethickness{1bp}
\qbezier(170.2714,72.4937)(170.2607,73.0217)(169.9807,73.3723)
\qbezier(169.9807,73.3723)(169.7006,73.723)(169.279,73.7363)
\qbezier(169.279,73.7363)(168.8573,73.723)(168.5772,73.3723)
\qbezier(168.5772,73.3723)(168.2972,73.0217)(168.2865,72.4937)
\qbezier(168.2865,72.4937)(168.2972,71.9656)(168.5772,71.6148)
\qbezier(168.5772,71.6148)(168.8573,71.2639)(169.279,71.2506)
\qbezier(169.279,71.2506)(169.7006,71.2639)(169.9807,71.6148)
\qbezier(169.9807,71.6148)(170.2607,71.9656)(170.2714,72.4937)
\put(89.663,53.9363){\rule{0.6381bp}{0.143bp}}
\put(89.4457,54.0683){\rule{1.073bp}{0.143bp}}
\put(89.3097,54.2003){\rule{1.3449bp}{0.143bp}}
\put(89.211,54.3323){\rule{1.5425bp}{0.143bp}}
\put(89.1364,54.4643){\rule{1.6916bp}{0.143bp}}
\put(89.0803,54.5963){\rule{1.8039bp}{0.143bp}}
\put(89.039,54.7283){\rule{1.8865bp}{0.143bp}}
\put(89.0107,54.8603){\rule{1.9432bp}{0.143bp}}
\put(88.9952,54.9923){\rule{1.9741bp}{0.143bp}}
\put(88.9903,55.1243){\rule{1.984bp}{0.143bp}}
\put(88.997,55.2563){\rule{1.9706bp}{0.143bp}}
\put(89.0149,55.3883){\rule{1.9349bp}{0.143bp}}
\put(89.0445,55.5203){\rule{1.8757bp}{0.143bp}}
\put(89.0883,55.6523){\rule{1.7879bp}{0.143bp}}
\put(89.1474,55.7843){\rule{1.6697bp}{0.143bp}}
\put(89.2253,55.9163){\rule{1.5138bp}{0.143bp}}
\put(89.3284,56.0483){\rule{1.3076bp}{0.143bp}}
\put(89.4732,56.1803){\rule{1.0179bp}{0.143bp}}
\put(89.6909,56.3123){\rule{0.5824bp}{0.1102bp}}
\linethickness{1bp}
\qbezier(90.9746,55.1794)(90.964,55.7075)(90.6838,56.0583)
\qbezier(90.6838,56.0583)(90.4037,56.4092)(89.9821,56.4225)
\qbezier(89.9821,56.4225)(89.5605,56.4092)(89.2806,56.0583)
\qbezier(89.2806,56.0583)(89.0006,55.7075)(88.99,55.1794)
\qbezier(88.99,55.1794)(89.0006,54.6514)(89.2806,54.3005)
\qbezier(89.2806,54.3005)(89.5605,53.9497)(89.9821,53.9363)
\qbezier(89.9821,53.9363)(90.4037,53.9497)(90.6838,54.3005)
\qbezier(90.6838,54.3005)(90.964,54.6514)(90.9746,55.1794)
\put(242,79){\makebox(0,0)[c]{$k$}}
\put(134.7815,118.2975){$j_{\nu+\ell,1}$}
\put(84,79){\makebox(0,0)[r]{$0$}}
\put(85,157){\makebox(0,0)[r]{$h_{\nu+\ell}(k)$}}
\put(140.7887,96.6495){$j_{\nu+\ell+1,1}$}
\put(226.1541,118.8255){$j_{\nu+\ell,2}$}
\put(85,36.9855){\makebox(0,0)[r]{$\sigma<-\ell$}}
\put(85,103.5135){\makebox(0,0)[r]{$\sigma>0$}}
\put(107.5911,78.5655){$k_{\nu+\ell,1}$}
\put(85,52.8255){\makebox(0,0)[r]{$-\ell$}}
\put(197.3828,78.5655){$k_{\nu+\ell,2}$}
\put(85,66.1575){\makebox(0,0)[r]{$0>\sigma>-\ell$}}
\end{picture}%
\endgroup%
}%
}
\hfill
\subfloat[{$\ell\geq1$, $\mu<0$}]{%

\resizebox{0.48\textwidth}{!}{%
\begingroup%
\setlength{\unitlength}{1bp}%
\fontsize{10.95}{13.14}\selectfont%
\begin{picture}(260,180)%
\PPWxaxis{68}{90}{244}%
\PPWyaxis{90}{16}{165}%
\linethickness{0.55bp}
\qbezier(89.4728,32.8988)(158.1356,32.8825)(226.7983,32.8662)
\linethickness{0.55bp}
\qbezier(89.4728,109.052)(157.9977,109.1015)(226.5226,109.151)
\put(185.9657,31.37){\rule{0.8712bp}{0.143bp}}
\put(185.6689,31.502){\rule{1.4648bp}{0.143bp}}
\put(185.4834,31.634){\rule{1.8361bp}{0.143bp}}
\put(185.3485,31.766){\rule{2.1057bp}{0.143bp}}
\put(185.2467,31.898){\rule{2.3092bp}{0.143bp}}
\put(185.17,32.03){\rule{2.4625bp}{0.143bp}}
\put(185.1138,32.162){\rule{2.5752bp}{0.143bp}}
\put(185.075,32.294){\rule{2.6527bp}{0.143bp}}
\put(185.0539,32.426){\rule{2.6949bp}{0.143bp}}
\put(185.0471,32.558){\rule{2.7083bp}{0.143bp}}
\put(185.0563,32.69){\rule{2.6901bp}{0.143bp}}
\put(185.0807,32.822){\rule{2.6413bp}{0.143bp}}
\put(185.1211,32.954){\rule{2.5605bp}{0.143bp}}
\put(185.181,33.086){\rule{2.4406bp}{0.143bp}}
\put(185.2617,33.218){\rule{2.2793bp}{0.143bp}}
\put(185.368,33.35){\rule{2.0665bp}{0.143bp}}
\put(185.5089,33.482){\rule{1.785bp}{0.143bp}}
\put(185.7066,33.614){\rule{1.3894bp}{0.143bp}}
\put(186.0038,33.746){\rule{0.7949bp}{0.1102bp}}
\linethickness{1bp}
\qbezier(187.7559,32.6131)(187.7414,33.1411)(187.3591,33.492)
\qbezier(187.3591,33.492)(186.9768,33.8429)(186.4013,33.8562)
\qbezier(186.4013,33.8562)(185.8258,33.8429)(185.4436,33.492)
\qbezier(185.4436,33.492)(185.0613,33.1411)(185.0468,32.6131)
\qbezier(185.0468,32.6131)(185.0613,32.0851)(185.4436,31.7342)
\qbezier(185.4436,31.7342)(185.8258,31.3834)(186.4013,31.37)
\qbezier(186.4013,31.37)(186.9768,31.3833)(187.3591,31.7342)
\qbezier(187.3591,31.7342)(187.7414,32.0851)(187.7559,32.6131)
\linethickness{0.8bp}
\qbezier(186.5324,89.6372)(186.526,86.8927)(186.5198,84.1482)
\qbezier(186.5109,80.2872)(186.5045,77.5427)(186.4982,74.7982)
\qbezier(186.4893,70.9372)(186.4831,68.1927)(186.4767,65.4482)
\qbezier(186.4679,61.5873)(186.4616,58.8428)(186.4553,56.0983)
\qbezier(186.4464,52.2373)(186.4401,49.4928)(186.4338,46.7483)
\qbezier(186.425,42.8873)(186.4186,40.1428)(186.4124,37.3983)
\qbezier(186.4035,33.5373)(186.4028,33.2622)(186.4022,32.9871)
\linethickness{0.55bp}
\qbezier(89.4728,69.7819)(157.9977,69.8315)(226.5226,69.8812)
\linethickness{1.4bp}
\qbezier(89.7003,51.4781)(104.6006,51.3461)(119.3927,49.7619)
\qbezier(119.3927,49.7619)(130.0827,48.6221)(140.6087,46.6407)
\qbezier(140.6087,46.6407)(151.1345,44.6595)(161.384,41.655)
\qbezier(161.384,41.655)(170.2266,39.0536)(178.7132,35.6698)
\qbezier(178.7132,35.6698)(187.1998,32.286)(195.2705,28.1512)
\qbezier(195.2705,28.1512)(203.341,24.0164)(210.9353,19.162)
\qbezier(210.9353,19.162)(218.5296,14.3076)(225.5876,8.765)
\put(89.2648,49.9819){\rule{0.871bp}{0.143bp}}
\put(88.968,50.1139){\rule{1.4646bp}{0.143bp}}
\put(88.7823,50.2459){\rule{1.836bp}{0.143bp}}
\put(88.6475,50.3779){\rule{2.1057bp}{0.143bp}}
\put(88.5457,50.5099){\rule{2.3092bp}{0.143bp}}
\put(88.469,50.6419){\rule{2.4627bp}{0.143bp}}
\put(88.4127,50.7739){\rule{2.5753bp}{0.143bp}}
\put(88.3739,50.9059){\rule{2.6529bp}{0.143bp}}
\put(88.3527,51.0379){\rule{2.6951bp}{0.143bp}}
\put(88.346,51.1699){\rule{2.7085bp}{0.143bp}}
\put(88.3551,51.3019){\rule{2.6903bp}{0.143bp}}
\put(88.3796,51.4339){\rule{2.6415bp}{0.143bp}}
\put(88.42,51.5659){\rule{2.5606bp}{0.143bp}}
\put(88.4799,51.6979){\rule{2.4407bp}{0.143bp}}
\put(88.5607,51.8299){\rule{2.2794bp}{0.143bp}}
\put(88.6671,51.9619){\rule{2.0665bp}{0.143bp}}
\put(88.8079,52.0939){\rule{1.7849bp}{0.143bp}}
\put(89.0057,52.2259){\rule{1.3893bp}{0.143bp}}
\put(89.3029,52.3579){\rule{0.7948bp}{0.1102bp}}
\linethickness{1bp}
\qbezier(91.055,51.225)(91.0404,51.753)(90.6581,52.1039)
\qbezier(90.6581,52.1039)(90.2757,52.4547)(89.7003,52.4681)
\qbezier(89.7003,52.4681)(89.1249,52.4547)(88.7425,52.1039)
\qbezier(88.7425,52.1039)(88.3603,51.753)(88.3456,51.225)
\qbezier(88.3456,51.225)(88.3603,50.697)(88.7425,50.3461)
\qbezier(88.7425,50.3461)(89.1249,49.9953)(89.7003,49.9819)
\qbezier(89.7003,49.9819)(90.2757,49.9953)(90.6581,50.3461)
\qbezier(90.6581,50.3461)(91.0404,50.697)(91.055,51.225)
\linethickness{0.55bp}
\qbezier(89.4728,51.61)(157.9977,51.665)(226.5226,51.72)
\put(242,79){\makebox(0,0)[c]{$k$}}
\put(84,79){\makebox(0,0)[r]{$0$}}
\put(85,157){\makebox(0,0)[r]{$h_{\nu+\ell}(k)$}}
\put(85,27.663){\makebox(0,0)[r]{$\sigma<-\ell$}}
\put(85,104.091){\makebox(0,0)[r]{$\sigma>0$}}
\put(173.563,94.191){$k_{\nu+\ell,1}$}
\put(85,46.275){\makebox(0,0)[r]{$-\ell$}}
\put(85,63.963){\makebox(0,0)[r]{$0>\sigma>-\ell$}}
\end{picture}%
\endgroup%
}%
}
\caption{{The equations for the positive and negative
radial eigenvalues. When $\sigma+\ell\leq0$, the first positive root
corresponds to $m=2$. In the negative-eigenvalue panels,
$k_{\nu+\ell,1}$ denotes $\widehat{k}_{\nu+\ell,1}$.}}
\label{fig:appendix-ball-roots}
\end{figure}

The regular radial factors above, multiplied by $G_\ell$, give
the eigenfunctions in each spherical harmonic space. For a fixed
$\ell$, the radial eigenvalues are simple. Their products with a
basis of spherical harmonics span a dense subspace of $L^2(B)$;
see \cite{Chavel}. The dimension of the degree-$\ell$ space is
\[
 d_{N,\ell}=\frac{(2\ell+N-2)(\ell+N-3)!}{\ell!(N-2)!}
 \quad(N\geq3),\qquad
 d_{2,0}=1,\quad d_{2,\ell}=2\quad(\ell\geq1).
\]
If radial eigenvalues from different degrees coincide, their
eigenspaces are added. Thus the multiplicity of a full eigenvalue
$\mu$ is $\sum_{\mu_{\ell,m}=\mu}d_{N,\ell}$; see
\cite[Remarks~3.3 and~3.4]{HassannezhadSher}.

These radial branches identify the first two eigenvalues used in
the main text; see \cite[Section~5]{Freitas}.
\medskip
\begin{samepage}
\begin{proposition}\leavevmode\label{prop:ball-first-two}
For $N\geq2$ and every real $\sigma$,
\[
 \mu_1(B,\sigma)=\mu_{0,1}(\sigma),\,\,\,
 \mu_2(B,\sigma)=\mu_{1,1}(\sigma).
\]
The first eigenvalue is simple and the second has multiplicity $N$.
Their signs are those of $\sigma$ and $\sigma+1$, respectively.
In particular, for $\sigma>0$,
\begin{equation}\label{eq:appendix-positive-first-two}
 \mu_1(B,\sigma)=k_{\nu,1}(\sigma)^2,\,\,\,
 \mu_2(B,\sigma)=k_{\nu+1,1}(\sigma)^2.
\end{equation}
\end{proposition}
\end{samepage}
\medskip
\noindent\textbf{Proof.}
The first radial functions $v_{\ell,1}$ are positive on $(0,1]$ for
every real $\sigma$. The Wronskian calculation used in
\eqref{eq:angular-gap-lower-bound} gives, for $\ell>i\geq0$,
\[
 \mu_{\ell,1}-\mu_{i,1}
 =(c_\ell-c_i)
 \frac{\displaystyle\int_0^1r^{N-3}v_{i,1}v_{\ell,1}\,\mathrm{d}r}
      {\displaystyle\int_0^1r^{N-1}v_{i,1}v_{\ell,1}\,\mathrm{d}r}
 >c_\ell-c_i>0.
\]
The boundary terms vanish by the common Robin condition and
regularity at the origin. The integrals converge because
$v_{j,1}(r)$ has order $r^j$ at zero. For each degree the radial
eigenvalues are strictly ordered, so the first full eigenvalue is
$\mu_{0,1}$, and the second is the smaller of $\mu_{1,1}$ and
$\mu_{0,2}$.

When $\sigma>0$, the root locations above give
\[
 k_{\nu+1,1}<j_{\nu+1,1}<k_{\nu,2}.
\]
When $-1<\sigma\leq0$, the Bessel recurrence relation gives
$\widetilde h_{\nu+1}(j_{\nu,1})=N-1>\sigma$. Strict increase of
$\widetilde h_{\nu+1}$ before its first pole then gives
\[
 k_{\nu+1,1}<j_{\nu,1}<k_{\nu,2}.
\]
For $\sigma\leq-1$, one has $\mu_{1,1}\leq0<\mu_{0,2}$.
Thus $\mu_{1,1}<\mu_{0,2}$ for every real $\sigma$.
The multiplicities are $d_{N,0}=1$ and $d_{N,1}=N$, and the sign
statements follow from the thresholds $\sigma=0$ and $\sigma=-1$.
\qed\\

Each fixed radial eigenvalue is simple within its angular degree.
The parameter-derivative formula \eqref{eq:angular-hf}, applied to
its radial eigenfunction $v_{\ell,m}$, gives
\[
 \frac{\mathrm{d}\mu_{\ell,m}}{\mathrm{d}\sigma}
 =\frac{|v_{\ell,m}(1)|^2}
 {\displaystyle\int_0^1r^{N-1}|v_{\ell,m}(r)|^2\,\mathrm{d}r}>0.
\]
The boundary value cannot vanish, since the Robin condition would
then give $v_{\ell,m}'(1)=0$ as well. Thus each radial eigenvalue
is strictly increasing in $\sigma$. The positive roots approach their next
Bessel pole as $\sigma\to+\infty$, giving
\[
 \lim_{\sigma\to+\infty}\mu_{\ell,m}(\sigma)=j_{\nu+\ell,m}^2.
\]
Every degree $\ell<-\sigma$ contributes one negative radial value,
with the multiplicity of its spherical harmonic space. More
precisely, if $-\ell\leq\sigma<1-\ell$ for an integer $\ell\geq1$,
the negative eigenvalues are $\mu_{0,1},\ldots,\mu_{\ell-1,1}$,
in strictly increasing order. Their total multiplicity is
\[
 \sum_{i=0}^{\ell-1}d_{N,i}
 =\binom{N+\ell-2}{N-1}+\binom{N+\ell-3}{N-1}.
\]
A binomial coefficient is zero when its upper nonnegative index
is smaller than its lower index. The count is $2\ell-1$ for $N=2$
and $\ell^2$ for $N=3$. At $\sigma=-\ell$, zero is also an eigenvalue,
with eigenfunctions $r^\ell G_\ell(\xi)$ and multiplicity $d_{N,\ell}$.
For $\sigma\geq0$, there are no negative eigenvalues.
Thus their number is finite for fixed $\sigma$, but increases
without bound as $\sigma\to-\infty$. For each fixed
$\ell$, $\mu_{\ell,1}\to-\infty$ in that limit. The scaling
law \eqref{eq:robin-scaling} gives the formulas on balls of any radius.

Table~\ref{tab:appendix-positive-spectrum} gives positive-parameter
values for $N=2,3$.
Values are rounded to five decimal places, with ratios calculated
before rounding.

\begin{table}[htbp]

\centering
\small
\setlength{\tabcolsep}{4pt}
\begin{tabular}{c|c|c|c|c|c|c|c|c}
\hline
$\ell$ & $\nu$ & $\sigma=1$ & $\sigma=2$ & $\sigma=3$ & $\sigma=4$ & $\sigma=5$
& $\sigma=100$ & $\sigma=1000$ \\ \hline
$0$ & $0$ & 1.25578 & 1.59945 & 1.78866 & 1.90808 & 1.98981 & 2.38090 & 2.40242 \\
$0$ & $1/2$ & 1.57080 & 2.02876 & 2.28893 & 2.45564 & 2.57043 & 3.11019 & 3.13845 \\
$1$ & $0$ & 2.40483 & 2.73462 & 2.94960 & 3.09890 & 3.20752 & 3.79360 & 3.82788 \\
$1$ & $1/2$ & 2.74371 & 3.14159 & 3.40561 & 3.59088 & 3.72638 & 4.44850 & 4.48892 \\ \hline
\multicolumn{2}{c|}{$\mu_2/\mu_1$, $N=2$} & 3.66722 & 2.92317 & 2.71939 & 2.63768 & 2.59844 & 2.53875 & 2.53873 \\
\multicolumn{2}{c|}{$\mu_2/\mu_1$, $N=3$} & 3.05095 & 2.39795 & 2.21373 & 2.13831 & 2.10166 & 2.04576 & 2.04575 \\ \hline
\end{tabular}
\caption{{Values of $k_{\nu+\ell,1}$ and
$\mu_2/\mu_1$ for $\ell=0,1$, $N=2,3$, and $\sigma>0$.}}
\label{tab:appendix-positive-spectrum}
\end{table}

For $\sigma=-1$,
$\frac{\mu_2}{\mu_1}\big|_B=0$ when $N\geq2$.
For $\sigma\in(-1,0)$, the ratio is given by
\[
 \frac{\mu_2}{\mu_1}\Big|_B
 =-\frac{k_{\nu+1,1}^2}{\widehat{k}_{\nu,1}^2}.
\]
Table~\ref{tab:appendix-negative-spectrum} gives values for $N=2,3$.
At $\sigma=0$, the ratio is undefined because $\mu_1(B,0)=0$.

\begin{table}[htbp]

\centering
\small
\setlength{\tabcolsep}{5pt}
\begin{tabular}{c|c|c|c|c|c|c}
\hline
 & $\nu$ & $\sigma=-0.1$ & $\sigma=-0.3$ & $\sigma=-0.5$ & $\sigma=-0.7$ & $\sigma=-0.9$ \\ \hline
$k_{\nu+1,1}$ & $0$ & 1.76104 & 1.57883 & 1.35660 & 1.06842 & 0.62721 \\
$k_{\nu+1,1}$ & $1/2$ & 1.98891 & 1.77934 & 1.52553 & 1.19873 & 0.70207 \\
$\widehat{k}_{\nu,1}$ & $0$ & 0.45286 & 0.80454 & 1.06569 & 1.29403 & 1.50599 \\
$\widehat{k}_{\nu,1}$ & $1/2$ & 0.55323 & 0.97767 & 1.28784 & 1.55477 & 1.79867 \\ \hline
\multicolumn{2}{c|}{$\mu_2/\mu_1$, $N=2$} & $-15.12188$ & $-3.85099$ & $-1.62049$ & $-0.68170$ & $-0.17345$ \\
\multicolumn{2}{c|}{$\mu_2/\mu_1$, $N=3$} & $-12.92444$ & $-3.31229$ & $-1.40319$ & $-0.59445$ & $-0.15236$ \\ \hline
\end{tabular}
\caption{{Values of $k_{\nu+1,1}$, $\widehat{k}_{\nu,1}$ and
$\mu_2/\mu_1$ for $N=2,3$ and $-1<\sigma<0$.}}
\label{tab:appendix-negative-spectrum}
\end{table}

\medskip
\noindent\textbf{The interval.}
For comparison with the ball, consider $-u''=\mu u$ on $(0,1)$
with $-u'(0)+\sigma u(0)=0$ and $u'(1)+\sigma u(1)=0$.
Its eigenvalues are simple and satisfy
$\mu_1<\mu_2<\cdots\to+\infty$; see \cite{Bucur,Laugesen}.

For $\mu=k^2>0$ and $\sigma\neq0$, the left boundary condition
gives the eigenfunction
\[
 u(x)=\frac{k}{\sigma}\cos(kx)+\sin(kx),
\]
up to a nonzero constant factor. The right boundary condition then
gives
\[
 (\sigma^2-k^2)\sin k+2\sigma k\cos k=0.
\]
Here $\sin k\neq0$, so equivalently
$\sigma^2+2\sigma k\cot k-k^2=0$. Solving for $\sigma$ gives
the two curves
\[
 \Sigma_\pm(k)=-k\cot k\pm\frac{k}{|\sin k|}.
\]
Positive eigenvalues are obtained from $\sigma=\Sigma_+(k)$ for
$\sigma>0$, and from $\sigma=\Sigma_-(k)$ for $\sigma<0$.
On $(0,\pi)$, their derivatives are
\[
 \Sigma_+'(k)=\frac{(1-\cos k)(k+\sin k)}{\sin^2 k}>0,
 \qquad
 \Sigma_-'(k)=\frac{(1+\cos k)(k-\sin k)}{\sin^2 k}>0.
\]
The first branch increases from $0$ to $+\infty$, while the second
increases from $-2$ to $0$. On each subsequent interval
$(m\pi,(m+1)\pi)$, $m\geq1$, $\Sigma_+$ increases from $0$ to
$+\infty$ and $\Sigma_-$ from $-\infty$ to $0$.
\begin{samepage}
These branches give the positive part of
Figure~\ref{fig:appendix-interval} and the eigenvalue locations
\[
 \begin{aligned}
 (m-1)^2\pi^2&<\mu_m<m^2\pi^2,
 &&\sigma>0,\quad m\geq1,\\
 (m-2)^2\pi^2&<\mu_m<(m-1)^2\pi^2,
 &&-2<\sigma<0,\quad m\geq2,\\
 (m-2)^2\pi^2&<\mu_m<(m-1)^2\pi^2,
 &&\sigma\leq-2,\quad m\geq3.
 \end{aligned}
\]
\end{samepage}

For $\mu=-k^2<0$, with $k>0$, the left boundary condition gives
\[
 u(x)=e^{-kx}-\frac{\sigma+k}{\sigma-k}e^{kx},
\]
up to a nonzero constant factor. Imposing the right boundary
condition gives the two negative-eigenvalue branches
\[
 \Sigma_+=-k\tanh(k/2),\,\,\,
 \Sigma_-=-k\coth(k/2).
\]
As $k$ increases from $0$ to $+\infty$, the first decreases from
$0$ to $-\infty$ and the second from $-2$ to $-\infty$.
The first negative branch therefore exists precisely for $\sigma<0$,
and the second for $\sigma<-2$.
At $\mu=0$, the solution is linear. Substituting it into the
boundary conditions gives $\sigma=0$, with a constant eigenfunction,
or $\sigma=-2$, with eigenfunction $1-2x$. At $\sigma=0$, the full
spectrum is $\mu_m=(m-1)^2\pi^2$, $m\geq1$, with eigenfunctions
$\cos((m-1)\pi x)$.
Thus $\mu_1<0<\mu_2$ for $-2<\sigma<0$, while
$\mu_1<\mu_2=0<\mu_3$ at $\sigma=-2$ and
$\mu_1<\mu_2<0<\mu_3$ for $\sigma<-2$.

In the diagram, the horizontal coordinate is signed:
$k=\sqrt{\mu}\geq0$ for $\mu\geq0$, and
$k=-\sqrt{-\mu}<0$ for $\mu<0$.

\begin{figure}[htbp]

\centering

\resizebox{0.48\textwidth}{!}{%
\begingroup%
\setlength{\unitlength}{1bp}%
\fontsize{10.95}{13.14}\selectfont%
%
\endgroup%
}%

\caption{{The schematic diagram of $\Sigma_+$ (solid lines) and $\Sigma_-$ (dotted lines).}}
\label{fig:appendix-interval}
\end{figure}

\FloatBarrier
Unlike the ball in dimensions $N\geq2$, an interval has at most
two negative Robin eigenvalues. Both tend to $-\infty$ as
$\sigma\to-\infty$. For example, at $\sigma=-3$ on $(0,1)$,
$\mu_1\approx-10.52118$ and $\mu_2\approx-6.63412$.
For $L>0$, \eqref{eq:robin-scaling} gives
$\mu_m((0,L),\sigma)=L^{-2}\mu_m((0,1),L\sigma)$.
Hence the two zero-eigenvalue thresholds on $(0,L)$ are $0$ and $-2/L$.
\\ \\
\textbf{Conflict of interest and data availability.}
\par\medskip
\indent The authors declare no conflicts of interest.
Data sharing is not applicable, as no datasets were generated or analysed.
\\ \\
\textbf{Declaration of AI-assisted writing.}
\par\medskip
\indent During the preparation of this work, the authors used ChatGPT and DeepSeek for language polishing, grammar checking, verification of some computations, and refining the exposition of certain technical arguments. After using these tools, the authors carefully reviewed and edited the content as needed and take full responsibility for the final version of this manuscript.
\\ \\
\textbf{Acknowledgment}
\par\medskip
\indent The authors thank Professor Richard Laugesen for helpful
discussions, comments and corrections on the initial draft, and references.

\bibliographystyle{amsplain}
\makeatletter
\def\@biblabel#1{#1.~}
\makeatother


\providecommand{\bysame}{\leavevmode\hbox to3em{\hrulefill}\thinspace}
\providecommand{\MR}{\relax\ifhmode\unskip\space\fi MR }
\providecommand{\MRhref}[2]{%
  \href{http://www.ams.org/mathscinet-getitem?mr=#1}{#2}
}
\providecommand{\href}[2]{#2}

\end{document}